\documentclass[12pt]{article}
\usepackage{amsmath,amssymb}
\newtheorem{tm}{Theorem}[section]
\newtheorem{lm}[tm]{Lemma}

\newtheorem{re}[tm]{Remark}

\newtheorem{pr}[tm]{Proposition}
 
 \newenvironment{demo}[1]{\par\smallskip\par\begin{trivlist}
\item[]{\bf #1}\ }{\end{trivlist}\par\smallskip\par}
\newcommand{\Proof}{\begin{demo}{{\it Proof.\ }}}
\newcommand{\qed}{\end{demo}}
\newcommand{\toy}{\ \rule[0em]{0.5ex}{1.8ex}}
\newcommand{\QED}{\toy\end{demo}}

\newcommand{\la}{\langle}
\newcommand{\ra}{\rangle}
\newcommand{\nn}{\nonumber}
\newcommand{\III}{{\vert \kern-.10em \vert \kern-.10em \vert}}
\newcommand{\ve}{\varepsilon}

\makeatletter
 
 \@addtoreset{equation}{section}
\makeatother
 
\begin{document}
\setlength{\baselineskip}{15pt} 
%
\bibliographystyle{plain}
\title{
Short time kernel asymptotics for Young SDE 
by means of  Watanabe distribution theory
\footnote
{
Revised on 7 Nov. 2013. \quad
{\bf Mathematics Subject Classification:}~ 60H07, 60F99, 60H10, 60G15.
{\bf Keywords:} fractional Brownian motion, Young integral, stochastic differential equation,
Malliavin calculus, Watanabe distribution, short time asymptotics.
}
}
\author{ Yuzuru INAHAMA
\footnote{
Graduate School of Mathematics,   Nagoya University,
Furocho, Chikusa-ku, Nagoya 464-8602, Japan.
Email:~\tt{inahama@math.nagoya-u.ac.jp} 
}
}
\date{   }
%
%
\maketitle

%

%
%

\begin{center}
{\bf  Abstract}
\end{center}
In this paper we study short time asymptotics of 
a density function of the solution of a stochastic differential equation 
driven by fractional Brownian motion with Hurst parameter $H~ (1/2 <H< 1)$
when the coefficient vector fields satisfy an ellipticity condition at the starting point.
We prove both on-diagonal and off-diagonal asymptotics under mild additional assumptions.
Our main tool is Malliavin calculus, in particular, 
Watanabe's theory of generalized Wiener functionals.
%
%
%

\section{Introduction}

Let $(w_t)_{t \ge 0}$ be the standard $d$-dimensional Brownian motion 
and let $V_i~(0 \le i \le d)$ be smooth vector fields on ${\bf R}^n$ with sufficient regularity.
Consider the following stochastic differential equation (SDE) of Stratonovich-type;
\begin{equation}
dy_t 
=  
\sum_{i=1}^d  V_i (y_t) \circ dw^i_t + V_0 (y_t) dt \qquad  \mbox{ with }  \qquad y_0 =a \in {\bf R}^n. 
\nn
\end{equation}
If the set of vector fields satisfies a hypoellipticity condition,
the solution $y_t =y_t(a)$ has a smooth density $p_t (a,a')$ with respect to Lebesgue measure on ${\bf R}^n$.
From an analytic point of view,  
$p_t (a,a')$ is a fundamental solution of the parabolic equation
$\partial u/\partial t = L u$, where $L = V_0 +(1/2) \sum_{i=1}^d V_i^2$,
and is also called a heat kernel of $L$.

In many fields of mathematics such as    
probability, analysis, mathematical physics, and differential geometry, 
short time asymptotic of $p_t (a,a')$ is a very important problem and 
has been studied extensively.
Although analytic methods are also well-known, we only discuss a probabilistic approach 
via Feynman-Kac formula in this paper.
Malliavin calculus is a very powerful theory 
and was used in many papers on this problem.

Among them, S. Watanabe's result seems to be one of the best.
(See \cite{wa} or Sections 5.8--5.10, \cite{iwbk}.)
His theory of distributional Malliavin calculus is not only very powerful,
but also user-friendly.  
Many heuristic operations are made rigorous in this theory 
and consequently the theory gives us a good view.
Moreover, this theory is quite self-contained in the sense that 
all the argument, from an explicit expression of the heat kernel 
to the final asymptotic result, 
is constructed without much help from other theories.

The theory goes as follows.
First, he constructed a theory of generalized Wiener functionals (i.e., Watanabe distributions)
in Malliavin calculus.
Then, he gave a representation of the heat kernel by using the pullback of Dirac's delta function;
$
p_t (a,a')  = {\mathbb E} [ \delta_{a'} (y_t(a))  ]$,
where the right hand side is the generalized expectation with respect to Wiener measure.
Finally, by establishing an asymptotic expansion theory in the spaces of 
generalized Wiener functionals, 
he obtained a short time expansion of $p_t (a,a')$ under very mild assumptions.
In this method,
an asymptotic expansion is actually obtained before taking the generalized expectation.

In this paper we consider the following problem.
Let $(w^H_t)_{t \ge 0}$ be $d$-dimensional fractional Brownian motion (fBm)
with Hurst parameter $H \in (1/2, 1)$.
Instead of the above SDE, we consider 
\begin{equation}
dy_t 
=  
\sum_{i=1}^d  V_i (y_t) dw^{H,i}_t + V_0 (y_t) dt \qquad  \mbox{ with }  \qquad y_0 =a \in {\bf R}^n. 
\nn
\end{equation}
This is an ordinary differential equation (ODE) in the sense of Young integral 
(see Lyons \cite{lyons}).
In fact, this is actually an ODE with a random driving path, but we call this SDE for simplicity.
Some researchers have studied the solution of the above SDE with Malliavin calculus. 
See \cite{nr, hn, ns, bh, ct} and references therein. 
Under the ellipticity or the hypoellipticity condition,
the solution $y_t=y_t(a)$ has a smooth density $p_t (a,a')$. See \cite{hn,ns,bh}.

In this paper, by using Malliavin calculus and, in particular, Watanabe distribution theory,
we will prove a short asymptotic expansion of this density 
in the elliptic case under mild assumptions.
This kind of asymptotics was already studied in \cite{bh, bo}, but without Malliavin calculus.
In \cite{bh}, they showed on-diagonal short time asymptotics when $V_0 \equiv 0$.
In \cite{bo}, by using Laplace's method, they showed off-diagonal short time asymptotics 
when $V_0 \equiv 0$ and the vector fields $V_i$'s satisfy a rather special condition.
Our results is a generalization of these preceding ones.
Notice that we do not assume the drift term $V_0$ is zero. 
One may think this is just a minor generalization,  
but this makes the asymptotic expansion much more complicated.

The organization of this paper is as follows:
In Section 2, we give settings, assumptions, and precise statements of two main theorems.
In Section 3, we recall basic properties of a Young ODE and 
its Jacobian process for later use.
In Section 4, we review Watanabe's theory of generalized Wiener functionals 
in Malliavin calculus.
In Section 5, we discuss the solution of Young ODE 
driven by fBm with Hurst parameter $H \in (1/2, 1)$ from the viewpoint of Malliavin calculus. 
We also prove uniform non-degeneracy of Malliavin covariance matrix of the solution
under the ellipticity condition.
In Section 6, we prove one of our main theorems, 
namely,  on-diagonal asymptotics of the kernel.
In section 7, we show the shifted solution of the Young SDE admits an asymptotic expansion 
in the sense of Watanabe distribution theory.
In Section 8, we prove the other of our main theorems, 
namely,  off-diagonal asymptotics of the kernel.
In Section 9, 
we prove that, under the ellipticity assumption at the stating point,
our main result (the off-diagonal asymptotics) holds 
when the end point is close enough to the starting point.
We also make sure that Baudoin and Ouyang's result in \cite{bo} is basically included in ours.

%
%
\section{Setting and main results}

\subsection{Setting}
In this subsection, we introduce a stochastic process that will play a main role in this paper.
From now on, dropping the superscript "$H$",
we denote by
$(w_t)_{ t \ge 0} =( w^1_t, \ldots, w^d_t)_{t \ge 0}$
 the $d$-dimensional fractional Brownian motion (fBm) with Hurst parameter
$H ~(1/2 <H<1)$.
It is a unique $d$-dimensional mean-zero Gaussian process with covariance
$$
{\mathbb E} [  w^i_s  w^j_t] = \frac{  \delta_{ij} }{2} ( |s|^{2H} +  |t|^{2H} -  |t-s|^{2H}), 
\qquad 
(s, t \ge 0).
$$
Note that, for any $c>0$, $(w_{ct})_{ t \ge 0}$ and $(c^{H} w_t)_{ t \ge 0}$ have the same law.
This property is called self-similarity or scale invariance.

Let $V_i : {\bf R}^n \to {\bf R}^n ~$ be $C_b^{\infty}$,
that is, $V_i$ is a bounded smooth function with 
bounded derivatives of all order ($0 \le i \le d$).
We consider the following stochastic 
ODE in the sense of Young;
\begin{equation}\label{main.ygODE.eq}
dy_t =  \sum_{i=1}^d  V_i (y_t) dw^i_t + V_0 (y_t) dt \qquad  \mbox{ with }  \qquad y_0 =a \in {\bf R}^n. 
\end{equation}
We will sometimes write $y_t =y_t(a) = y_t(a, w)$ etc. to make explicit the dependence on $a$ and $w$.

\subsection{Assumptions}

In this subsection we introduce assumptions of the main theorems.
First, we assume the ellipticity of the coefficient of (\ref{main.ygODE.eq}) at the starting point $a \in {\bf R}^n$.
\\
\\
{\bf (A1):}~ 
The set of vectors $\{V_1(a), \ldots, V_d(a)\}$ linearly spans ${\bf R}^n$.
\\
\\
It is known that, under Assumption {\bf (A1)}, 
the law of the solution $y_t$ has a density $p_t(a, a')$ with respect to the Lebesgue measure on ${\bf R}^n$
for any $t >0$ (see \cite{bh, ns}).
Hence, for any measurable set $U \subset {\bf R}^n$, 
${\mathbb P} (y_t \in U) = \int_U p_t(a, a') da'$.

Let ${\cal H} = {\cal H}^H$ be the Cameron-Martin space of fBm $(w_t)$.
For $\gamma \in {\cal H}$, we denote by $\phi^0_t  =\phi^0_t (\gamma)$ be the solution
of the following Young ODE;
\begin{equation}\label{main.phi0.eq}
d\phi^0_t 
=  \sum_{i=1}^d  V_i (  \phi^0_t) d\gamma^i_t
 \qquad  \mbox{ with }  \qquad \phi^0_0 =a \in {\bf R}^n. 
\end{equation}
Set,  for $a \neq a'$, 
\[
K_a^{a'} = \{ \gamma \in {\cal H} ~|~  \phi^0_1(\gamma) =a'\}.
\]
If we assume {\bf (A1)} for all $a$, this set $K_a^{a'} $ is not empty.
If $K_a^{a'} $ is not empty, it is a Hilbert submanifold of ${\cal H}$.
From the Schilder-type large deviation theory, it is easy to see that
$\inf\{ \|\gamma\|_{\cal H} ~|~ \gamma \in  K_a^{a'}\} 
= \min\{ \|\gamma\|_{\cal H} ~|~ \gamma \in  K_a^{a'}\}$.
Now we introduce  the following assumption;
\\
\\
{\bf (A2):} $\bar{\gamma} \in K_a^{a'}$ which minimizes ${\cal H}$-norm exists uniquely.
\\
\\
In the sequel, $\bar{\gamma}$ denotes the minimizer in Assumption {\bf (A2)}.
We also assume that 
$\| \,\cdot\, \|^2_{{\cal H}} /2$ is not so degenerate at $\bar\gamma$ in the following sense.
\\
\\
{\bf (A3):} At $\bar{\gamma}$, the Hessian of the functional  
$K_a^{a'} \ni \gamma \mapsto \| \gamma \|^2_{{\cal H}} /2 $ is strictly positive in the form sense.
More precisely, 
if $(- \ve_0, \ve_0) \ni u \mapsto f(u) \in K_a^{a'}$ is a smooth curve in $K_a^{a'}$
such that $f(0) = \bar\gamma$ and $f^{\prime} (0)  \neq 0$, then
$
(d/du)^2 \vert_{u=0} \|   f(u) \|^2_{{\cal H}} /2   > 0.$
\\
\\
Later we will give a more analytical condition {\bf (A3)'}, which is equivalent to {\bf (A3)} 
under {\bf (A2)}.
In \cite{wa}, Watanabe used {\bf (A3)'}.
We will also use {\bf (A3)'} in the proof.
In order to state {\bf (A3)'}, however, we have to introduce a lot of notations.
So, we presented {\bf (A3)} here for ease of presentation.

\subsection{Index sets}

In this  subsection we introduce several index sets for the exponent of the small parameter $\ve >0$, 
which will be used in the asymptotic expansion.
Unlike in the preceding papers,  index sets in this paper are not the set of natural numbers 
and are rather complicated.
Set
\[
\Lambda_1 = \{   n_1 + \frac{n_2}{H}  ~|~  n_1, n_2 \in {\bf N}  \},
\]
where ${\bf N}= \{0,1,2,\ldots\}$.
We denote by $0=\kappa_0 <\kappa_1 < \kappa_2 < \cdots$ all the elements of $\Lambda_1$ in increasing order.
Several  smallest elements are explicitly given as follows;
\[
\kappa_1 = 1, \quad   \kappa_2 =  \frac{1}{H},  \quad   \kappa_3 = 2, \quad \kappa_4 = 1+\frac{1}{H},  
\quad  \kappa_5 = 3 \wedge \frac{2}{H}, \ldots
\]
As usual,  using the scale invariance (i.e., self-similarity) of fBm,
we will consider the scaled version of (\ref{main.ygODE.eq}).
(See the scaled Young ODE (\ref{sc2_ygSDE.eq}) below).
From its explicit form, one can easily see why $\Lambda_1$ appears. 

We also set 
\[
\Lambda_2 = \{ \kappa -1 ~|~ \kappa \in \Lambda_1 \setminus \{0\} \}
=
\Bigl\{ 
0,  \, \frac{1}{H} -1,  \, 1,   \, \frac{1}{H}, 
  \, \bigl(  3 \wedge \frac{2}{H}  \bigr)-1,  \ldots
 \Bigr\}
\]
and 
\[
\Lambda'_2 = \{ \kappa -2 ~|~ \kappa \in \Lambda_1 \setminus \{0, 1, 1/H\} \}
=
\Bigl\{ 
0,   \, \frac{1}{H} -1, 
  \, \bigl(  3 \wedge \frac{2}{H}  \bigr)-2,  \ldots
 \Bigr\}.
\]
Next we set 
\[
\Lambda_3 = 
\{  a_1 +a_2 + \cdots + a_m ~|~  \mbox{$m \in {\bf N}_+$ and $a_1 ,\ldots, a_m \in \Lambda_2$} \}.
\]
In the sequel, $\{ 0=\nu_0 <\nu_1<\nu_2 <\cdots \}$ stands for all the elements of 
$\Lambda_3$ in increasing order.
Similarly, 
\[
\Lambda'_3 = 
\{  a_1 +a_2 + \cdots + a_m ~|~  \mbox{$m \in {\bf N}_+$ and $a_1 ,\ldots, a_m \in \Lambda'_2$} \}.
\]
In the sequel, $\{ 0=\rho_0 <\rho_1<\rho_2 <\cdots \}$ stands for all the elements of 
$\Lambda'_3$ in increasing order.
Finally,  
$$\Lambda_4 =\Lambda_3  +\Lambda'_3 =\{ \nu + \rho ~|~ \nu   \in  \Lambda_3  , \rho \in  \Lambda'_3\}.$$
We denote by $\{ 0 = \lambda_0 < \lambda_1 < \lambda_2 < \cdots \}$ all the elements of $\Lambda_4$  in increasing order.
%

\subsection{Statement of the main results}

In this subsection we state two main results of ours,
which are basically analogous to the corresponding ones in Watanabe \cite{wa}.
However, there are some differences.
First,  the exponents of $t$ are not (a constant multiple of) natural numbers.
Second, cancellation of "odd terms" as in p. 20 and p. 34, \cite{wa}
does not happen in general in our case. 
(If the drift term in Young ODE (\ref{main.ygODE.eq}) is zero, 
then this kind of cancellation takes place as in \cite{bh, bo}).

The following is a short time asymptotic expansion of the diagonal of the kernel function.
This is much easier than the off-diagonal case.


\begin{tm}\label{thm.MAIN.on}
Assume {\bf (A1)}.
Then, the diagonal of the kernel $p(t, a, a)$ admits the following 
asymptotics as $t \searrow 0$;
\[
p(t, a, a) \sim
 \frac{1}{t^{n H}} \bigl( c_0 + c_{\nu_1} t^{\nu_1 H} +  c_{\nu_2} t^{\nu_2 H} +  \cdots \bigr)
\]
for certain real constants $c_0, c_{\nu_1}, c_{\nu_2}, \ldots$.
Here, $\{ 0=\nu_0 <\nu_1<\nu_2 <\cdots \}$ are all the elements of 
$\Lambda_3$ in increasing order.
\end{tm}


We also have off-diagonal short time  asymptotics of the kernel function.

\begin{tm} \label{thm.MAIN.off}
Assume $a \neq a'$ and {\bf (A1)}--{\bf (A3)}. 
Then, we have the following asymptotic expansion as $t \searrow 0$;
    \[
p(t , a, a')  \sim
    \exp \Bigl(  -\frac{ \|  \bar\gamma \|^2_{{\cal H}}}{2  t^{2H} }  
    + \frac { \beta }  { t^{2 H - 1} }   \Bigr)
  \frac{1}{t^{n H}} 
     \bigl\{
      \alpha_{  \lambda_0}  +  \alpha_{\lambda_1}  t^{ \lambda_1 H} + \alpha_{\lambda_2}  t^{\lambda_2  H  }  +\cdots
            \bigr\}
\]
 for certain real constants $\beta, \alpha_{\lambda_j}  ~ (j =0,1,2,\ldots)$. 
 Here, $\{ 0=\lambda_0 <\lambda_1<\lambda_2 <\cdots \}$ are all the elements of 
$\Lambda_4$ in increasing order.
 \end{tm}

\begin{re}\label{re.const}
{\bf (i)}~
Consider the following simplest case; $n=d =1$ and $y_t =a+ w_t + b t$ with $b \in {\bf R}$.
Then, for each $t>0$, this induces a Gaussian measure with mean $a +b t$ and variance $t^{2H}$.
Hence, the kernel is given by
\begin{align}
p(t,a, a') &= \frac{1}{\sqrt{2 \pi} t^H }  \exp \Bigl( -\frac{ (a +bt -a')^2 }{ 2 t^{2H}} \Bigr)
\nn\\
&=
\frac{1}{\sqrt{2 \pi} t^H } e^{ - (a - a')^2 /(2t^{2H}) } e^{ - b (a - a') /t^{2H -1} }
e^{- b^2 t^{2-2H} /2}
\nn\\
&=
e^{ - (a - a')^2 /(2t^{2H})  - b (a - a') /t^{2H -1} }
\frac{1}{\sqrt{2 \pi} t^H } 
\Bigl(
1 - \frac{b^2}{2} t^{2( H^{-1}-1)H } +\frac{ b^4}{2^2 2!}   t^{4( H^{-1}-1)H } - \cdots
\Bigr).
\nn
\end{align}
This example may illustrates that the asymptotics in Theorem \ref{thm.MAIN.off} are not so strange.
\\
\noindent
{\bf (ii)}~
Some of the constants in Theorems \ref{thm.MAIN.on} and \ref{thm.MAIN.off}
can be obtained explicitly.
For example, in Theorems \ref{thm.MAIN.on},
$c_0 = [(2\pi)^{n/2} \det(\sigma(a) \sigma(a)^*) ]^{-1}$ and
$$
c_{\nu_1}= c_{(1/H)-1} = \sum_{j=1}^n 
\partial_j \delta_0 \bigl(  V_1 (a) w^1_1 +\cdots +  V_d (a) w^d_1 \bigr) \cdot V_0 (a)^j =0.
$$
Here, $\sigma(a) \sigma(a)^*$ is the covariance matrix of the $n$-dimensional
Gaussian random variable $\sum_{j=1}^d V_j (a) w^j_1$.
In Theorems \ref{thm.MAIN.off}, $\beta =  \la \bar\nu ,  \phi_1^{1/H} \ra$.
The notations in this remark will be given later.
\end{re}


\subsection{Outline of proof of off-diagonal asymptotics}
In this subsection we outline the proof of Theorem \ref{thm.MAIN.off} in a heuristic way
so that the reader would not get lost in technical details.
The argument in this subsection is not rigorous.
For $\ve \in (0,1]$ and $\bar{\gamma}$ as in {\bf (A2)}, consider the following SDE;
\[
d\tilde{y}^{\ve}_t 
=  \sum_{i=1}^d  V_i (\tilde{y}^{\ve}_t) ( \ve dw^i_t + d  \bar{\gamma}_t )+ 
 V_0 (\tilde{y}^{\ve}_t) \ve^{1/H}  dt \qquad  \mbox{ with }  \qquad \tilde{y}^{\ve}_0 =a 
\]
(We denote by $y^{\ve}$ 
the solution of the above ODE with $\bar{\gamma} =0$.)

From the scaling property of fBm and  a routine argument in Watanabe's theory, 
\begin{align}
p(\ve^{1/H}, a, a') 
=
{\mathbb E} \bigl[  
\delta_{a'} ( y_{\ve^{1/H} })   
 \bigr]
 =
{\mathbb E} \bigl[  
\delta_{a'} ( y_1^{\ve})   
 \bigr]
 = 
   {\mathbb E} \bigl[  
\delta_{a'} ( y_1^{\ve} )   \chi_{\eta} (\ve, w)
 \bigr]
   + \mbox{(a small term).}
   \nn
\end{align}
Here, $ \chi_{\eta} (\ve, w)$ is a ${\bf D}_{\infty}$-functional  
which looks like the indicator of a small ball of a certain radius $\eta >0$ centered at $\bar{\gamma}$.
By Schilder-type large deviations, the second term above is negligible. 
By Cameron-Martin theorem, 
the fisrt term is equal to
\[
\exp \bigl(  -\frac{ \|  \bar\gamma \|^2_{{\cal H}}}{2\ve^2} \bigr)
{\mathbb E} \bigl[  
   \exp \bigl(  - \frac{1}{\ve} \la \bar\gamma, w \ra \bigr)
\delta_{a'} ( \tilde{y}_1^{\ve} )   \chi_{\eta} (\ve, w + \frac{\bar\gamma}{\ve})
 \bigr].
 \]
Here, $\chi_{\eta} (\ve, w + \bar\gamma /\ve)$ 
 does not contribute to the asymptotic expansion since it is of the form $1 + O(\ve^N)$
for any large $N \in {\bf N}$.
So, it is sufficient to consider the two factors; $\delta_{a'} ( \tilde{y}_1^{\ve} )$
and $\exp ( -\la \bar\gamma, w \ra /\ve)$.

We will prove in Section \ref{sec.taylor} that
$\tilde{y}_1^{\ve}$ admits the following expansion
for certain $\phi^{\kappa_j}$'s
both in ${\bf D}_{\infty} ({\bf R}^n)$-sense and the deterministic sense.
\[
\tilde{y}^{\ve}_1 \sim \phi^0_1 + \ve^{\kappa_1} \phi^{\kappa_1}_1 + \ve^{\kappa_2} \phi^{\kappa_2}_1 
+ \cdots    \qquad \mbox{as $\ve \searrow 0$,}
\qquad
(\kappa_i \in \Lambda_1={\bf N} +\frac{1}{H} {\bf N})
\]
From the SDE for $\tilde{y}^{\ve}$, one can easily see that 
the index set for this Taylor expansion of It\^o map should be $ \Lambda_1$.
Set $R^{2,\ve} 
= \tilde{y}^{\ve} -( \phi^0 + \ve \phi^{1} + \ve^{1/H} \phi^{1/H} )$.
In fact, $\phi^0, ~ \phi^{1/H}$ do not depend on $w$.  
Then, we see from $\phi^0_1 =a'$ that
\[
\delta_{a'} ( \tilde{y}_1^{\ve} ) 
= \delta_{0} \bigl( \ve \cdot \frac{\tilde{y}_1^{\ve} -a'}{\ve} \bigr) 
=
\ve^{-n} \delta_{0} (  \phi^{1}_1 + \ve^{(1/H) -1} \phi^{1/H}_1 + \ve^{-1}  R^{2,\ve}_1 ). 
\]
Since 
$ (\tilde{y}_1^{\ve} -a')/\ve =  \phi^{1}_1 + \ve^{(1/H) -1} \phi^{1/H}_1 + \ve^{-1}  R^{2,\ve}_1$
is uniformly non-degenerate in $\ve$ in the sense of Malliavin under {\bf (A1)}
and indexed by $\Lambda_2$,
its composition with the Dirac measure $\delta_0$ is well-defined and 
admits a Taylor-like expansion with the index set $\Lambda_3$.

Next we consider the other factor.
We will show that there exists $\bar{\nu} \in {\bf R}^n$
such that $\la \bar\gamma, w\ra = \la \bar{\nu}, \phi^1_1\ra$,
where the right hand side is the inner product of ${\bf R}^n$.
Under the condition that
$ \phi^{1} + \ve^{(1/H) -1} \phi^{1/H} + \ve^{-1}  R^{2,\ve} =0$,
we have
\[
 \exp \bigl(  - \frac{1}{\ve} \la \bar\gamma, w \ra \bigr)
 =
 \exp \bigl(  \frac{\la \bar\nu, \phi^{1/H}_1 \ra }{ \ve^{2-1/H} } \bigr)
 \cdot
 \exp \bigl(  \frac{\la \bar\nu, R^{2,\ve} \ra }{ \ve^{2} } \bigr).
   \]
It is obvious that the index set for $R^{2,\ve}/ \ve^{2}$ is $\Lambda'_2$,
which implies 
that the index set for $\exp (\la \bar\nu, R^{2,\ve}/ \ve^{2} \ra)$ is $\Lambda'_3$.
From this heuristic explanation, we see that 
$p(\ve^{1/H}, a, a') $ admits an  asymptotic expansion 
and why $\Lambda_4 =\Lambda_3+ \Lambda'_3$ appears as the index set of the asymptotics.
By setting $\ve =t^H$, we have the desired  short time expansion.

When we try to make the above argument rigorous, 
the most difficult part is to prove integrability of various 
Wiener functionals of exponential-type.
This is highly non-trivial and we will prove a few lemmas 
for that purpose in Subsection \ref{subsec.int}.
Assumption {\bf (A3)} is actually a sufficient condition 
for those lemmas to hold.


\section{Basic properties of Young ODE and 
$L^q$-integrability of Jacobian process}
\label{sec.yint}

In this section we recall the basic properties of a Young ODE and 
its Jacobian process (i.e., derivative process).
There is no new result in this section.
These facts are scattered across many literatures and it is not so easy to find a suitable one.
(In this sense, Lejay \cite{le} may be useful.)
Here, we summarize some results, in particular, $L^q$-integrability 
of
the Jacobian process driven by fBm with Hurst parameter $H >1/2$ for later use.
(Z\"ahle \cite{za} generalized 
Young integral and ODE by using fractional calculus, 
but we do not use it in this paper.)


We always assume that $1/2 <\alpha \le 1$ and the time interval is $[0,1]$. 
Let $C^{\alpha- hld} ([0,1] ; {\bf R}^d)$ 
be the spaces of ${\bf R}^d$-valued $\alpha$-H\"older continuous paths.
The Banach norms are defined by
\begin{align}
\|x\|_{ \alpha - hld } &= |x_0| + \sup_{0 \le s <t \le 1} \frac{|x_t -x_s|}{(t-s)^{\alpha}},
\nn
\end{align}
The closed subspaces of paths that starts at the origin is denoted by 
$C_0^{\alpha- hld} ([0,1] ; {\bf R}^d)$.

Let $\sigma : {\bf R}^n \to {\rm Mat} (n, d)$ 
and 
$b: {\bf R}^n \to {\bf R}^n$
be sufficiently regular.
Consider the following ODE in the Young sense;
\begin{equation}\label{ygODE.eq}
dy_t = \sigma(y_t) dx_t 
+ b(y_t) dt 
\qquad  \mbox{ with }  \qquad y_0 =a. 
\end{equation}
Here, $x \in C_0^{\alpha- hld} ([0,1] ; {\bf R}^d)$  and $a \in {\bf R}^n$ is the initial value.
Let $V_i : {\bf R}^n \to {\bf R}^n$ be the $i$th column vector of $\sigma$ ($1 \le i \le d$)
and set $V_0 =b$.
%
%
Then, ODE (\ref{ygODE.eq}) can be rewritten equivalently as follows;
\begin{equation}\label{ygODE2.eq}
dy_t =  \sum_{i=1}^d  V_i (y_t) dx^i_t 
+ V_0 (y_t) dt 
\qquad  \mbox{ with }  \qquad y_0 =a. 
\end{equation}
Some researchers prefer this style.
In this paper we will use both (\ref{ygODE.eq}) and (\ref{ygODE2.eq}).

Assume $\sigma$  and $b$ are $C^2_b$, 
that is, $ \max_{0 \le i \le 2} (\| \nabla^i \sigma \|_{\infty} +   \| \nabla^i b \|_{\infty}) <\infty$,
where $\| \, \cdot \,\|_{\infty} $ stands for the sup-norm. 
Then the above  ODE has a unique solution for any given $x$ and $a$
in $\alpha$-H\"older setting.
Moreover, the map
\begin{equation}\label{locL.eq}
  C_0^{\alpha-hld} ([0,1] ; {\bf R}^d) \times {\bf R}^n  
   \ni (x,a)  \mapsto y \in C^{\alpha-hld} ([0,1] ; {\bf R}^n)
\end{equation}
is locally Lipschitz continuous (i.e., Lipschitz continuous on any bounded set).
We will sometimes write $y =I(x,\lambda)$, where $\lambda_t =t$.
(In this paper $a$ is fixed.)


Now we discuss the Jacobian process (i.e., the derivative process) $J$ of the ODE (\ref{ygODE.eq}), or equivalently (\ref{ygODE2.eq}).
$J_t$ is a (formal) derivative of the solution flow $a \mapsto y_t=y_t(a)$ of the Young ODE (\ref{ygODE.eq}).

For $v \in {\bf R}^n$, we denote the directional derivative along $v$
by $\nabla_v \sigma (y) = \nabla \sigma (y) \la  v, \,\cdot\, \ra$, etc.
So, $ \nabla \sigma$ takes its values in 
$
L^{(2)} ( {\bf R}^n,  {\bf R}^d; {\bf R}^n) =   ( {\bf R}^{n} )^* \otimes  ( {\bf R}^{d} )^* \otimes  {\bf R}^n,
$
which is equipped with the usual Hilbert-Schmidt norm.
Notations such as $\nabla^i  V_j$, 
$ \nabla^2  \sigma = \nabla\nabla  \sigma$, $\nabla^2  b$, etc. should be understood in a similar way.

The Jacobian process $J$ takes its values in ${\rm Mat }(n,n) = L({\bf R}^n,  {\bf R}^n )$ and satisfies  
\begin{equation}\label{jODE1.eq}
dJ_t  =  \nabla \sigma(y_t) \la  J_t, dx_t \ra
+ \nabla b(y_t) \la  J_t \ra dt 
\qquad  \mbox{ with }  \qquad   J_0={\rm Id}_n.
\end{equation}
More precisely,  by setting 
$M_t = \int_0^t \{ \nabla \sigma(y_s) \la  \,\cdot\, , dx_s \ra 
+ \nabla b(y_s) \la  \,\cdot\,  \ra  ds \}$,
we may rewrite this equation  as follows;
\begin{equation}\label{jODE2.eq}
dJ_t  =  dM_t  \cdot J_t  \qquad  \mbox{ with }  \qquad   J_0={\rm Id}_n.
\end{equation}
The dot on the right hand denotes the matrix multiplication.
When we need to specify the driving path, 
we will write $J(x, \lambda)$, where $\lambda_t =t$. 
The equivalent equation for $J$ that corresponds to (\ref{ygODE2.eq}) is as follows;
\begin{equation}\label{jODE3.eq}
dJ_t  = \sum_{i=1}^d  \nabla V_i (y_t)  \la  J_t \ra  dx^i_t 
+ \nabla V_0 (y_t)    \la  J_t \ra  dt 
\qquad  \mbox{ with }  \qquad  J_0={\rm Id}_n.
\end{equation}

Assume for safety that $\sigma$ and $b$ are $C_b^3$.
It is known that the system of Young ODEs  (\ref{ygODE.eq}) and (\ref{jODE1.eq}) has 
a unique solution $(y, J)$ for given 
$x \in C_0^{\alpha- hld} ([0,1] ; {\bf R}^d)$  and $a \in {\bf R}^n$
in $\alpha$-H\"older setting
and local Lipschitz continuity of $(x,a) \mapsto (y,J)$ also holds in this case.

Now let us consider the moment estimate for H\"older norms of 
 $J$ and $J^{-1}$, when the driving path $x$ is 
the $d$-dimensional fBm $w=(w_t)_{0 \le t \le 1}$ with Hurst parameter $H \in (1/2, 1)$.
Take any $\alpha \in ( 1/2, H)$. 
Then, almost surely, $\|  w\|_{\alpha -hld} <\infty$.
(By the way, $\|  w\|_{1/H -var} =\infty$, a.s. See \cite{ckr, pr}.
Hence, $\|  w\|_{H -hld} =\infty$, a.s.)

The differential equations are given as follows; 
\begin{align}
dy_t  = \sigma(y_t) dw_t + b(y_t) dt  \,\, \mbox{ with }  \, y_0 =a 
\quad\mbox{ and }\quad
dJ_t  =  dM_t  \cdot J_t   \,\, \mbox{ with }  \, J_0={\rm Id}_n,
\label{Jdrfbm.eq}
\end{align}
where $M_t = 
\int_0^t \{ \nabla \sigma(y_s) \la  \,\cdot\, , dw_s   \ra + \nabla b(y_s) \la  \,\cdot\,   \ra ds\}$.
For simplicity we call them SDEs, though they are just deterministic Young ODEs 
driven by a random input $w$ (and $\lambda$).

\begin{pr}\label{pr.momJ}
Let $1/2 <\alpha <H$ and assume that the coefficients $\sigma$ and $b$ are $C^3_b$. 
Let $J$ be as in (\ref{Jdrfbm.eq}) above.
Then,  $\|  J\|_{\alpha -hld} $ and  $\|  J^{-1}\|_{\alpha -hld}$ have moments of all order, i.e.,  
$\|  J\|_{\alpha -hld}, \|  J^{-1}\|_{\alpha -hld} \in \cap_{1 \le q <\infty} L^q$.
\end{pr}

\Proof
This is already known. Here, we give a sketch of proof only. 

Since (\ref{jODE1.eq}) is linear, the solution can be written explicitly as follows. 
\begin{equation}\label{jexp1.eq}
J_t   = \Bigl( {\rm Id}_n + \sum_{k=1}^{\infty} M^{[k]}_{s,t}   \Bigr) J_s
\qquad   \quad  
(0 \le s \le t  \le 1),
\end{equation}
where 
\begin{equation}\label{Mexp1.eq}
 M^{[k]}_{s,t}  =  \int_{ s \le t_1 \le \cdots \le t_k \le t}  dM_{t_k} \cdots  dM_{t_2} dM_{t_1}.
 \end{equation}

We can apply the same argument as in the proof of 
Lyons' extension theorem  (p.35, \cite{lq}) to  obtain
\begin{align}
  \| J\|_{\alpha -hld}   \le  1 + c' (1 +  \| w\|^{1/\alpha}_{\alpha -hld} )
   \exp( c \| w \|^{1/\alpha}_{\alpha -hld}  ).
  \label{Jvar2.ineq}
\end{align}
Here, positive constants
$c, c'$  depend only on $\alpha, \sigma, b$.
Since $1/\alpha <2$, we can apply Fernique's square exponential integrability theorem
for Gaussian measures.

$J^{-1}$ has a series expansion similar to (\ref{jexp1.eq})--(\ref{Mexp1.eq})
and can be dealt with in the same way.

It is also possible to prove Proposition \ref{pr.momJ} by using 
Hu and Nualart's result on integrability of
$\sup_{0 \le t \le 1}|J_t|$ in \cite{hn}
plus 
a cutoff argument.
\QED

\begin{re}\label{re.supint}
This kind integrability problem for Jacobian process becomes very difficult when $H<1/2$.
Cass, Litterer, and Lyons \cite{cll} recently proved it in rough path setting 
for Gaussian rough path including 
fractional Brownian rough path with $1/4 <H \le 1/2$.
\end{re}


\section{Preliminaries from Watanabe's asymptotic theory of generalized Wiener functionals}

We recall Watanabe's theory of generalized Wiener functionals 
in Malliavin calculus.
Most of the contents and the notations
in this section are borrowed from \cite{wa} or Sections 5.8--5.10, Ikeda and Watanabe \cite{iwbk} 
with trivial modifications.
Shigekawa \cite{shbk} and Nualart \cite{nbk} are also good textbooks of Malliavin calculus
and we will sometimes refer to them.
There is no new result in this section.

Let $(W, {\cal H}, \mu)$ be an abstract Wiener space.  
(The results in \cite{wa}
or Sections 5.8--5.10, \cite{iwbk} also holds on any abstract Wiener space.)
The following are of particular importance in this paper:
\\
\\
{\bf (a)}~ Basics of Sobolev spaces ${\bf D}_{q,r} ({\cal K})$ of ${\cal K}$-valued 
(generalized) Wiener functionals, 
where $q \in (1, \infty)$, $r \in {\bf R}$, and ${\cal K}$ is a real separable Hilbert space.
As usual, we will use the spaces ${\bf D}_{\infty} ({\cal K})$, $\tilde{{\bf D}}_{\infty} ({\cal K})$ of test functions 
and  the spaces ${\bf D}_{-\infty} ({\cal K})$, $\tilde{{\bf D}}_{-\infty} ({\cal K})$ of 
generalized Wiener functionals (i.e., Watanabe distributions) as in \cite{iwbk}.
\\
{\bf (b)}~ Meyer's equivalence of Sobolev norms. 
(Theorem 8.4, \cite{iwbk}. A stronger version can be found in Theorem 4.6, \cite{shbk})
\\
{\bf (c)}~Pullback $T \circ F$ of tempered Schwartz distribution $T \in {\cal S}^{\prime}({\bf R}^n)$
on ${\bf R}^n$
by a non-degenerate Wiener functional $F \in {\bf D}_{\infty} ({\bf R}^n)$. (see Sections 5.9, \cite{iwbk}.)
\\
{\bf (d)}~A generalized version of integration by parts formula in the sense 
of Malliavin calculus
 for Watanabe distribution. (p. 7, \cite{wa} or p. 377, \cite{iwbk})
\\
\\

Now we consider a family of Wiener functionals indexed by a small parameter $\ve \in (0,1]$.
When the index set of asymptotics is ${\bf N}$, it is explained in  Sections 5.9, \cite{iwbk}.
This is just a slight generalization of it.

Consider a family of ${\cal K}$-valued  Wiener functionals $\{ F(\ve, w)\}_{0<\ve \le 1}$
and assume $F(\ve, \,\cdot\,) \in {\bf D}_{\infty} ({\cal K})$ for each $\ve$.
We say $F(\ve, \,\cdot\,) =O (\ve^{\kappa})$ in ${\bf D}_{q,k}({\cal K})$, $\kappa \in {\bf R}$, as $\ve \searrow 0$,
if  $\| F(\ve, \,\cdot\,)  \|_{q, k}= O (\ve^{\kappa})$.
We say 
$F(\ve, \,\cdot\,) =O (\ve^{\kappa})$ in ${\bf D}_{\infty}({\cal K})$ as $\ve \searrow 0$,
if 
$F(\ve, \,\cdot\,) =O (\ve^{\kappa})$ in ${\bf D}_{p,k}({\cal K})$ for any $1<q<\infty$ and  $k \in {\bf N}$.

Let $0=\kappa_0 <\kappa_1 <\kappa_2<\cdots \nearrow \infty$ and
$f_0, f_{\kappa_1}, f_{\kappa_2}, \ldots \in {\bf D}_{\infty} ({\cal K})$. 
We write
\[
F(\ve, \,\cdot\,) \sim  f_0 + \ve^{\kappa_1} f_{\kappa_1} +\ve^{\kappa_2} f_{\kappa_2} + \cdots   
\qquad
\mbox{in ${\bf D}_{\infty} ({\cal K})$ as $\ve \searrow 0$,}
\]
if, for any $m \in {\bf N}$, it holds that
\[
F(\ve, \,\cdot\,) - ( f_0 + \ve^{\kappa_1} f_{\kappa_1} + \cdots +\ve^{\kappa_m} f_{\kappa_m})
 = O (\ve^{\kappa_{m+1}  })
\qquad
\mbox{in ${\bf D}_{\infty} ({\cal K})$ as $\ve \searrow 0$.}
\]

In a similar way, we can define asymptotic expansions in 
${\bf D}_{- \infty} ({\cal K})$, $\tilde{{\bf D}}_{ \infty} ({\cal K})$, 
$\tilde{{\bf D}}_{- \infty} ({\cal K})$ for a general index set, too,
but we omit them.

We recall basic facts for such asymptotic expansions in the Sobolev spaces.
Let  $0=\kappa_0 < \kappa_1 < \kappa_2 < \cdots \nearrow \infty$ be as above.
In Proposition \ref{pr.Wasym.prod} below, $0=\nu_0 < \nu_1 < \nu_2 < \cdots \nearrow \infty$ are 
all the elements of $\{  \kappa_i +\kappa_j ~|~ i, j \in {\bf N} \}$ in increasing order.
The fundamental case $\kappa_j =j$ is treated
in Proposition 9.3, Section 5.9, \cite{iwbk}.
The following is a straight forward modification of it.

\begin{pr}\label{pr.Wasym.prod}
{\bf (i)}~ Suppose that  $F (\ve, \,\cdot\,) \in {\bf D}_{\infty} ({\cal K})$ 
admits an expansion such as
\[
F(\ve, \,\cdot\,) \sim  f_0 + \ve^{\kappa_1} f_{\kappa_1} +\ve^{\kappa_2} f_{\kappa_2} + \cdots   
\qquad
\mbox{in ${\bf D}_{\infty} ({\cal K})$ as $\ve \searrow 0$,}
\]
with $f_{\kappa_j} \in {\bf D}_{\infty} ({\cal K})$ for all $j \in {\bf N}$.  
Suppose also that $G (\ve, \,\cdot\,) \in {\bf D}_{\infty}  ~\mbox{(or $\tilde{{\bf D}}_{\infty}$)}$  admits an expansion such as 
\[
G(\ve, \,\cdot\,) \sim  g_0 + \ve^{\kappa_1} g_{\kappa_1} +\ve^{\kappa_2} g_{\kappa_2} + \cdots   
\qquad
\mbox{in ${\bf D}_{\infty}$   ~ (or resp.  $\tilde{{\bf D}}_{\infty}$) as $\ve \searrow 0$,}
\]
with $g_{\kappa_j} \in {\bf D}_{\infty}  ~ \mbox{(or resp. $\tilde{{\bf D}}_{\infty} $)} $ for $j \in {\bf N}$.  
Then,  $H(\ve, w)= F (\ve, w)G (\ve, w)$ satisfies that
\[
H(\ve, \,\cdot\,) \sim  h_0 + \ve^{\nu_1} h_{\nu_1} +\ve^{\nu_2} h_{\nu_2} + \cdots   
\qquad
\mbox{in ${\bf D}_{\infty}({\cal K})$   ~ (or resp.  $\tilde{{\bf D}}_{\infty}({\cal K})$) 
as $\ve \searrow 0$,}
\]
where $h_{\nu_n} \in {\bf D}_{\infty}({\cal K})~ \mbox{(or resp. $\tilde{{\bf D}}_{\infty} ({\cal K})$)} $ are given by 
the following formal multiplication;
\[
h_{\nu_n} = \sum_{(i,j); \kappa_i + \kappa_j=   \nu_n}  g_{\kappa_i} f_{\kappa_j}.
\]
\\
\\
{\bf (ii)}~
Suppose that $G (\ve, \,\cdot\,) \in {\bf D}_{\infty}  ~\mbox{(or $\tilde{{\bf D}}_{\infty}$)}$  
admits an expansion such as 
\[
G(\ve, \,\cdot\,) \sim  g_0 + \ve^{\kappa_1} g_{\kappa_1} +\ve^{\kappa_2} g_{\kappa_2} + \cdots   
\qquad
\mbox{in ${\bf D}_{\infty}$   ~ (or resp.  $\tilde{{\bf D}}_{\infty}$) as $\ve \searrow 0$,}
\]
with $g_{\kappa_j} \in {\bf D}_{\infty}  ~ \mbox{(or resp. $\tilde{{\bf D}}_{\infty} $)} $ for all $j \in {\bf N}$.  
Suppose also that $\Phi (\ve , \,\cdot\,) \in \tilde{{\bf D}}_{-\infty} ({\cal K})$ 
admits an expansion such as
\[
\Phi(\ve, \,\cdot\,) \sim  \phi_0 + \ve^{\kappa_1} \phi_{\kappa_1} +\ve^{\kappa_2} \phi_{\kappa_2} + \cdots   
\qquad
\mbox{in $\tilde{{\bf D}}_{-\infty} ({\cal K})$ as $\ve \searrow 0$,}
\]
with $\phi_{\kappa_j} \in \tilde{{\bf D}}_{- \infty} ({\cal K})$ for all $j \in {\bf N})$.  
Then, $\Psi (\ve, w)= G (\ve, w)\Phi(\ve, w)$ satisfies that
\begin{equation}\label{eq.psi.mult1}
\Psi(\ve, \,\cdot\,) \sim  \psi_0 + \ve^{\nu_1} \psi_{\nu_1} +\ve^{\nu_2} \psi_{\nu_2} + \cdots   
\qquad
\mbox{in $\tilde{{\bf D}}_{-\infty} ({\cal K}) ~\mbox{(or resp. ${\bf D}_{-\infty} ({\cal K}))$}$ as $\ve \searrow 0$,}
\end{equation}
where $\psi_{\nu_n} \in  \tilde{{\bf D}}_{-\infty}({\cal K})~ 
\mbox{(or resp. ${\bf D}_{-\infty} ({\cal K})$)} $ are given by 
the following formal multiplication;
\begin{equation}\label{eq.psi.mult2}
\psi_{\nu_n} = \sum_{(i,j); \kappa_i + \kappa_j=   \kappa_n}  g_{\kappa_i} \phi_{\kappa_j}.
\end{equation}
\\
\\
{\bf (iii)}~
Suppose that $G (\ve, \,\cdot\,) \in {\bf D}_{\infty}$  admits an expansion such as 
\[
G(\ve, \,\cdot\,) \sim  g_0 + \ve^{\kappa_1} g_{\kappa_1} +\ve^{\kappa_2} g_{\kappa_2} + \cdots   
\qquad
\mbox{in ${\bf D}_{\infty}$ as $\ve \searrow 0$,}
\]
with $g_{\kappa_j} \in {\bf D}_{\infty}$ for all $j \in {\bf N}$.  
Suppose also that $\Phi (\ve , \,\cdot\,) \in {\bf D}_{-\infty} ({\cal K})$ 
admits an expansion such as
\[
\Phi(\ve, \,\cdot\,) \sim  \phi_0 + \ve^{\kappa_1} \phi_{\kappa_1} +\ve^{\kappa_2} \phi_{\kappa_2} + \cdots   
\qquad
\mbox{in ${\bf D}_{-\infty} ({\cal K})$ as $\ve \searrow 0$,}
\]
with $\phi_{\kappa_j} \in {\bf D}_{- \infty} ({\cal K})$ for all $j \in {\bf N}$.  
%
Then, (\ref{eq.psi.mult1}) and  (\ref{eq.psi.mult2}) hold in ${\bf D}_{- \infty} ({\cal K})$. 
\end{pr}

\begin{re}\label{re.soroe.ind}
In {\bf (i)} of the above Proposition, 
the index sets $\{\kappa_j \}_{j=0,1,2, \ldots}$ for the asymptotic expansions for 
$F(\ve, \,\cdot\,) $ and $G(\ve, \,\cdot\,)$ are the same.
However, these index sets for $F$ and $G$ may differ, 
because  the union of the two index sets can be regarded as a new index set.
Similar remarks hold for  {\bf (ii)}  and {\bf (iii)}, too.
\end{re}

Next we consider asymptotic expansions for  the pullback.
Let $F(\ve , \,\cdot\,) \in {\bf D}_{\infty} ({\bf R}^n)$ for $0<\ve \le 1$.
We say $F$ is uniformly non-degenerate in the sense of Malliavin if 
\[
\sup_{0< \ve \le 1}  \|  \det \bigl(  \la DF^i (\ve, \,\cdot\,),DF^j (\ve, \,\cdot\,)\ra_{{\cal H}} \bigr)_{1 \le i,j \le n}^{-1} \|_q  <\infty
\qquad
\mbox{ for all $1<q<\infty$.}
\]
Here, $D$ stands for the ${\cal H}$-derivative.

The following is a straight forward  modification of Theorem 9.4, \cite{iwbk}.
 In this theorem, $0=\nu_0 < \nu_1 < \nu_2 < \cdots \nearrow \infty$ are 
all the elements of 
$$
\{  \kappa_{j_1}+ \cdots + \kappa_{j_n} ~|~ \, n =1,2,\ldots, \mbox{ and }  j_1, \ldots, j_n  \in {\bf N} \}
$$
 in increasing order.
\begin{tm}\label{tm.asym.plbk}
Let $F(\ve, \,\cdot\,) \in {\bf D}_{\infty} ({\bf R}^n) ~(0 < \ve \le 1)$ satisfy the following;
\[
F(\ve, \,\cdot\,) \sim  f_0 + \ve^{\kappa_1} f_{\kappa_1} +\ve^{\kappa_2} f_{\kappa_2} + \cdots   
\qquad
\mbox{in ${\bf D}_{\infty} ({\bf R}^n)$ as $\ve \searrow 0$,}
\]
with $f_{\kappa_j} \in {\bf D}_{\infty} ({\bf R}^n)$ for all $j \in {\bf N}$.  
We also assume that $F$ is uniformly non-degenerate in the sense of Malliavin.
Then, for any $T \in {\cal S}' ({\bf R}^n)$,
$\Phi(\ve,w) :=T \circ F(\ve , w)$ has the following asymptotic expansion;
\[
\Phi(\ve, \,\cdot\,) \sim  \phi_0 + \ve^{\kappa_1} \phi_{\kappa_1} +\ve^{\kappa_2} \phi_{\kappa_2} + \cdots   
\qquad
\mbox{in $\tilde{{\bf D}}_{-\infty}$ as $\ve \searrow 0$,}
\]
where $\phi_{\kappa_j} \in \tilde{{\bf D}}_{- \infty}$ is determined by a formal Taylor expansion 
as follows;  
\[
\Phi(\ve, \,\cdot\,) 
=
\sum_{\alpha} \frac{1}{\alpha!} (\partial^{\alpha} T) (f_0) 
[ \ve^{\kappa_1} f_{\kappa_1} +\ve^{\kappa_2} f_{\kappa_2} + \cdots ]^{\alpha}
=
  \phi_0 + \ve^{\nu_1} \phi_{\nu_1} +\ve^{\nu_2} \phi_{\nu_2} + \cdots,
\]
where the (formal) summation is over all multi-indexes $\alpha=(\alpha_1, \ldots ,\alpha_n) \in {\bf N}^n$.
(We set $\partial^{\alpha}= \prod_j (\partial / \partial x^j )^{\alpha_j}$ 
and $b^{\alpha} = \prod_j b_j^{\alpha_j}$ for $b=(b_1 , \ldots, b_n) \in {\bf R}^n$ as usual.)
For instance, $\phi_0 = T(f_0)$ and $\phi_{\kappa_1} =\sum_{j=1}^n f_{\kappa_1}^j 
\cdot (\partial T / \partial x^j )(f_0)$ and so on.
\end{tm}


Unlike in the usual stochastic analysis, 
almost every Wiener functional  in this paper is  continuous with respect to the topology of an abstract Wiener space,
because we work in the framework of Young integration.
Therefore, the following proposition will be very useful.  
For Banach spaces ${\cal X}_1, \ldots, {\cal X}_m, {\cal Y}$, $L^{(m)} ({\cal X}_1, \ldots, {\cal X}_m ; {\cal Y})$
denotes the space of bounded $m$-multilinear maps from ${\cal X}_1\times \cdots \times {\cal X}_m$ to ${\cal Y}$.

\begin{pr}\label{pr.muti.kuo}
Let $(W, {\cal H}, \mu)$ be an abstract Wiener space. 
Then, we have the following bounded inclusions;
\begin{eqnarray*}
L^{(m)} ( \underbrace{W, \ldots, W }_{m}; {\bf  R})   \hookrightarrow L^{(m)} (\underbrace{W, \ldots, W }_{m-1} , {\cal H}; {\bf  R})  
  \hookrightarrow  ({\cal H}^* )^{\otimes m}. 
\end{eqnarray*}
Here, the tensor product on the right hand side is Hilbert-Schmidt as usual.
\end{pr}

\Proof
The left bounded inclusion is obvious.
The right one is in p. 103, Kuo \cite{kuo}.
\QED


\section{Some results on  Malliavin calculus
for the solution of Young ODE driven by fBm with $H>1/2$}

In this section we discuss the solution of Young ODE 
driven by fBm with Hurst parameter $H \in (1/2, 1)$. 
We give moment estimates for the derivatives of the solution 
and prove uniform non-degeneracy of Malliavin covariance matrix of the solution.


Take $\alpha \in (1/2, H)$.
We denote by $\mu =\mu^H$ the law of $d$-dimensional fBm starting at $0$.
This Gaussian measure is supported in $C_0^{\alpha-hld} ([0,1] ; {\bf R}^d)$.
Cameron-Martin space is denoted by ${\cal H}={\cal H}^H$.
We set ${\cal W}$ to be the closure of ${\cal H}$ in $C_0^{\alpha-hld} ([0,1] ; {\bf R}^d)$.
Then, $({\cal W}, {\cal H}, \mu)$ becomes an abstract Wiener space.
(Note that the separable Hilbert space  ${\cal H}$ is not dense in $C_0^{\alpha-hld} ([0,1] ; {\bf R}^d)$, which is not separable.)
We denote by $(w_t)_{0 \le t \le 1}= (w_t^H)_{0 \le t \le 1}$ the canonical realization of fBm.

From now on,  we assume that $\sigma: {\bf R}^n \to {\rm Mat} (n,d)$ and $b: {\bf R}^n \to {\bf R}^n$ 
are  $C_b^{\infty}$.
We recall Young SDE (\ref{main.ygODE.eq}) driven by fBm $(w_t)$ in the following form;
\begin{equation}\label{ygSDE.eq}
dy_t = \sigma(y_t) dw_t + b(y_t) dt \qquad  \mbox{ with }  \qquad y_0 =a. 
\end{equation}
Then $y(w) =I (w, \lambda)$, where  $\lambda_t =t$ and $I$ is the It\^o map 
corresponding to the coefficients $[\sigma; b]=[ V_1, \ldots,V_d; V_0]$.
$I$ is everywhere-defined  and  continuous from $C_0^{\alpha-hld} ([0,1] ; {\bf R}^{d+1})$ to $C^{\alpha-hld} ([0,1] ; {\bf R}^d)$,
as we have explained in Section \ref{sec.yint}.

Moreover, $I$ is smooth in Fr\'echet sense (See Li and Lyons \cite{ll}) and, 
in particular,  $y=I ( \,\cdot \, , \lambda)$
is infinitely differentiable in
${\cal H}$-direction (see Nualart and Saussereau \cite{ns}). 
These are deterministic results.
In the sense of Malliavin calculus,  it is shown in  Hu and Nualart \cite{hn}  
that  $y_T : {\cal W} \to {\bf R}^n$ is ${\bf D}_{\infty}$ for any $T \in [0,1]$.

We can obtain 
an explicit form of the directional derivative $\xi_t^{h} := D_h y_t ~(h \in {\cal H})$ 
by differentiating (\ref{ygSDE.eq}); 
\begin{equation}\label{der1_SDE.eq}
d\xi_t^{h} - \nabla \sigma (y_t) \la \xi_t^{h}, dw_t  \ra  -   \nabla b (y_t)\la\xi_t^{h}\ra dt =   \sigma (y_t)dh_t
 \qquad  \mbox{ with }  \qquad
\xi_0^{h} =  0,
\end{equation}
or equivalently, 
\begin{equation}\label{der2_SDE.eq}
\xi_T^{h} = J(w,\lambda)_T \int_0^T   J(w, \lambda)^{-1}_t   \sigma (y_t) dh_t.
\end{equation}
Note that all the integrations above are in the Young sense.
An ODE for $J=J(w,\lambda)$ is given in (\ref{jODE1.eq}).
Let $h,k \in {\cal H}$.
By differentiating the above ODE, we see that  $\xi_t^{k,h} := D_k D_h y_t $ 
satisfies the following ODE;
\begin{align}
\lefteqn{
d\xi_t^{k,h} - \nabla \sigma (y_t) \la \xi_t^{k,h}, dw_t  \ra  -   \nabla b (y_t)\la\xi_t^{k,h}\ra dt 
=
 \nabla^2  \sigma (y_t) \la \xi_t^{k},   \xi_t^{h},   dw_t  \ra  
 }
\nn\\
&+  \nabla \sigma (y_t) \la \xi_t^{h}, dk_t  \ra 
 + \nabla \sigma (y_t) \la \xi_t^{k}, dh_t  \ra + \nabla^2  b (y_t) \la \xi_t^{k},   \xi_t^{h}\ra  dt
  \quad  \mbox{ with }  
\xi_0^{k,h} =  0.    \label{der3_SDE.eq}
\end{align}
Equivalently, we have
\begin{align}
\xi_T^{k,h}  &= J(w,\lambda)_T \int_0^T   J(w,\lambda)^{-1}_t   
  \bigl\{
   \nabla^2  \sigma (y_t) \la \xi_t^{k},   \xi_t^{h},   dw_t  \ra 
   \nn\\
   & \qquad  
   +  \nabla \sigma (y_t) \la \xi_t^{h}, dk_t  \ra 
 + \nabla \sigma (y_t) \la \xi_t^{k}, dh_t  \ra + \nabla^2  b (y_t) \la \xi_t^{k},   \xi_t^{h}\ra  dt
\bigr\}.
\label{der4_SDE.eq}
\end{align}
We can also obtain higher order directional derivatives in a similar way, but we omit them.


In a proof for the main theorem, we need to consider $\tilde{y}^{\ve} (w) = I (\ve w + \gamma, \ve^{1/H}  \lambda)$, 
where $\gamma \in {\cal H}$ is a fixed element and $\ve \in (0,1]$.
This process satisfies the following Young SDE; 
\begin{equation}\label{sc1_ygSDE.eq}
d\tilde{y}_t^{\ve}  = \sigma( \tilde{y}_t^{\ve})  \ve dw_t +  \sigma( \tilde{y}_t^{\ve})  d\gamma_t +  b( \tilde{y}_t^{\ve}) \ve^{1/H} dt 
\qquad  \mbox{ with }  \qquad \tilde{y}_0^{\ve} =a. 
\end{equation}
When $\gamma =0$,  we write $y^{\ve}$ for $\tilde{y}^{\ve}$.
In that case, 
self-similarity of $(w_t)$ implies that the two processes
 $(y_{  \ve^{1/H} t })_{0 \le t \le 1}$ and $(y_t^{\ve} )_{0 \le t \le 1}$ have the same law.

In the next proposition
we give  estimates for the derivatives $D^k  \tilde{y}_T^{\ve}$.
As we stated above, it is known that $y_T$ (and hence $\tilde{y}_T^{\ve}$)
is ${\bf D}_{\infty}$.
In that sense, this proposition is not new.
But, the estimate in powers of $\ve$  in (\ref{D^ky.ineq}) may be new.
Also, the proof is slightly different from the preceding papers,
because Proposition \ref{pr.muti.kuo} is used.
%
%
%
\begin{pr}\label{pr.est.Dy}
Take any $\gamma \in {\cal H}$ and fix it.
Then, for any $q \in (1, \infty)$ and $k =0,1,2,\ldots,$
there exists a positive constant $C_{q,k}$ such that 
\begin{equation}\label{D^ky.ineq}
 {\mathbb E}  [ \| D^k   \tilde{y}_T^{\ve}  \|^q_{( {\cal H}^*)^{\otimes k}}  ]^{1/q}  \le C_{q,k} \ve^k    \qquad \mbox{  for all $\ve \in (0,1]$ and $T \in [0,1]$. }
\end{equation}
\end{pr}

\Proof
In this proof, an unimportant positive constant $C$ may change from line to line.
First, consider the case $k=0$.
Since $\omega(s,t) = ( \|  w\|^p_{\alpha -hld} +  \|  \gamma\|^p_{\alpha -hld}  +1)  (t-s)$ satisfies
\[
| (\ve w_t  + \gamma_t)  - (\ve w_s  + \gamma_s)  | +  | \ve^{1/H} t -  \ve^{1/H}  s|  \le \omega(s,t)^{1/p},  \qquad (0 \le s \le t \le 1,~ p=1/\alpha),
\]
we can use a well-known estimate for the solutions of Young ODEs  to obtain that
\begin{equation}
| \tilde{y}^{\ve}_T | \le  \|   \tilde{y}^{\ve}  \|_{\alpha -hld}   \le |a|+  C (1+   \|  w\|^p_{\alpha -hld} +  \|  \gamma\|^p_{\alpha -hld} )
 \label{hoel1.ineq}
\end{equation}
 for some constant $C=C_K$.
 Fernique's theorem immediately implies (\ref{D^ky.ineq}) for $k=0$.

 %
 
Next let us consider the case $k=1$.  
By slightly modifying (\ref{der1_SDE.eq})--(\ref{der2_SDE.eq}), we can easily see that 
$\tilde{\xi}^{\ve; h}_t := D_h \tilde{y}^{\ve}_t$ satisfies the following (\ref{der5_SDE.eq})--(\ref{der6_SDE.eq}); 
\begin{equation}\label{der5_SDE.eq}
d \tilde{\xi} _t^{\ve; h} - \nabla \sigma (  \tilde{y}^{\ve}_t ) \la \tilde{\xi}_t^{\ve; h}, d(\ve w_t +\gamma_t)  \ra  
-   \nabla b  ( \tilde{y}^{\ve}_t)\la  \tilde{\xi}_t^{\ve; h}\ra    \ve^{1/H} dt
 =  
  \sigma (\tilde{y}^{\ve}_t  ) \ve dh_t
 \quad  \mbox{ with }  \quad
\tilde{\xi}_0^{\ve; h} =  0,
\end{equation}
or equivalently, 
\begin{equation}\label{der6_SDE.eq}
 \tilde{\xi} _T^{\ve; h}
 = \tilde{J}_T \int_0^T   \tilde{J}^{-1}_t   \sigma (   \tilde{y}^{\ve}_t)  \ve dh_t,
\end{equation}
where $\tilde{J} = J(\ve w +\gamma,  \ve^{1/H} \lambda )$.
From this, we can easily see that
\begin{align}
|  \tilde{\xi} _T^{\ve; h}| \le   \|  \tilde{\xi}^{\ve; h} \|_{\alpha -hld}  
&\le 
C \ve \| \tilde{J} \|_{\infty} \|  \tilde{J}^{-1}_{\cdot}  \sigma (   \tilde{y}^{\ve}_{\cdot})   \|_{\alpha -hld} \|h \|_{\alpha -hld}
\nn\\
&
\le
C \ve  \| \tilde{J} \|_{\infty} \|  \tilde{J}^{-1}  \|_{\alpha -hld}  
(1+\|  \tilde{y}^{\ve}   \|_{\alpha -hld} )\|h \|_{{\cal H}} 
\label{hoel2.ineq}
\end{align}
and, hence, 
\begin{align}
 \| D  \tilde{y}_T^{\ve}  \|_{ {\cal H}^*}  \le 
 C \ve  \| \tilde{J} \|_{\infty} \|  \tilde{J}^{-1}  \|_{\alpha -hld}  
(1+ \|  \tilde{y}^{\ve}   \|_{\alpha -hld}). 
\nn
  \end{align}
 By a slight modification of Proposition \ref{pr.momJ}, 
 $L^q$-norm of $ \|  \tilde{J}^{\pm 1}  \|_{\alpha -hld}$ is finite and  bounded 
 in $\ve$ for any fixed $q \in (1,\infty)$.
 (Just replace $w$ and $\lambda$ in Proposition \ref{pr.momJ}  
 by $\ve w +\gamma$ and $\ve^{1/H} \lambda$, respectively.)
 Hence, using H\"older's inequality, we obtain  (\ref{D^ky.ineq}) for $k=1$.

 We prove the case $k=2$. Set $D_k D_h  \tilde{y}_T^{\ve}  =\tilde{\xi} _T^{\ve; k, h}$
for simplicity.
Then, in the same way as in  (\ref{der3_SDE.eq})--(\ref{der4_SDE.eq}),
we have
\begin{align}
\tilde{\xi}_T^{\ve; k,h}  &= \tilde{J}_T \int_0^T   \tilde{J}^{-1}_t   
  \bigl\{
   \nabla^2  \sigma (\tilde{y}^{\ve}_t) \la \tilde{\xi}_t^{\ve; k},   \tilde{\xi}_t^{\ve; h},   d(\ve w_t +\gamma_t) \ra 
   \nn\\
   & \qquad  
   +  \nabla \sigma ( \tilde{y}^{\ve}_t ) \la \tilde{\xi}_t^{\ve; h}, \ve dk_t  \ra 
 + \nabla \sigma ( \tilde{y}^{\ve}_t ) \la \tilde{\xi}_t^{\ve ; k}, \ve dh_t  \ra 
+ \nabla^2  b (\tilde{y}^{\ve}_t ) \la \tilde{\xi}_t^{\ve; k},   \tilde{\xi}_t^{\ve; h}\ra  dt
\bigr\}.
\label{der7_SDE.eq}
\end{align}
From this, we have 
\begin{align}
\|  \tilde{\xi}^{\ve; k,h}\|_{\alpha -hld}   
&\le 
C \| \tilde{J} \|_{\infty} \|  \tilde{J}^{-1}  \|_{\alpha -hld} 
\bigl\{
\|\nabla^2 \sigma (   \tilde{y}^{\ve}_{\cdot})   \|_{\alpha -hld} 
\|   \tilde{\xi} _T^{\ve; k}  \|_{\alpha -hld} \|   \tilde{\xi} _T^{\ve; h}  \|_{\alpha -hld} 
\| \ve w +\gamma \|_{\alpha -hld} 
\nn\\
&\quad
+
\|\nabla \sigma (   \tilde{y}^{\ve}_{\cdot})   \|_{\alpha -hld} 
 (\|   \tilde{\xi}^{\ve; k}  \|_{\alpha -hld}   \|\ve  h \|_{\alpha -hld}  
+
\|   \tilde{\xi}^{\ve; h}  \|_{\alpha -hld}   \| \ve k \|_{\alpha -hld}  )
\nn\\
&\quad
+
\|\nabla^2 \sigma (   \tilde{y}^{\ve}_{\cdot})   \|_{\alpha -hld} 
\|   \tilde{\xi}^{\ve; k}  \|_{\alpha -hld} \|   \tilde{\xi}^{\ve; h}  \|_{\alpha -hld} 
\bigr\}
\nn\\
&
\le  C\ve^2 \| \tilde{J} \|_{\infty}^3   \|  \tilde{J}^{-1}  \|_{\alpha -hld}^3
 (1+\| w \|_{\alpha -hld}^p +\|\gamma \|_{\alpha -hld}^p  )^4
\|h \|_{\alpha -hld} \|k \|_{ \alpha -hld}.
\label{hoel3.ineq}
\end{align}
Here, we used (\ref{hoel1.ineq}) and (\ref{hoel2.ineq}) 
From Proposition \ref{pr.muti.kuo}, we see 
\begin{align}
\|  \tilde{\xi}_T^{\ve; k,h}\|_{ {\cal H}^* \otimes{\cal H}^*}
\le   
 C\ve^2 \| \tilde{J} \|_{\infty}^3   \|  \tilde{J}^{-1}  \|_{\alpha -hld}^3
 (1+\| w \|_{\alpha -hld}^p +\|\gamma \|_{\alpha -hld}^p  )^4.
\end{align}
Using the moment estimate for $ \|  \tilde{J}^{\pm 1}  \|_{\alpha -hld}$ again, 
we show (\ref{D^ky.ineq}) for $k=2$.

Finally, we briefly explain the higher order cases $(k \ge 3)$. 
We can  show it in a similar way by induction.
(We assume $\alpha$-H\"older norm of $D^{m}_{h_1, \ldots, h_{m}} \tilde{y}^{\ve}$
is dominated by  $\prod_{j=1}^{m} \|h_j\|_{ \alpha -hld} \times O(\ve^m)$ 
in any $L^q$-sense for $m \le k-1$ (as in (\ref{hoel3.ineq}) for $m=2$)
and then we will prove that $\alpha$-H\"older norm of
$D^{k}_{h_1, \ldots, h_{k}} \tilde{y}^{\ve}$ also does.)

For simplcity,  set
$\tilde{\eta}^{\ve}_t = D^k_{h_1, \ldots, h_k} \tilde{y}^{\ve}_t$.
It satisfies the following simple linear ODE similar to (\ref{der5_SDE.eq});
\[
d \tilde{\eta} _t^{\ve} - \nabla \sigma (  \tilde{y}^{\ve}_t ) \la \tilde{\eta}_t^{\ve}, d(\ve w_t +\gamma_t)  \ra  
-   \nabla b  ( \tilde{y}^{\ve}_t)\la  \tilde{\eta}_t^{\ve}\ra    \ve^{1/H} dt
 =  
  dG_t^{\ve}
 \quad  \mbox{ with }  \quad
\tilde{\eta}_0^{\ve} =  0.
\]
Here, $G^{\ve}$ is of the form 
\[
G_t^{\ve} = G^{\ve} \bigl(  \tilde{y}^{\ve}, D_{h_{j_1}}  \tilde{y}^{\ve} ,  \ldots, 
 D^{k-1}_{h_{j_1}, \ldots, h_{j_{k-1}} }   \tilde{y}^{\ve},  w, \gamma, h_1, \ldots, h_k  \bigr)_t
\]
and is of order $k$ in $\ve$.
Note that there is no derivative of order $k$ on the right hand side.
As in (\ref{der7_SDE.eq}), we have 
$
\tilde{\eta}_T^{\ve}  = \tilde{J}_T \int_0^T   \tilde{J}^{-1}_t dG_t^{\ve}.$
Using this we can estimate $\alpha$-H\"older norm of $\tilde{\eta}^{\ve}$ for $k$ in the same way as in (\ref{hoel3.ineq}).
  \QED

\begin{re}
We already have {\rm (i)} Fr\'echet smoothness of $\tilde{y}^{\ve}_T$ in the deterministic sense
 and
{\rm (ii)} $L^q$-estimates for derivatives as in this proposition.
From these, we can easily verify that $\tilde{y}^{\ve}_T \in {\bf D}_{\infty}$ as follows.
(For simplicity of notations, we only consider the case $\gamma=0, \ve=1$.)
By using Taylor expansion, we have  
\[
\frac{ y_T(w + r h) -y_T(w) }{r} -  D_hy_T(w)
=
r \int_0^1 d\theta (1-\theta) D^2 y_T(w + r \theta h) \la h,h\ra 
\]
for all $w \in {\cal W}$, $h \in {\cal H}^H$, and $r \in {\bf R}$.
Note that the derivative $D$ on the both sides of
 the above equation is in the deterministic sense.
By Proposition \ref{pr.est.Dy} and Cameron-Martin formula, 
the right hand side is $O(r)$ as $r \to 0$
in $L^q$-norm for any $q \in (1, \infty)$.
This implies that $y_T \in {\bf D}_{q,1}$ for any $q \in (1, \infty)$
and the derivative $D y_T$ in  (deterministic) Fr\'echet sense
is also the derivative in the sense of Malliavin calculus.
(See Proposition 4.21, \cite{shbk} for instance.)
The higher order derivatives can be dealt with in the same way.
\end{re}


Now we show that,  under the ellipticity condition {\bf (A1)}
for $\sigma$ (i.e., for $V_1 , \ldots , V_d$), 
the Malliavin covariance matrix for 
\begin{equation}\label{undeg.eq}
   \frac{ \tilde{y}^{\ve}_1 - a' }{\ve}  
\end{equation}
is uniformly non-degenerate in the sense of Malliavin as $\ve \searrow 0$.
Here,  
we set  $a' =\phi_1^0$ 
for the solution of the following ODE; 
$d\phi^0_t = \sigma( \phi^0_t) d \gamma_t$ with $\phi^0_0=a$.
\\
\\
{\bf (A1):}~ 
The set of vectors $\{V_1(a), \ldots, V_d(a)\}$ linearly spans ${\bf R}^n$.
\\
%
%
%
%
%

Nualart and Saussereau \cite{ns} showed
non-degenarcy of Malliavin covariance matrix of $y_T$ under 
{\bf (A1)}.
Baudoin and Hairer \cite{bh},
  proved non-degeneracy under H\"ormander's hypoellipticity condition for vector fields 
$\{ V_1 , \ldots , V_d ; V_0\}$.

In the next proposition,
we will prove uniform non-degeneracy of (\ref{undeg.eq}) under {\bf (A1)}
 by slightly modifying Baudoin-Hairer's argument.
(The special case $\gamma =0$ has already appeared in Baudoin and Ouyang \cite{bo2}.)
%
%
%
%
%
%
\begin{pr}\label{pr.unif.nondeg}
Let $\tilde{y}^{\ve} = ( \tilde{y}^{\ve, 1}, \ldots , \tilde{y}^{\ve, n})$ be the solution of (\ref{sc1_ygSDE.eq})
and assume {\bf (A1)}.
Then, $(\tilde{y}^{\ve}_1 - a' )/\ve$ is uniformly non-degenerate in the sense of Malliavin 
as $\ve \searrow 0$.
\end{pr}

\Proof
Let $y=(y_t)$ be the solution of (\ref{ygSDE.eq}).
In p. 388-389, \cite{bh},  
an explicit form of  the Malliavin covariance matrix for $y_1$ is given.
By replacing the vector fields $[\sigma; b]=[ V_1, \ldots,V_d; V_0]$ with $[\ve \sigma;  \ve^{1/H }b]=[ \ve V_1, \ldots, \ve V_d; \ve^{1/H } V_0]$, 
we can easily see that
\begin{align}
\lefteqn{
\frac{1}{\ve^2} 
\bigl( \la Dy^{\ve; i}_1 , Dy^{\ve; j}_1 \ra_{{\cal H}}  \bigr)_{1 \le i,j \le n} 
=
H(2H-1)  J(\ve w,  \ve^{1/H} \lambda )_1  
}
\nn\\
&
\times
\int_0^1  \int_0^1        
 J(\ve w,  \ve^{1/H} \lambda )_u^{-1}
 \sigma (y^{\ve}_u)   \sigma (y^{\ve}_v)^* 
   J(\ve w,  \ve^{1/H} \lambda )_v^{-1,*}  
   |u-v|^{2H-2}  dudv
J(\ve w,  \ve^{1/H} \lambda )_1^*.
\nn
\end{align}
Here, $\lambda_t =t$ and  $A^*$ denotes the transposed matrix of $A$.
By shifting $w \mapsto w + (\gamma /\ve)$, we have 
\begin{align}
\frac{1}{\ve^2} 
\bigl( \la D\tilde{y}^{\ve; i}_1 , D\tilde{y}^{\ve; j}_1 \ra_{{\cal H}}  \bigr)_{1 \le i,j \le n} 
=
H(2H-1)  \tilde{J}_1    
  C
 \tilde{J}_1^*,   
  \end{align}
where $\tilde{J}_t = J(\ve w +\gamma,  \ve^{1/H} \lambda )_t$ as before and  we set 
\[
C = \int_0^1  \int_0^1        
 \tilde{J}_s^{-1}
 \sigma (\tilde{y}^{\ve}_s)      \sigma (\tilde{y}^{\ve}_t)^* 
   \tilde{J}_t^{-1,*}  
   |s-t|^{2H-2}  dsdt.
    \]
Since  
$\sup_{0<\ve \le 1 } \| \tilde{J}^{\pm1}_1 \|_q <\infty$ for any $q \in (1, \infty)$,  
it is sufficient to prove 
\begin{equation} \label{suff.hypo.ineq}
\sup_{0<\ve \le 1 } \| |\det C|^{-1} \|_q <\infty   \qquad \mbox{ for any $1< q <\infty$}.
\end{equation}

We will follow the argument in pp. 387--340, \cite{bh}.
In order to show (\ref{suff.hypo.ineq}) above, 
it is sufficient to prove that, for any $1< q <\infty$,  there exists $\rho_0 (q)$, 
which is independent of $\ve$ and satisfies that, 
\begin{equation} \label{suff2.hypo.ineq}
\sup_{ {\bf a} \in {\bf R}^n,   \| {\bf a}\| =1  }
\mu \bigl( \la  {\bf a}, C {\bf a} \ra    \le \rho \bigr)
\le \rho^q  
\qquad \mbox{ for any $\rho \in (0, \rho_0 (q))$ and $\ve \in (0,1]$}.
\end{equation}
(For a proof, see Lemma 2.3.1,  Nualart \cite{nbk}).
As in \cite{bh},  
\begin{align}
\la {\bf a}, C {\bf a} \ra  
&= 
 \sum_{j=1}^d  
\int_0^1  \int_0^1        
\la   {\bf a}, \tilde{J}_s^{-1}  V_j  (\tilde{y}^{\ve}_s)     \ra
\la   {\bf a}, \tilde{J}_t^{-1}  V_j  (\tilde{y}^{\ve}_t)     \ra
   |s-t|^{2H-2}  dsdt
   \nn
   \\
   &= 
    \sum_{j=1}^d 
      \bigl\|  
       \la   {\bf a}, \tilde{J}_{\cdot}^{-1}  V_j  (\tilde{y}^{\ve}_{\cdot}  )     \ra 
        \bigr\|^2_{ \hat{\cal H}}.
       \end{align}
By a Norris-type lemma (Corollary 4.5, \cite{bh}),   there exists $0< \beta <1/2$ such that 
for any $r < H-(1/2)$ and $0< \rho \le 1$, 
the following inequalities hold;
\begin{align}
\lefteqn{
\mu \bigl( \la  {\bf a}, C {\bf a} \ra    \le \rho \bigr)
\le 
\min_{1 \le j \le d}  
\mu \bigl(   \|  \la   {\bf a}, \tilde{J}_{\cdot}^{-1}  V_j  (\tilde{y}^{\ve}_{\cdot}  )  
   \ra  \bigr\|_{\hat{\cal H}}    \le  \rho^{1/2} \bigr)
}
\nn\\
&\le
\min_{1 \le j \le d}  
\Bigl[
\mu \bigl(   \|  \la   {\bf a}, \tilde{J}_{\cdot}^{-1}  V_j  (\tilde{y}^{\ve}_{\cdot}  )     \ra  \bigr\|_{L^{\infty}}   <\rho^{\beta /2} \bigr)
+
\mu \bigl(   \|  \la   {\bf a}, \tilde{J}_{\cdot}^{-1}  V_j  (\tilde{y}^{\ve}_{\cdot}  )     \ra  \bigr\|_{r -hld}   >\rho^{ - \beta /2} \bigr)
\Bigr]
\nn\\
&\le
\min_{1 \le j \le d}  
\Bigl[
\mu \bigl(   |  \la   {\bf a},   V_j  (a  )    \ra  |   <\rho^{\beta /2} \bigr)
+
\mu \bigl(   \|  \la   {\bf a},      \tilde{J}_{\cdot}^{-1}  V_j  (\tilde{y}^{\ve}_{\cdot}  )       \ra  \bigr\|_{\alpha -hld}   >\rho^{ - \beta /2} \bigr)
\Bigr].
\label{suff3.hypo.ineq}
\end{align}
Here, in the last inequality, we evaluated at $t=0$ and used $r <\alpha$.
Note that the set in the first term on the right hand side 
is already independent of $\ve$ and non-random 
(i.e., either $\emptyset$ or the whole set ${\cal W}$).

Recall  that, for any $q$,
$
E[  \|  \tilde{J}^{-1}  \|^q_{\alpha -hld}  +  \|  \tilde{y}^{\ve} \|^q_{\alpha -hld}   ]   \le  c_1
$
for some constant $c_1=c_1(q)$ which is independent of $\ve$.
Then, using  Chebyshev's inequality, we have
\[
\mu \bigl(   \|  \la   {\bf a}, \tilde{J}_{\cdot}^{-1}  V_j  (\tilde{y}^{\ve}_{\cdot}  )     \ra
  \bigr\|_{\alpha -hld}   >\rho^{ - \beta /2} \bigr)
\le 
c_2  \rho^{ q}
\]
for some constant $c_2=c_2 (q)$ which is independent of $\ve$.

Let us consider the first term on the right hand side of (\ref{suff3.hypo.ineq}).
By {\bf (A1)}, there exists $c' >0$ such that $\sigma(a) \sigma(a)^* \ge c' {\rm Id}_n$
in the form sense.
We have
\begin{align}
n \max_{1 \le j \le d}  |  \la   {\bf a},   V_j  (a  )     \ra  |^2 
\ge 
\sum_{1 \le j \le d}  |  \la   {\bf a},   V_j  (a  )     \ra  |^2
= \la  {\bf a},   \sigma(a) \sigma(a)^*  {\bf a}\ra \ge c' >0.
\nn
\end{align}
Hence, if $\rho^{\beta /2} \le \sqrt{ c'/n}$,  
then $\min_{1 \le j \le d}  \mu \bigl(   |  \la   {\bf a},   V_j  (a  )    \ra  |   <\rho^{\beta /2} \bigr)=0$
and
\[
\mu \bigl( \la  {\bf a}, C {\bf a} \ra    \le \rho \bigr)
\le 
c_2 (q) \rho^{ q}.
\]
From this, we can easily see (\ref{suff2.hypo.ineq}) holds with 
$\rho_0 (q) = c_2 (q+1)^{-1} \wedge  (c' /n)^{1/ \beta}$.
This completes the proof.
\QED

\begin{re}
In the above proof, $\hat{\cal H}$ is 
another Hilbert space that is unitarily isometric to ${\cal H}$.
Loosely speaking, it is defined as follows:
Denote by ${\cal E}$ the set of ${\bf R}^d$-valued step functions on $[0,1]$.
Let $\hat{\cal H}$ be the Hilbert space defined as the closure of ${\cal E}$
by the inner product 
\[
\la I_{[0,s]}  v , I_{[0,t]} w \ra_{ \hat{\cal H} }  = R(s,t)  \la v,w\ra_{{\bf R}^d },
\qquad (t,s \in [0,1], \,  v,w\in {\bf R}^d), 
\]
where we set $R(s,t) ={\mathbb E}[w_s^i w_t^i]$.
(For inastance, see Section 5.1, Nualart \cite{nbk} or \cite{bh, ns}
for more information on $\hat{\cal H}$.)
\end{re}


\section{ On-diagonal short time asymptotics}

The aim of this section 
is to prove Theorem \ref{thm.MAIN.on}, namely, 
on-diagonal short time asymptotic expansion of  the density of the solution
of the Young SDE (\ref{main.ygODE.eq})
(or equivalently (\ref{ygSDE.eq})) under the ellipticity assumption {\bf (A1)}.

Let us consider the solution $(y_t) = (y_t (a))$
of  Young differential equation (\ref{main.ygODE.eq}) with an initial condition $y_0 =a \in {\bf R}^n$ 
driven by fBm $(w_t)$ with $H >1/2$.
It is shown in \cite{ns, bh}
that, under {\bf (A1)}, the law of the solution has a smooth density $p(t, a, a')$, i.e., 
\[
{\mathbb P} \bigl( y_t (a)  \in A   \bigr)  = \int_A p(t, a, a') da'
\qquad
\mbox{( for any Borel set $A \subset {\bf R}^n$).}
\]
For $t>0$, 
$y_t =y_t (a)$ is ${\bf D}_{\infty}$ and  non-degenerate in the sense of Malliavin. 
By the same argument as in Ikeda and Watanabe \cite{iwbk}, 
we have the following expression; $p(t, a, a') = {\mathbb E}[ \delta_{a'}  ( y_t (a)) ] =
{}_{{\bf D}_{- \infty} } \la   \delta_{a'}  ( y_t (a))  , {\bf 1} \ra_{{\bf D}_{\infty} }$.
By the self-similarity of fBm,  
$(y_{ \ve^{1/H} t } )_{t \ge 0}$ and $(y^{\ve}_{ t } )_{t \ge 0}$ have the same law,
where $y^{\ve}$ is given by (\ref{sc1_ygSDE.eq}) with $\gamma =0$.
From this,  we see that 
$ p( \ve^{1/H}, a, a') =  {\mathbb E}[ \delta_{a'}  ( y^{\ve}_1 (a)) ]$.
%
%
%

The most important part of the proof 
 is an asymptotic expansion of $ y^{\ve}_1 $ in $\ve \in (0,1]$ in ${\bf D}_{\infty} $-topology.
For that purpose,  we introduce the following index set for exponent of $\ve$.  
Set
\[
\Lambda_1 = \{   n_1 + \frac{n_2}{H}  ~|~  n_1, n_2 \in {\bf N}  \}.
\]
We denote by $0=\kappa_0 <\kappa_1 < \kappa_2 < \cdots$ the elements of $\Lambda_1$ in increasing order.
Several smallest elements are explicitly given as follows;
\[
\kappa_1 = 1, \quad   \kappa_2 =  \frac{1}{H},  \quad   \kappa_3 = 2, \quad \kappa_4 = 1+\frac{1}{H},  
\quad  \kappa_5 = 3 \wedge \frac{2}{H}, \ldots
\]


\begin{pr} \label{pr.yve.dia}
The family of  Wiener functional $y^{\ve}_1 ~(0<\ve \le 1)$ admits 
the following asymptotic expansion as $\ve \searrow 0$; 
\[
y^{\ve}_1
 \sim  a + \ve f_{1} +\ve^{\kappa_2} f_{\kappa_2} +   \ve^{\kappa_3} f_{\kappa_3}  +\cdots   
\qquad
\mbox{in ${\bf D}_{\infty} ({\bf R}^n)$}
\]
for certain $ f_{\kappa_1}, f_{\kappa_2}, \ldots \in {\bf D}_{\infty} ({\bf R}^n)$.
\end{pr}

\Proof
For ${\bf j} = (j_1, \ldots, j_m) \in \{ 0,1, \ldots, d \}^m$, we set 
$|{\bf j}| = m$ and 
$$
\|{\bf j}\| =  \frac{ \sharp \{ 1 \le k  \le m ~|~ j_k =0\} }{H}  +  \sharp \{ 1 \le k  \le m ~|~ j_k \neq 0\}.
$$
We denote by ${\cal I}_m$  the totality of such ${\bf j}$'s with $|{\bf j}| = m$ and set ${\cal I} =\cup_{m=1}^{\infty}  {\cal I}_m$.

We will use  the following convention.  We set $t = w^0_t$. 
Then, the ODE for $y^{\ve}$ (that is, (\ref{sc1_ygSDE.eq}) with $\gamma =0$) reads;
\begin{equation}
dy^{\ve}_t
=
\ve^{1/H} V_0 ( y^{\ve}_t) dw^0_t  + \sum_{j=1}^d  \ve  V_j ( y^{\ve}_t) dw^j_t
\qquad
\mbox{with $y^{\ve}_0 =a$.}
\label{sc2_ygSDE.eq}
\end{equation}

It is easy to see that 
\begin{align}
y^{\ve}_1 -a 
&=
 \ve^{1/H}    \int_0^1   V_0 ( y^{\ve}_t) dw^0_t  + \sum_{j=1}^d  \ve \int_0^1   V_j ( y^{\ve}_t) dw^j_t
 \nn\\
  &= \ve^{1/H}    \int_0^1   V_0 ( a) dw^0_t  + \sum_{j=1}^d  \ve \int_0^1   V_j (a ) dw^j_t
       \nn\\
        & \quad +  \ve^{1/H}    \int_0^1   \{ V_0 ( y^{\ve}_t)  -   V_0 ( a)  \}  dw^0_t  
          + \sum_{j=1}^d  \ve \int_0^1   \{  V_j ( y^{\ve}_t)  - V_j ( a) \} dw^j_t
  \nn\\
  &= \ve^{1/H}   V_0 ( a)  + \sum_{j=1}^d  \ve    V_j (a ) w^j_1
  \nn\\
  & \quad
    +  \ve^{1/H}    \int_0^1    \int_0^{t_1}  \{   \ve^{1/H}  \hat{V}_0 V_0 ( y^{\ve}_{ t_2} )  dt_{2}
    + \sum_{j' =1}^d   \ve   \hat{V}_{j'} V_0 ( y^{\ve}_{ t_2} )  dw^j_{t_2}   \}  dw^0_{t_1}  
     \nn\\
     &\quad
     +
      \sum_{j =1}^d  \ve^{1+(1/H)}     \int_0^1 \int_0^{t_1}     \hat{V}_{0} V_j ( y^{\ve}_{ t_2} )  dw^{j'}_{t_2}    dw^j_{t_1} 
     +  \sum_{j, j' =1}^d  \ve^2     \int_0^1    \int_0^{t_1}     \hat{V}_{j'} V_j ( y^{\ve}_{ t_2} )  dw^{j'}_{t_2}    dw^j_{t_1}  
                    \nn\\
                    & = 
                    \ve^{1/H}   V_0 ( a)  + \sum_{j=1}^d  \ve    V_j (a ) w^j_1
                    + 
                    \sum_{ |{\bf j}| = 2 }    \ve^{   \|{\bf j}\| }    
                    \int_0^1    \int_0^{t_1}     \hat{V}_{j_2} V_{j_1} ( y^{\ve}_{ t_2} )  dw^{j_2}_{t_2}    dw^{j_1}_{t_1}.                                      
                 \label{iter.tenk1.eq}
                   \end{align}
Here, $\hat{V}_{i}  V_j $ denotes a vector field $V_i$ (as a first order differential operator) acting on a ${\bf R}^n$-valued 
function $V_j$.

Repeating the same argument for the last term on the right hand side of (\ref{iter.tenk1.eq}), we have 
      \begin{align}
y^{\ve}_1 -a 
&=    
 \ve^{1/H}   V_0 ( a)  + \sum_{j=1}^d  \ve    V_j (a ) w^j_1
                    + 
                       \sum_{ |{\bf j}| = 2 }    \ve^{   \|{\bf j}\| }    
                     \hat{V}_{j_2} V_{j_1} ( a )    \int_0^1    \int_0^{t_1}    dw^{j_2}_{t_2}    dw^{j_1}_{t_1}    
                    \nn\\
                    & \quad 
                    +
   \sum_{ |{\bf j}| = 3 }    \ve^{   \|{\bf j}\| }    
                    \int_0^1    \int_0^{t_1}    \int_0^{t_2}   
                          \hat{V}_{j_3}    \hat{V}_{j_2} V_{j_1} ( y^{\ve}_{ t_2} )  dw^{j_3}_{t_3}  dw^{j_2}_{t_2}    dw^{j_1}_{t_1}.     
                    \label{iter.tenk2.eq}
                   \end{align}
     Here, $   \hat{V}_{j_3}    \hat{V}_{j_2} V_{j_1} =    \hat{V}_{j_3}  (   \hat{V}_{j_2} V_{j_1} )$.     
          In general,  we have 
                \begin{align}
y^{\ve}_1 -a 
&=    
 \sum_{ 1 \le |{\bf j}| \le n-1}    \ve^{   \|{\bf j}\| }    
                  \hat{V}_{j_{n-1}}\cdots \hat{V}_{j_2} V_{j_1} ( a )  
  \int_0^1    \int_0^{t_1}   \cdots   \int_0^{t_{n-2}}   dw^{j_{n-1}}_{t_{n-1}} \cdots  dw^{j_2}_{t_2}    dw^{j_1}_{t_1}    
                    \nn\\
                    & \quad 
                    + 
\sum_{ |{\bf j}| = n }    \ve^{   \|{\bf j}\| }    
                    \int_0^1    \int_0^{t_1}   \cdots  \int_0^{t_{n-1}}   
                          \hat{V}_{j_n} \cdots  \hat{V}_{j_2} V_{j_1} ( y^{\ve}_{ t_n} )  dw^{j_n}_{t_n}  \cdots dw^{j_2}_{t_2}    dw^{j_1}_{t_1}.     
    \label{iter.tenk3.eq}
                   \end{align}
Let us observe the first term.
From  basic properties of Young integral,  
we easily see that, for any $m$,
the real-valued functional 
$ \int_0^1    \int_0^{t_1}   \cdots   \int_0^{t_{m-1}}   dw^{j_{m}}_{t_{m}} \cdots  dw^{j_2}_{t_2}    dw^{j_1}_{t_1}$ 
is in $m$th (inhomogeneous) Wiener chaos and hence it is in any ${\bf D}_{q,k}~(1<q<\infty, k \in {\bf N})$.

Next we consider the last term in (\ref{iter.tenk3.eq}).
We set 
\[
Q_{\ve}(w) = \int_0^1    \int_0^{t_1}   \cdots  \int_0^{t_{n-1}}   
                          \hat{V}_{j_n} \cdots  \hat{V}_{j_2} V_{j_1} ( y^{\ve}_{ t_n} )  dw^{j_n}_{t_n}  \cdots dw^{j_2}_{t_2}    dw^{j_1}_{t_1} 
\]
and will prove  $Q_{\ve} = O(1)$ as $\ve \searrow 0$ in ${\bf D}_{q,k}({\bf R}^n)$ for 
any $1<q<\infty, \, k \in {\bf N}$.
(For simplicity of notation, we denote $G=\hat{V}_{j_n} \cdots  \hat{V}_{j_2} V_{j_1}$ 
and assume $j_i \neq 0$ for all $i$.
The other case is actually easier.)

Since $\|y^{\ve}\|_{\alpha -hld}$ is $O(1)$ in any $L^q$ as $\ve\searrow 0$,
 $Q_{\ve}(w) $ is $O(1)$ in any $L^q$, too.
 Now we estimate the derivatives.
For $h \in {\cal H}$, we have
\begin{align}
D_h Q_{\ve}(w) 
&=
\int_0^1  \cdots  \int_0^{t_{n-1}}   
     \nabla G ( y^{\ve}_{ t_n} ) \la D_h y^{\ve}_{ t_n}\ra  dw^{j_n}_{t_n}  \cdots dw^{j_2}_{t_2}    dw^{j_1}_{t_1} 
\nn\\
& \quad
+\sum_{l=1}^n 
\int_0^1  \cdots  \int_0^{t_{n-1}} 
G ( y^{\ve}_{ t_n} )   dw^{j_n}_{t_n}  \cdots dh^{j_l}_{t_l}  \cdots  dw^{j_1}_{t_1}.
\nn
\end{align}
H\"older norms of $y^{\ve}$ and $ D_h y^{\ve}$ were estimated in (\ref{hoel1.ineq})--(\ref{hoel2.ineq}).
From these, we see that $\|D Q_{\ve} \|_{{\cal H}^*} =O(1)$ in any $L^q$.

Similarly,  $h,k  \in {\cal H}$, we have
\begin{align}
D_kD_h Q_{\ve}(w) 
&=
\int_0^1  \cdots  \int_0^{t_{n-1}}   
     \nabla G ( y^{\ve}_{ t_n} ) \la D_kD_h y^{\ve}_{ t_n}\ra  dw^{j_n}_{t_n}  \cdots dw^{j_2}_{t_2}    dw^{j_1}_{t_1} 
\nn\\
& \quad
+
\int_0^1  \cdots  \int_0^{t_{n-1}}   
     \nabla G ( y^{\ve}_{ t_n} ) \la D_ky^{\ve}_{ t_n},  D_h y^{\ve}_{ t_n}\ra 
 dw^{j_n}_{t_n}  \cdots dw^{j_2}_{t_2}    dw^{j_1}_{t_1} 
\nn\\
& \quad
+\sum_{l=1}^n 
\int_0^1  \cdots  \int_0^{t_{n-1}} 
\nabla G ( y^{\ve}_{ t_n} )  \la D_h y^{\ve}_{ t_n}\ra
 dw^{j_n}_{t_n}  \cdots dk^{j_l}_{t_l}  \cdots  dw^{j_1}_{t_1}
\nn\\
& \quad
+\sum_{l=1}^n 
\int_0^1  \cdots  \int_0^{t_{n-1}} 
\nabla G ( y^{\ve}_{ t_n} )  \la D_k y^{\ve}_{ t_n}\ra
 dw^{j_n}_{t_n}  \cdots dh^{j_l}_{t_l}  \cdots  dw^{j_1}_{t_1}
\nn\\
& \quad
+\sum_{l \neq m}^n 
\int_0^1  \cdots  \int_0^{t_{n-1}} 
G ( y^{\ve}_{ t_n} ) 
  dw^{j_n}_{t_n}  \cdots dh^{j_l}_{t_l} \cdots dh^{j_m}_{t_m} \cdots  dw^{j_1}_{t_1}.
\nn
\end{align}
H\"older norm of $D_k D_h y^{\ve}$ was estimated in (\ref{hoel3.ineq}).
Combined with Proposition \ref{pr.muti.kuo},
the above implies that  $\|D^2 Q_{\ve} \|_{{\cal H}^* \otimes {\cal H}^*} =O(1)$ in any $L^q$.
Higher order derivatives can be done in the same way.

Now we prove the proposition.  
In order to get the asymptotic expansion up to order $\kappa_m$ (i.e., the remainder is of order $\kappa_{m+1}$),  it is sufficient 
(i)~ to consider the expansion (\ref{iter.tenk3.eq}) 
with $n-1$ being the smallest integer which is not less than  $\kappa_m$
and (ii)~ to set 
\[
f_{\kappa_l}  (w)= \sum_{  \|{\bf j} \| =\kappa_l }   
                  \hat{V}_{j_{n}}\cdots \hat{V}_{j_2} V_{j_1} ( a )  
  \int_0^1    \int_0^{t_1}   \cdots   \int_0^{t_{n-1}}   dw^{j_{n}}_{t_{n}} \cdots  dw^{j_2}_{t_2}    dw^{j_1}_{t_1} 
\]
for all $1\le l \le m$.
\QED

Before we prove on-diagonal short time kernel asymptotics, we define two more index sets
for exponent of $\ve$.
Set 
$
\Lambda_2 = \{ \kappa -1 ~|~ \kappa \in \Lambda_1 \setminus \{0\} \}$.
Smallest elements of  $\Lambda_2$ are 
$$0,  \quad \frac{1}{H} -1,  \quad  1,   \quad \frac{1}{H}, 
  \quad \bigl(  3 \wedge \frac{2}{H}  \bigr)-1,  \ldots
$$
Next we set 
$
\Lambda_3 = 
\{  a_1 +a_2 + \cdots + a_m ~|~  \mbox{$m \in {\bf N}_+$ and $a_1 ,\ldots, a_m \in \Lambda_2$} \}$.
In the sequel, $\{ 0=\nu_0 <\nu_1<\nu_2 <\cdots \}$ stands for all the elements of 
$\Lambda_3$ in increasing order.


%
%

\vspace{4mm}
{\it Proof of Theorem \ref{thm.MAIN.on}}~
First, note that 
$$ 
p( \ve^{1/H}, a, a) =  {\mathbb E}[ \delta_{a}  ( y^{\ve}_1 (a)) ]
=
 {\mathbb E}[ \delta_{0}    \bigl( \ve \frac{ y^{\ve}_1 (a) -a}{\ve}  \bigr) ]
=
\ve^{-n} {\mathbb E}[ \delta_{0}    \bigl(  \frac{ y^{\ve}_1 (a) -a}{\ve}  \bigr) ].
$$

By Proposition \ref{pr.unif.nondeg}, $( y^{\ve}_1 (a) -a)/\ve$
is uniformly non-degenerate.
It admits asymptotic expansion in ${\bf D}_{\infty} ({\bf R}^n)$ as in Proposition \ref{pr.yve.dia}.
Then, by Theorem \ref{tm.asym.plbk}, the following asymptotic expansion holds in 
$\tilde{{\bf D}}_{- \infty}$ as $\ve \searrow 0$;
\[
\delta_{0}    \bigl(  \frac{ y^{\ve}_1 (a) -a}{\ve}  \bigr)
\sim
\phi_0 + \ve^{\nu_1} \phi_{\nu_1} + \ve^{\nu_2} \phi_{\nu_2} +\cdots
\qquad
\mbox{as $\ve \searrow 0$.}
\]
By taking the generalized expectation and setting $c_{\nu_k} = {\mathbb E}[\phi_{\nu_k} ]$,  we have
\[
p( \ve^{1/H}, a, a) 
\sim
\ve^{-n}
\bigl(
c_0 + 
c_{\nu_1}  \ve^{\nu_1} 
+ c_{\nu_2} \ve^{\nu_2}  +\cdots \bigr)
\qquad
\mbox{as $\ve \searrow 0$.}
\]
Putting $\ve = t^H$, we complete the proof of Theorem \ref{thm.MAIN.on}.
\toy
%


\section{Taylor expansion of It\^o map around a Cameron-Martin path}
\label{sec.taylor}

In this section we prove an asymptotic expansion for $\tilde{y}^{\ve} =I (\ve w +\gamma, \ve^{1/H} \lambda)$, 
which was defined in (\ref{sc1_ygSDE.eq}).
The base point $\gamma \in {\cal H}$ of the expansion is arbitrary, but  fixed.
First,  we prove that $\tilde{y}^{\ve}$ admits the following expansion in $C^{\alpha -hld} ([0,1]; {\bf R}^n)$;
\begin{equation}\label{asyphi1.eq}
\tilde{y}^{\ve} \sim \phi^0 + \ve^{\kappa_1} \phi^{\kappa_1} + \ve^{\kappa_2} \phi^{\kappa_2} 
+ \cdots    \qquad \mbox{as $\ve \searrow 0$,}
\qquad
(\kappa_i \in \Lambda_1={\bf N} +\frac{1}{H} {\bf N}),
\end{equation}
for some $C^{\alpha -hld} ([0,1]; {\bf R}^n)$-valued Wiener functional 
$\phi^0 ,   \phi ^{\kappa_1},  \phi^{\kappa_2}, \ldots$.
Since the It\^o map $I$ in the sense of Young integral equation is smooth in Fr\'echet sense (see \cite{ll}), 
this kind expansion holds in deterministic sense.
In this paper, however, we need to prove this expansion in $L^q$-sense.

Before we state the proposition precisely, 
 we now give a heuristic argument to find  an explicit form of $\phi^{\kappa_m}$.
 To find an ODE for $\phi^0$ is easy.
\begin{align}
d\tilde{y}_t^{\ve}  &= \sigma( \tilde{y}_t^{\ve})   (\ve dw_t + d\gamma_t )+  b( \tilde{y}_t^{\ve}) \ve^{1/H} dt 
 &\mbox{ with }  \qquad \tilde{y}_0^{\ve} =a,
\nn\\
d\phi^0_t  &= \sigma( \phi^0_t)  d\gamma_t    \qquad  & \mbox{ with }  \qquad \phi_0^{0} =a. 
\label{phi0.def}
\end{align}
Set $\triangle \phi : =  \tilde{y}^{\ve} - \phi^0$ 
and put it in the above ODE for $\tilde{y}^{\ve}$. Then we have 
\begin{align}
d(\phi^0 + \triangle \phi) 
&= \sigma( \phi^0 + \triangle \phi )   (\ve dw + d\gamma)+  b(\phi^0 + \triangle \phi) \ve^{1/H} dt
\nn\\
&= 
\sum_{k=0}^{\infty} \frac{ \nabla^k \sigma ( \phi^0) }{k! }
\la  \underbrace{  \triangle \phi, \ldots, \triangle  \phi }_{k} ; \ve dw + d\gamma\ra
+
\sum_{k=0}^{\infty}
\frac{  \nabla^k b ( \phi^0)  }{k!}
 \la\underbrace{   \triangle \phi, \ldots, \triangle \phi}_{k} \ra \ve^{1/H} dt.
\nn
\end{align}
Assume $\triangle \phi$ admits the asymptotic expansion  (\ref{asyphi1.eq}).
Then, by putting it in the above equation and picking up the terms of order $\ve^{\kappa_m}$,  
we find  an ODE for $\phi^{\kappa_m}$.
Note that $\phi^{\kappa_m}_0 =0$ for all $m \ge 1$.

For $\kappa_m = 1, 1/H, 2$,  we can write down the ODEs explicitly as follows;
\begin{align}
d\phi^{1}_t   -  \nabla \sigma( \phi^0_t) \la\phi^{1 }_t,     d\gamma_t \ra  
&= \sigma ( \phi^0_t) dw_t,
\label{phiODE2.def}
\\
d\phi^{1/H}_t   -  \nabla \sigma( \phi^0_t) \la\phi^{1 /H}_t  ,  d\gamma_t \ra  
&= b ( \phi^0_t) dt,
\label{phiODE3.def}
\\
d\phi^{2}_t   -  \nabla \sigma( \phi^0_t) \la\phi^{2 }_t   , d\gamma_t \ra  
&=   \nabla \sigma ( \phi^0_t)\la \phi^1_t,  dw_t \ra
+
  \frac12  \nabla^2 \sigma ( \phi^0_t)\la \phi^1_t,   \phi^1_t, d\gamma_t \ra.
\label{phiODE4.def}
\end{align}
Note that $\phi^{1/H}$ is independent of $w$, i.e, non-random with respect to $\mu$.

%
For $\kappa_m \ge 2$, 
\begin{align}\label{phiODE.def}
d\phi^{\kappa_m}_t   -  \nabla \sigma( \phi^0_t) \la\phi^{\kappa_{m} }_t    , d\gamma_t \ra  
&=
\sum_{k=1}^{\infty}  \sum_{ \kappa_{i_1} +\cdots + \kappa_{i_k} = \kappa_{m}  -1} 
 \frac{  \nabla^k \sigma ( \phi^0_t)  }{k!}  \la\phi^{\kappa_{i_1} }_t,   \ldots,  \phi^{\kappa_{i_k} }_t ;  dw_t \ra 
   \nn\\
   &\quad 
    + 
    \sum_{k=2}^{\infty}  \sum_{ \kappa_{i_1} +\cdots + \kappa_{i_k} = \kappa_{m}  } 
 \frac{  \nabla^k \sigma ( \phi^0_t)   }{k!}
\la \phi^{\kappa_{i_1} }_t,   \ldots,  \phi^{\kappa_{i_k} }_t ;  d \gamma_t \ra 
   \nn\\
   &\quad 
    + \sum_{k=1}^{\infty}  \sum_{ \kappa_{i_1} +\cdots + \kappa_{i_k} = \kappa_{m}  -(1/H)} 
     \frac { \nabla^k b ( \phi^0_t)  }{k!}
     \la\phi^{\kappa_{i_1} }_t,   \ldots,  \phi^{\kappa_{i_k} }_t ;  dt \ra.
\end{align}
%
%
%
The summations in the first term on  the right hand side is taken over all $ \kappa_{i_1},  \ldots, \kappa_{i_k} \in \Lambda_1 \setminus \{0\}$
such that $ \kappa_{i_1} +\cdots + \kappa_{i_k} = \kappa_{m}  -1$ hold.
$ \kappa_{i_j} =0$ is not allowed. So, the sum  is actually  a finite sum.
The second and the third terms should be understood in the same way.
An important observation is that the right hand side of  (\ref{phiODE.def}) does not involve $\phi^{\kappa_m}$, 
but only $\phi^0, \phi^1, \ldots, \phi^{\kappa_{m-1}}$.
These ODEs have a rigorous meaning.  So, we inductively 
define $\phi^{\kappa_m}$ as a unique solution of (\ref{phiODE2.def})--(\ref{phiODE.def}).

If the right hand side of(\ref{phiODE2.def})--(\ref{phiODE.def}) is denoted by $d Q^{\kappa_m}_t$, 
then $\phi^{\kappa_m}$ can be written explicitly as follows;
\begin{align}\label{phiODE5.def}
\phi^{\kappa_m}_T
= 
\tilde{J}(\gamma)_T \int_0^T   \tilde{J}(\gamma)_t^{-1}  d Q^{\kappa_m}_t,   %
\end{align}
where we set $\tilde{J}(\gamma) = J(\gamma, 0) =J (0 w +\gamma , 0^{1/H} \lambda) $.
See (\ref{jODE1.eq}) for the definition of $J$.

Define the remainder term $R^{\kappa_{m+1},  \ve}$ by 
\[
R^{\kappa_{m +1}, \ve}_t  =   \tilde{y}^{\ve}_t 
-   \bigl(  \phi^0_t   + \ve    \phi^1_t + \cdots +  \ve^{\kappa_m}    \phi^{\kappa_m}_t  \bigr).
 \]
We will estimate this remainder term in $L^q$-sense.
\begin{pr}\label{pr.LqestR}
For any $m \in {\bf N}$ and $q \in (1, \infty)$, 
$\|  \phi^{\kappa_m} \|_{\alpha -hld} \in L^q (\mu)$ and 
\[
{\mathbb E}  \bigl[    \|   R^{\kappa_{m +1}, \ve}  \|_{\alpha -hld}^q  \bigr]^{1/q}=O (\ve^{\kappa_{m +1}})  
\qquad
\quad
\mbox{as $\ve \searrow 0$.}
\]
\end{pr}

\Proof
From the expression (\ref{phiODE5.def}) and induction,  it is easy to see that 
$\|  \phi^{\kappa_m} \|_{\alpha -hld} \in \cap_{1 <q < \infty} L^q$ for any $m$.
Let us consider $R^{1, \ve}_t = \triangle \phi=  \tilde{y}^{\ve} - \phi^0 = I (\ve w +\gamma, \ve^{1/H} \lambda)-I (\gamma, {\bf 0})$.
Here, $I$ stands for the It\^o map and ${\bf 0}$ stands for one-dimensional constant path staying at $0$.

Define $\omega(s,t) = ( \|  w\|^p_{\alpha -hld} +  \|  \gamma\|^p_{\alpha -hld}  +1)  (t-s)$ with $\alpha =1/p$.
This control function satisfies
\begin{align}
&
| (\ve w_t  + \gamma_t)  - (\ve w_s  + \gamma_s)  | +  | \ve^{1/H} t -  \ve^{1/H}  s|  &\le \omega(s,t)^{1/p} 
\nn\\
& 
|  \{ (\ve w_t  + \gamma_t)  - (\ve w_s  + \gamma_s) \} - \{  \gamma_t  -  \gamma_s     \}  |
+
|\ve^{1/H} t -  \ve^{1/H}  s|  &\le    \ve \omega(s,t)^{1/p}
\nn
\end{align}
for all $0 \le s \le t \le 1$ and $\ve \in [0,1]$.
Hence, by the local Lipschitz continuity of It\^o map $I$, 
\[
| R^{1, \ve}_t -R^{1, \ve}_s | \le \ve  C (1+ \omega(0,1))^{(p-1)/p}  \exp ( C \omega(0,1) ) \omega(s,t)^{1/p}
\]
for some positive constant $C$.
Since $p<2$, we can use Fernique's theorem to obtain
the desired estimate holds when $\kappa_{m+1} =1$.

Before we prove the higher order cases, let us observe the concrete expression
 for several  $R^{\kappa_{m +1}, \ve}$'s.
In the sequel, we write $\kappa_{m +1} =: \kappa_{m } +$ for simplicity of notation.
First we consider $R^{1+, \ve} = R^{1/H, \ve}=   \tilde{y}^{\ve} -   \phi^0  -\ve    \phi^1$.
A straight forward computation yields;
\begin{align}
dR^{1+, \ve}_t
&=
\ve \{ \sigma(\tilde{y}^{\ve}_t) -  \sigma( \phi^0_t) \} dw_t 
\nn\\
&+
\Bigl[
\{       \sigma(\tilde{y}^{\ve}_t)   -  \sigma( \phi^0_t)  \}d\gamma_t  - \nabla \sigma (\phi^0_t ) \la \ve \phi^1_t, d\gamma_t  \ra  \Bigr]
+\ve^{1/H}    b(\tilde{y}^{\ve}_t)dt.
\label{Rem0.eq}
\end{align}
From this,  we immediately have
\begin{align}
dR^{1+, \ve}_t    &-   \nabla \sigma (\phi^0_t ) \la  R^{1+, \ve}_t   , d\gamma_t  \ra
=
\ve \{ \sigma(\tilde{y}^{\ve}_t) -  \sigma( \phi^0_t) \} dw_t 
\nn\\
&+\frac12 \int_0^1 d\theta  
 \nabla \sigma (\phi^0_t + \theta  R^{1, \ve}_t ) \la   R^{1, \ve}_t ,R^{1, \ve}_t , d\gamma_t  \ra  
 +\ve^{1/H}    b(\tilde{y}^{\ve}_t)dt 
 \quad (=: dL^{1+, \ve}_t).
\label{Rem1.eq}
\end{align}
Observe that, on the right hand side,  there are only $R^{1, \ve}, \tilde{y}^{\ve}, \phi^0, \gamma, w$, 
which are known quantities, 
but no $R^{1+, \ve}$.
Since $R^{1+, \ve}_T =   \tilde{J}(\gamma)_T \int_0^T   \tilde{J}(\gamma)_t^{-1}  d L^{1+, \ve}_t$ as before,
it suffices to show that $\| L^{1+, \ve} \|_{\alpha -hld} =O(\ve^{1/H})$ for any $L^q$.

Since  $\| \ve^{1/H}   \int_0^{\cdot} b(\tilde{y}^{\ve}_t)dt  \|_{\alpha -hld} \le  C \ve^{1/H} \| \tilde{y}^{\ve} \|_{\alpha -hld} $, 
the third term of $L^{1+, \ve}$ is $O(\ve^{1/H})$ in any $L^q$.
Similarly, 
$ \ve \|   \int_0^{\cdot}\{ \sigma(\tilde{y}^{\ve}_t) -  \sigma( \phi^0_t) \} dw_t  \|_{\alpha -hld} \le  C \ve \|  R^{1, \ve} \|_{\alpha -hld} \| w \|_{\alpha -hld} $,
the first 
term of $L^{1+, \ve}$ is $O(\ve^{2})$ in any $L^q$.
For any $\theta$, 
$\|   \nabla \sigma (\phi^0_{\cdot}+ \theta  R^{1, \ve}_{\cdot} )   \|_{\alpha -hld} \le  
C( \|  \phi^0   \|_{\alpha -hld} +  \|  R^{1, \ve}   \|_{\alpha -hld} )$.
Hence, we have
$$
\|  \int_0^{\cdot}   \int_0^1 d\theta  
 \nabla \sigma (\phi^0_t + \theta  R^{1, \ve}_t ) \la   R^{1, \ve}_t ,R^{1, \ve}_t , d\gamma_t  \ra   \|_{\alpha -hld} 
 \le 
  C( \|  \phi^0   \|_{\alpha -hld} +  \|  R^{1, \ve}   \|_{\alpha -hld} )\|  R^{1, \ve} \|_{\alpha -hld}^2.
  $$
We see from the above inequality that the second 
term of $L^{1+, \ve}$ is $O(\ve^{2})$ in any $L^q$ 
and hence $\| L^{1+, \ve} \|_{\alpha -hld} =O(\ve^{1/H})$ in any $L^q$.
Thus, we have obtained the desired estimate for $R^{1+, \ve}=R^{1/H, \ve}$.

The estimate for  $R^{(1/H)+, \ve} =   \tilde{y}^{\ve} -   \phi^0  -\ve    \phi^1-\ve^{1/H}    \phi^{1/H}$
can easily be obtained as follows.
We can immediately see from (\ref{phiODE4.def}) and (\ref{Rem1.eq}) that
\begin{align}
dR^{(1/H)+, \ve}_t    &-   \nabla \sigma (\phi^0_t ) \la  R^{(1/H)+, \ve}_t   , d\gamma_t  \ra
=
\ve \{ \sigma(\tilde{y}^{\ve}_t) -  \sigma( \phi^0_t) \} dw_t 
\nn\\
& +\frac12 \int_0^1 d\theta  
 \nabla \sigma (\phi^0_t + \theta  R^{1, \ve}_t ) \la   R^{1, \ve}_t ,R^{1, \ve}_t , d\gamma_t  \ra  
\nn\\
& +\ve^{1/H} \{   b(\tilde{y}^{\ve}_t)    - b(\phi^0_t) \}dt
 \qquad (=: dL^{(1/H)+, \ve}_t).
\label{Rem1/H.eq}
\end{align}
Notice that we have essentially shown that $\| L^{(1/H)+, \ve} \|_{\alpha -hld} =O(\ve^{2})$ in any $L^q$.
Thus, we have obtained the desired estimate for $R^{(1/H)+, \ve}=R^{2, \ve}$.

Next,  we will estimate 
$R^{2+, \ve} =   \tilde{y}^{\ve} -   \phi^0  -\ve    \phi^1-\ve^{1/H}    \phi^{1/H} -  \ve^2 \phi^2$.
From (\ref{phiODE3.def}), (\ref{phiODE4.def}), 
and (\ref{Rem0.eq}),  we see that
\begin{align}
dR^{2+, \ve}_t
&=
\bigl[ 
 \{ \sigma(\tilde{y}^{\ve}_t) -  \sigma( \phi^0_t) \} \ve dw_t  
 -  \nabla \sigma (\phi^0_t ) \la  \ve \phi^1_t, \ve dw_t  \ra
\bigr]
\nn\\
&+
\bigl[
\{       \sigma(\tilde{y}^{\ve}_t)   -  \sigma( \phi^0_t)  \}d\gamma_t 
 - \nabla \sigma (\phi^0_t ) \la \ve \phi^1_t  +\ve^{1/H}    \phi^{1/H} + \ve^2 \phi^2  , d\gamma_t  \ra  \bigr]
\nn\\
&
-\frac12  \nabla^2 \sigma (\phi^0_t ) \la \ve \phi^1_t ,   \ve \phi^1_t ,  d\gamma_t\ra
+\ve^{1/H}    \{   b(\tilde{y}^{\ve}_t)    - b(\phi^0_t) \}dt.
\label{Rem2.eq}
\end{align}
The second term on the right hand side is equal to
\begin{align}
 \nabla \sigma (\phi^0_t ) \la   dR^{2+, \ve}_t,   d\gamma_t  \ra
  &+
     \frac12  \nabla^2 \sigma (\phi^0_t ) \la  R^{1, \ve}_t ,    R^{1,\ve}_t ,  d\gamma_t\ra 
     \nn\\
     &+
      \int_0^1 \frac{(1-\theta)^2d\theta}{2!} \nabla^3 \sigma (\phi^0_t  + \theta R^{1, \ve}_t ) \la  R^{1, \ve}_t ,R^{1, \ve}_t ,   R^{1,\ve}_t ,  d\gamma_t\ra.     
      \nn
        \end{align}
Hence, (\ref{Rem2.eq}) is equivalent to the following;
\begin{align}
dR^{2+, \ve}_t -  \nabla \sigma (\phi^0_t ) \la   dR^{2+, \ve}_t,   d\gamma_t  \ra
&=
\bigl[ 
 \{ \sigma(\tilde{y}^{\ve}_t) -  \sigma( \phi^0_t) \} \ve dw_t  
 -  \nabla \sigma (\phi^0_t ) \la  \ve \phi^1_t, \ve dw_t  \ra
\bigr]
\nn\\
&+ 
 \frac12  
 \bigl[ 
  \nabla^2 \sigma (\phi^0_t ) \la  R^{1, \ve}_t ,    R^{1,\ve}_t ,  d\gamma_t\ra 
   -\nabla^2 \sigma (\phi^0_t ) \la \ve \phi^1_t ,   \ve \phi^1_t ,  d\gamma_t\ra  
   \bigr]
   \nn\\
   &+
      \int_0^1 \frac{(1-\theta)^2d\theta}{2!} \nabla^3 \sigma (\phi^0_t  + \theta R^{1, \ve}_t ) \la  R^{1, \ve}_t ,R^{1, \ve}_t ,   R^{1,\ve}_t ,  d\gamma_t\ra     
      \nn\\
      & +\ve^{1/H}    \{   b(\tilde{y}^{\ve}_t)    - b(\phi^0_t) \}dt   \qquad (=: dL^{2+, \ve}_t).
      \label{Rem3.eq}
\end{align}
  Then, $R^{2+, \ve}_T =   \tilde{J}(\gamma)_T \int_0^T   \tilde{J}(\gamma)_t^{-1}  d L^{2+, \ve}_t$.

Let us observe the right hand side of (\ref{Rem3.eq}). 
There are no $R^{2+, \ve}$ or $\phi^2$.
By the assumption of induction, we may only use the relation
$R^{2, \ve}=R^{(1/H)+, \ve}  =   \tilde{y}^{\ve} -   \phi^0  -\ve    \phi^1-\ve^{1/H}    \phi^{1/H}$ and the estimates of 
$R^{\kappa, \ve}$ for $\kappa=1,1/H, 2$ (and of $\phi^{\kappa}$'s).
 In the same way as above,  by using the Taylor expansion, 
 we can prove that  $\| L^{2+, \ve} \|_{\alpha -hld} =O(\ve^{1 +(1/H ) } )$ in any $L^q$.    
Cancellation of the terms of order $\le 2$ on the right hand side is no mystery 
because of the way $\phi^{\kappa}$'s are defined.
    Thus, we have obtained the desired estimate for $R^{2+, \ve}=R^{1 +(1/H), \ve}$.

Higher order remainder terms can be dealt with in a similar way.
We give a sketch of proof.
There exists 
\[
L^{\kappa_{m +1}, \ve}_t 
= L^{\kappa_{m +1}, \ve} [ \phi^0, \ldots, \phi^{\kappa_{m-1}} ;  R^{1,\ve}_t ,  \ldots,   R^{\kappa_m , \ve} ; w, \gamma]_t 
\]
such that
$ dR^{\kappa_{m +1}, \ve}_t -  \nabla \sigma (\phi^0_t ) \la   dR^{\kappa_{m +1}, \ve}_t,   d\gamma_t  \ra = d L^{\kappa_{m +1}, \ve}_t $.
Due to cancellation
 $\| L^{\kappa_{m +1}, \ve} \|_{\alpha -hld} =O(\ve^{\kappa_{m +1}} )$ holds in any $L^q$.
 This proves the assertion.
 \QED


The next proposition shows that, when evaluated at $t=1$,
Eq. (\ref{asyphi1.eq}) gives an asymptotic expansion 
in ${\bf D}_{\infty} ({\bf R}^n)$.

\begin{pr}\label{pr.tay.Dinf}
We have the following asymptotic expansion in ${\bf D}_{\infty} ({\bf R}^n)$.
\begin{equation}\label{8asy.eq}
\tilde{y}^{\ve}_1  
\sim \phi^0_1 + \ve^{\kappa_1} \phi^{\kappa_1}_1  + \ve^{\kappa_2} \phi^{\kappa_2}_1 
+ \cdots   
\qquad\quad
\mbox{as $\ve \searrow 0$. }
\end{equation}
Here, $0=\kappa_0 < \kappa_1 < \kappa_2 < \cdots$ are all the  elements of $\Lambda_1={\bf N} +\frac{1}{H} {\bf N}$ in increasing order.
\end{pr}

\Proof
By using induction and basic properties of Young integral,  we can easily see that
$\phi^{\kappa_m}_1$ is in $[\kappa_m]$-th inhomogeneous Wiener chaos 
for each $t$ and $m$. 
In particular, $\phi^{\kappa_m}_1 \in {\bf D}_{\infty}$. 
If $k \ge [\kappa_m] +1$, then $D^k R^{\kappa_{m+1}, \ve }_1 = D^k \tilde{y}^{\ve}_1$.
From Proposition \ref{pr.est.Dy}, this is $O (\ve^{k})$, and hence $O (\ve^{ \kappa_{m+1}})$
in any $L^q$.
A stronger version of  Meyer's equivalence (e.g., Theorem 4.6, \cite{shbk})
implies that $R^{\kappa_{m+1}, \ve }_1 $ is $O (\ve^{ \kappa_{m+1}})$ 
in ${\bf D}_{q,k}$  for any $q$ and sufficiently large $k$.
Since ${\bf D}_{q,k}$-norm is increasing in $k$, 
the proof is completed.
\QED


We now recall the following Taylor expansion of It\^o map around $\gamma$
in the deterministic sense.
\begin{lm}\label{lm.tay.deter}
{\rm (i)}~For each $m$, there exists $c=c(\kappa_m)$ such that
\[
\| \phi^{\kappa_m}\|_{\alpha -hld}  \le c (1+ \| w \|_{\alpha -hld} )^{\kappa_m }
\qquad
\mbox{ for all $w \in C_0^{\alpha -hld} ([0,1], {\bf R}^d)$.}
\]
\noindent
{\rm (ii)}~For each $m$ and $r>0$, 
there exists $c'=c' (\kappa_m, r)$ such that
\[
\| R^{\kappa_{m+1}, \ve }\|_{\alpha -hld}  \le c' (\ve + \| \ve w \|_{\alpha -hld} )^{\kappa_{m+1} },
\qquad
\mbox{ if $\| \ve w \|_{\alpha -hld} \le r$.}
\]
\end{lm}

\Proof
This is immediate since $\tilde{y}^{\ve} =I (\ve w +\gamma, \ve^{1/H} \lambda)$
and It\^o map $I$ is Fr\'echet smooth by Li and Lyons's result \cite{ll}.
It is also possible to prove this lemma by using the explicit expression of $R^{\kappa_{m+1}, \ve }$ 
and mathematical induction
as in the proof of Proposition \ref{pr.LqestR} above.
\QED


\section{Off-diagonal short time asymptotics}

In this section we prove the short time asymptotics of kernel function
$p_t(a,a')$ when $a \neq a'$.
We basically follow Watanabe \cite{wa}.
In this paper, however,  we can localize around the energy minimizing path in the abstract Wiener space
since It\^o map is continuous in our setting.
This makes the proof slightly simpler.

\subsection{Localization around energy minimizing path}

For $\gamma \in {\cal H}$,
let $\phi^0= \phi^0(\gamma)$ be a unique solution of  (\ref{phi0.def}), which starts at $a \in {\bf R}^n$.
Set,  for $a \neq a'$, 
\[
K_a^{a'} = \{ \gamma \in {\cal H} ~|~  \phi^0_1(\gamma) =a'\}.
\]
We only consider the case that $K_a^{a'} $ is not empty.
For example, if {\bf (A1)} is satisfied for any $a$, then $K_a^{a'} $ is not empty for any $a'$.
From the Schilder-type large deviation theory, it is easy to see that
$\inf\{ \|\gamma\|_{\cal H} ~|~ \gamma \in  K_a^{a'}\} 
= \min\{ \|\gamma\|_{\cal H} ~|~ \gamma \in  K_a^{a'}\}$.

We continue to assume {\bf (A1)}.
Now we introduce another assumption;
\vspace{3mm}
\\
{\bf (A2):} $\bar{\gamma} \in K_a^{a'}$ which minimizes ${\cal H}$-norm exists uniquely.
\vspace{3mm}
\\
In the sequel, $\bar{\gamma}$ denotes the minimizer in Assumption {\bf (A2)}
and we use the results of the previous section for this $\bar{\gamma}$.

Note that {\rm (i)}~ the mapping $\gamma \in {\cal H} \hookrightarrow {\cal W} \mapsto  \phi_1^0 (\gamma) \in {\bf R}^n$
is Fr\'echet differentiable
and 
{\rm (ii)}~ its Jacobian is a surjective linear mapping from ${\cal H}$ to ${\bf R}^n$ for any $\gamma$,
because there exists a positive constant $c=c(\gamma)$ such that 
\begin{equation}\label{detuki.eq}
\Bigl(  \la D  \phi_1^{0,i } (\gamma),   D  \phi_1^{0, j} (\gamma) \ra_{{\cal H}^*}  \Bigr)_{ 1 \le i,j \le n}  \ge  c \cdot {\rm Id}_n.
\end{equation}
This can be shown in the same way as in the proof of Proposition \ref{pr.unif.nondeg}.
(Actually, it is easier since $\gamma$ is non-random and fixed here.)

Therefore, by the Lagrange multiplier method,  there exists 
$\bar{\nu} =(\bar{\nu}_1, \ldots, \bar{\nu}_n) \in {\bf R}^n$ uniquely such that
the map
\begin{equation}\label{Lmul1.eq}
{\cal H} \times {\bf R}^n \ni (\gamma, \nu) \mapsto   
\frac12   \|\gamma\|_{\cal H}^2 - \la  \nu, \phi^0_1(\gamma) -a' \ra_{{\bf R}^n} \in {\bf R}
\end{equation}
attains extremum at $(\bar{\gamma}, \bar{\nu})$.
By differentiating in the direction of $k \in {\cal H}$, 
we have
\begin{equation}\label{Lmul2.eq}
\la \bar{\gamma} , k \ra_{\cal H} = \la  \bar{\nu},  D_k \phi^0_1( \bar{\gamma} ) \ra_{{\bf R}^n}
=
\bigl\la  \bar{\nu},  
\tilde{J}(\bar{\gamma})_1 \int_0^1   \tilde{J}(\bar{\gamma})_t^{-1}     \sigma(   \phi^0_t(\bar{\gamma} ) )  dk_t
 \bigr\ra_{{\bf R}^n}.
\end{equation}
Here, the definition of $\tilde{J}(\bar{\gamma})$ was given just below (\ref{phiODE5.def}) 
and the integral on the right hand side is  Young integral.
Hence, $\la \bar{\gamma} , \,\cdot\, \ra_{\cal H}$ extends to a continuous linear 
functional on ${\cal W}$.


Let us introduce Besov-type norms.
In the context of Malliavin calculus, these norms are often more useful than
 H\"older norms and $p$-variation norms 
 since (a power of) these norms become ${\bf D}_{\infty}$-functionals.
For $m >0$, $0 < \theta <1$, and $x \in C_0 ([0,1], {\bf R}^d)$,
we set 
\[
\| x \|_{m,\theta -B} :=  \Bigl(  \iint_{0 \le s \le t \le 1}  
\frac{ | x_t -x_s |^m }{ |t-s|^{2 +m\theta} }dsdt \Bigr)^{1/m}.
\]
and 
$
 C^{m,\theta -B}_0 ([0,1], {\bf R}^d) 
= \{  x \in C_0 ([0,1], {\bf R}^d)   ~|~  \| x \|_{m,\theta -B} < \infty \}.
$
It is known that $\| x\|_{\theta -hld} \le c \| x \|_{m,\theta -B}$
for some constant $c= c_{m, \theta} >0$.
Hence, $ C^{m,\theta -B}_0 ([0,1], {\bf R}^d) \subset C^{\theta -hld}_0 ([0,1], {\bf R}^d)$.

Let $(w_t)$ be fBm with Hurst parameter $H \in (1/2,1)$
and let $\alpha (=1/p) <H$ as before.
Since ${\mathbb E}[ | w_t -w_s|^2] =d |t-s|^{2H}$,
we can easily see 
${\mathbb E}[  \| x \|_{m,\alpha -B}^m ]< \infty$ if  $m >1/(H - \alpha)$.
Therefore, the law of fBm, $\mu=\mu^H$, is supported in 
$C^{m,\alpha -B}_0 ([0,1], {\bf R}^d)$ if  $m >1/(H - \alpha)$.
Set ${\cal W}_B$ to be the closure of  Cameron-Martin space ${\cal H}={\cal H}^H$
in $C^{m,\alpha -B}_0 ([0,1], {\bf R}^d)$.
Then, $( {\cal W}_B, {\cal H}, \mu)$ is also an abstract Wiener space.
%
%


Now we recall Schilder-type large deviation principle for 
scaled Gaussian measures.
For $\ve >0$, let $\mu_{\ve}$ be the law of the law of the process 
$(\ve w_t)_{0 \le t \le 1}$. 
This is a measure on ${\cal W}_B$.
Set ${\cal I}(w) = \|w\|_{{\cal H}}^2/2 ~(\mbox{if $w \in {\cal H}$})$
and ${\cal I} (w) =\infty~ (\mbox{otherwise})$.
It is well-known that  ${\cal I}: {\cal W}_B \to [0,\infty]$ is lower semicontinuous
and that ${\cal I}$ is good, i.e.,
the level set $\{ w ~|~ {\cal I}(w) \le r \}$ is compact in ${\cal W}_B$ for any $r \in [0, \infty)$.

The family $\{\mu_{\ve} \}_{\ve >0}$ satisfies large deviation principle 
as $\ve \searrow 0$ with a good rate function ${\cal I}$, that is,  for any 
measurable set $A \subset {\cal W}_B$
\begin{align}
-\inf_{w \in A^{\circ}}  {\cal I}(w)  &\le \liminf_{ \ve \searrow 0 }  \ve^2 \log \mu_{\ve}(A^{\circ})  
\le 
 \limsup_{ \ve \searrow 0 }  \ve^2 \log \mu_{\ve}(\bar{A}) 
\le  
-\inf_{w \in \bar{A}} {\cal I}(w).
\label{ldp1.ineq}
\end{align}

Next, set $\hat\mu_{\ve} = \mu_{\ve} \otimes \delta_{ \ve^{1/H} \lambda }$,
where $\lambda$ is a one-dimensional path defined by $\lambda_t =t$
and $\otimes$ stands for the product of probability measures. 
In other words, $\hat\mu_{\ve}$ is the law of the $(d+1)$-dimensional  process 
$(\ve w_t , \ve^{1/H} t )_{0 \le t \le 1}$ under $\mu$.
This measure is supported on 
${\cal W}_B \oplus {\bf R} \la \lambda \ra \subset C_0^{m,\alpha -B}( [0,1]; {\bf R}^{d+1})$.
Define $\hat{ {\cal I}} (w ; l) = \|w\|_{{\cal H}}^2/2 ~(\mbox{if $w \in {\cal H}$ and $l_t \equiv 0$})$
and $\hat{ {\cal I}} (w, l) =\infty~ (\mbox{otherwise})$.
Here, $l$ is a one-dimensional path.

From (\ref{ldp1.ineq}) we can easily show that
 $\{\hat\mu_{\ve} \}_{\ve >0}$ satisfies large deviation principle 
as $\ve \searrow 0$ with a good rate function $\hat{ {\cal I}}$,
that is, 
for any measurable set $A \subset {\cal W}_B \oplus {\bf R} \la \lambda \ra$,
\begin{align}
-\inf_{w \in A^{\circ}} \hat{ {\cal I} }(w)  
\le \liminf_{ \ve \searrow 0 }  \ve^2 \log \hat\mu_{\ve}(A^{\circ})  
\le 
 \limsup_{ \ve \searrow 0 }  \ve^2 \log \hat\mu_{\ve}(\bar{A}) 
\le  
-\inf_{w \in \bar{A}} \hat{ {\cal I}}(w).
\label{ldp2.ineq}
\end{align}
We will use  (\ref{ldp2.ineq})  in Lemma \ref{lm.ldpcut} below to show that
only a neighborhood of 
the minimizer $\bar{\gamma}$ contributes to the asymptotic expansion.

From now on, 
we will fix an even integer $m>0$ such that $m >1/(H - \alpha)$.
Then, it is easy to check $\| w\|_{m,\alpha -B}^m \in {\bf D}_{\infty}$. 
In fact, this functional is an element of $m$th inhomogeneous Wiener chaos,
i.e., $D^{m+1} \| w\|_{m,\alpha -B}^m =0$.



Now we introduce a cut-off function. 
Let $\psi : {\bf R} \to [0,1]$ be a smooth function such that 
$\psi (u) =1$ if $|u| \le 1/2$ and $\psi (u) =0$ if $|u| \ge 1$.
For each  $\eta >0$ and $\ve >0$, we set 
\[
\chi_{\eta} (\ve, w) 
= \psi \Bigl(  
\frac{1}{\eta^m}     \| \ve w - \bar{\gamma}\|_{m,\alpha -B}^m
\Bigr).
\]
We can easily see that
$\chi_{\eta} (\ve, \,\cdot\,) \in {\bf D}_{\infty}$.
Shifting by $\bar{\gamma}/\ve$, we have
$$\chi_{\eta} (\ve, w +\frac{\bar{\gamma}}{\ve}) 
=\psi \Bigl(  
\frac{\ve^m}{\eta^m}     \|  w \|_{m,\alpha -B}^m
\Bigr).
$$
It is easy to see from Taylor expansion for $\psi$ 
that, for any $\eta >0$ and any $M \in {\bf N}$, 
the following asymptotics holds;
\begin{equation}\label{chi_asy.eq}
\chi_{\eta} (\ve, w +\frac{\bar{\gamma}}{\ve}) = 1+O (\ve^M)
\qquad
\mbox{in ${\bf D}_{\infty}$ as $\ve \searrow 0$.}
\end{equation}


%
%
The following lemma states that only the paths sufficiently close to
 the minimizer $\bar{\gamma}$ contribute to the asymptotics.
\begin{lm}\label{lm.ldpcut}
Assume {\bf (A1)} and {\bf (A2)}.
Then, for any $\eta >0$, there exists $c=c_{\eta}>0$ such that 
\[
0 \le {\mathbb E} \bigl[ (1 - \chi_{\eta} (\ve, w) ) \cdot  \delta_{a'} (y^{\ve}_1 )
\bigr]
=
O \Bigl(
\exp \bigl\{  - \frac{   \|  \bar{\gamma}\|_{\cal H}^2 +c }{2 \ve^2}  \bigr\}
\Bigr)
\qquad
\mbox{as $\ve \searrow 0$.}
\]
\end{lm}

\Proof
We take $\eta' >0$ arbitrarily and we will fix it for a while.  
It is obvious that
\begin{align}
0 \le
 {\mathbb E} \bigl[ (1 - \chi_{\eta} (\ve, w) ) \cdot  \delta_{a'} (y^{\ve}_1 )
\bigr]
=
 {\mathbb E} 
\Bigl[ (1 - \chi_{\eta} (\ve, w) ) 
\psi \Bigl(   
\frac{ | y^{\ve}_1 -  a' |^2  }{\eta'^2 }
\Bigr)
 \cdot  \delta_{a'} (y^{\ve}_1 )
\Bigr].
\label{ldp_pf1}
\end{align}

Set $g(u) = u \vee 0$ for $u \in {\bf R}$.
Then, in the sense of distributional derivative, $g^{\prime \prime} (u) =\delta_0$.
Take a bounded continuous function $C: {\bf R}^n \to {\bf R}$ such that  
$C(u_1, \ldots, u_n) = g(u_1- a'_1) g(u_2- a'_2)\cdots g(u_n- a'_n)$ 
if $|u -a'| \le 2\eta'$.
Then, the right hand side of (\ref{ldp_pf1}) is equal to 
\begin{align}
 {\mathbb E} 
\Bigl[ (1 - \psi ) \bigl(\frac{1}{\eta^m}     \| \ve w - \bar{\gamma}\|_{m,\alpha -B}^m \bigr) 
\cdot
\psi \Bigl(   
\frac{ | y^{\ve}_1 -  a' |^2  }{\eta'^2 }
\Bigr)
 \cdot   ( \partial_1^2 \cdots \partial_n^2 C)    (y^{\ve}_1 )
\Bigr].
\label{ldp_pf2}
\end{align}

Now, we use integration by parts for (generalized) Wiener functionals
as in pp. 6--7, \cite{wa}
 to see that 
(\ref{ldp_pf2}) is equal to a finite sum of the following form;
\begin{align}
\sum_{j,k }  {\mathbb E} 
\Bigl[ 
F_{j,k} (\ve, w) \cdot (1 - \psi )^{(j)} 
\bigl(\frac{1}{\eta^m}     \| \ve w - \bar{\gamma}\|_{m,\alpha -B}^m \bigr) 
\cdot
\psi^{(k)} \Bigl(   
\frac{ | y^{\ve}_1 -  a' |^2  }{\eta'^2 }
\Bigr)
 \cdot  
C (y^{\ve}_1 )
\Bigr].
\label{ldp_pf2.5}
\end{align}
Here, $F_{j,k} (\ve, w)$ is a polynomial in components of 
{\rm (i)}~ $y^{\ve}_1$ and its derivatives, {\rm (ii)}~ $ \| \ve w - \bar{\gamma}\|_{m,\alpha -B}^m$
and its derivatives, 
{\rm (iii)}~$\tau (\ve)$, which is Malliavin covariance matrix of $y^{\ve}_1$,  and its derivatives,
and {\rm (iv)}~ $\kappa (\ve):= \tau (\ve)^{-1}$.
Note that the derivatives of $\kappa (\ve)$ do not appear.

From  Proposition \ref{pr.unif.nondeg},
there exists $r'>0$ such that 
 $|\kappa^{ij} (\ve)| = O(\ve^{-r'})$ in $L^q$ 
as $\ve \searrow 0$ for all $1<q<\infty$.
(Recall a well-known formula to obtain the inverse matrix $A^{-1}$ with the adjugate matrix of $A$
divided by $\det A$.)
Therefore, there exists $r>0$ such that 
$|F_{j,k} (\ve)| = O(\ve^{-r})$ in $L^q$ 
as $\ve \searrow 0$ for all $1<q<\infty$.

By H\"older's inequality, (\ref{ldp_pf2.5}) is dominated by 
\begin{align}
& \quad
\frac{c}{ \ve^r}  \sum_{j,k } 
{\mathbb E} 
\Bigl[ 
\bigl| (1 - \psi )^{(j)} \bigl(\frac{1}{\eta^m}     \| \ve w - \bar{\gamma}\|_{m,\alpha -B}^m \bigr) \bigr|^{q'}
\,
\Bigl|
\psi^{(k)} \Bigl(   
\frac{ | y^{\ve}_1 -  a' |^2  }{\eta'^2 }
\Bigr)\Bigr|^{q'}
\Bigr]^{1/q'}
\nn\\
&
\le
\frac{c}{ \ve^r}
\mu \Bigl[ 
  \| \ve w - \bar{\gamma}\|_{m,\alpha -B}^m \ge \frac{\eta^m} {2},
\quad
| y^{\ve}_1 -  a'  | \le \eta'
\Bigr]^{1/q'}.
\label{ldp_pf3}
\end{align}
Here, $1/q +1/q' =1$ and $c=c(q,q', \eta, \eta')$ is a positive constant,
which may change from line to line.

Since we may let $q' \searrow 1$ after taking $\limsup$, we obtain the following; 
\begin{align}
& \quad
\limsup_{\ve \searrow 0}  
\ve^2 \log  {\mathbb E} \bigl[ (1 - \chi_{\eta} (\ve, w) ) \cdot  \delta_{a'} (y^{\ve}_1 )
\bigr]
\nn\\
&\le
\limsup_{\ve \searrow 0}  
\ve^2 \log 
\mu \Bigl[ 
  \| \ve w - \bar{\gamma}\|_{m,\alpha -B}^m \ge \frac{\eta^m}{2},
\quad
| y^{\ve}_1 -  a'  | \le \eta'
\Bigr]
\nn\\
&=
\limsup_{\ve \searrow 0}  
\ve^2 \log 
\hat\mu^{\ve} \Bigl[ \bigl\{   ( w, l ) \in {\cal W}_B \oplus {\bf R}\la \lambda\ra ~|~
  \|  w - \bar{\gamma}\|_{m,\alpha -B}^m \ge \frac{\eta^m}{2},
\,
| I(w,l)_1 -  a'  | \le \eta'
\bigr\}
\Bigr]
\nn\\
&\le
- \inf \Bigl\{ \frac{ \|\gamma \|^2_{{\cal H}} }{2}~|~   
\|  \gamma - \bar{\gamma}\|_{m,\alpha -B}^m \ge \frac{\eta^m}{2},
\quad
| \phi^0[\gamma]_1 -  a'  | \le \eta'   \Bigr\}.
\label{ldp_pf4}
\end{align}
Here, $I$ denotes the It\^o map corresponding to ODE (\ref{ygSDE.eq})
and we have used the large deviation for the last inequality.
(Note that  continuity of It\^o map is used.)
Recall that
$\phi^0[\gamma]=  I(\gamma, {\bf 0})$ is given by ODE (\ref{phi0.def}).

Now let $\eta'$ tend to $0$.
As $\eta'$ decreases, the right hand side of (\ref{ldp_pf4}) decreases.
The proof is finished if the limit is strictly smaller than $- \|\bar{\gamma} \|^2_{{\cal H}}  /2$.
Assume otherwise.
Then, there exists $\{\gamma_k\}_{k=1}^{\infty} \subset {\cal H}$ such that
\[
\|  \gamma_k - \bar{\gamma}\|_{m,\alpha -B}^m \ge \frac{\eta^m}{2}, \quad 
| \phi^0[\gamma_k]_1 -  a'  | \le \frac{1}{k},  \mbox{ and,  }
\liminf_{k \to \infty} ( - \frac{\|\gamma_k \|^2_{{\cal H}}  }{2}) 
\ge - \frac{\|\bar{\gamma} \|^2_{{\cal H}}  }{2}.
\]
In particular, 
$\{\gamma_k\}$ is bounded in ${\cal H}$ and, hence, precompact in ${\cal W}_B$.  
Let $\gamma_{\infty}$ be any limit point.
For simplicity, 
a subsequence that converges to $\gamma_{\infty}$ is again denoted by $\{\gamma_k\}$.
Since $\gamma \mapsto \phi^0[\gamma]_1$ is continuous with respect to the topology 
of ${\cal W}_B$, we see that $\phi^0[\gamma_{\infty}]_1 = a' $ holds.
Also, we have $\|  \gamma_{\infty} - \bar{\gamma}\|_{m,\alpha -B}^m \ge \eta^m/2$.
So, $\gamma_{\infty} \neq \bar{\gamma}$.
From the lower semicontinuity of the rate function, 
we see that $\gamma_{\infty} \in {\cal H}$ and 
$\|\gamma_{\infty} \|^2_{{\cal H}}  /2 \le \| \bar{\gamma} \|^2_{{\cal H}}  /2$.
This clearly contradicts Assumption  {\bf (A2)}.
\QED


\subsection{Integrability lemmas}\label{subsec.int}

In this subsection, we prove a few lemmas for integrability of Wiener functionals
of exponential type which will be used in the short time asymptotic expansion.

Throughout this subsection we assume {\bf (A2)}.
Let $\bar{\gamma}$ be as in {\bf (A2)} and let $\phi^{\kappa_j}$ and 
$R^{\kappa_j +, \ve} =R^{\kappa_{j +1}, \ve}$ ($j=0,1,2,\ldots$)
be as in Section \ref{sec.taylor}
 with $\gamma = \bar{\gamma}$. 
First we consider 
\begin{align}
\frac{ R^{2+, \ve} }{\ve^2}
&=  
\frac{1}{\ve^2} (\tilde{y}^{\ve} -   \phi^0  -\ve    \phi^1-\ve^{1/H}    \phi^{1/H} -  \ve^2 \phi^2)
\nn\\
&= \ve^{ \kappa_4 -2 } \phi^{ \kappa_4} + \ve^{ \kappa_5 -2} \phi^{ \kappa_5} +\cdots.
\nn
\end{align}
Here, $ \kappa_4=1+(1/H)$ and $ \kappa_5= 3 \wedge (2/H)$.

\begin{lm}\label{lm.integ1}
Assume {\bf (A2)}.
For any $M>0$, there exists $\eta >0$ such that 
\[
\sup_{0 < \ve \le 1} 
{\mathbb E} \bigl[
\exp \bigl(   
M \la \bar\nu,  R^{2+, \ve}_1 \ra /\ve^2
\bigr)
I_{ \{  \|\ve w\|_{m,\alpha -B} \le \eta\} }
\bigr]  <\infty.
\] 
\end{lm}

\Proof
By Lemma \ref{lm.tay.deter}, if $\|\ve w\|_{ \alpha -hld} \le 1$, then 
there exists a constant $c_1, c_2 >0$ such that
$$
\| R^{2+, \ve} \|_{\alpha -hld} \le c_1 (\ve +  \|\ve w\|_{ \alpha -hld} )^{1+(1/H) }
\le 
c_2 (\ve +  \|\ve w\|_{m,\alpha -B})^{1+(1/H) }.
$$
Hence, if  $ \|\ve w\|_{m,\alpha -B} \le \eta \le 1$, then 
$$
\| R^{2+, \ve} \|_{\alpha -hld} /\ve^2 
\le 
c_2 (1 +  \| w\|_{m,\alpha -B})^{2 }(\ve + \eta)^{(1/H) -1} .
$$
Recall that, by Fernique's theorem, there exists a positive constant $\beta>0$ such that 
$
{\mathbb E} [ \exp (\beta  (1 +  \| w\|_{m,\alpha -B})^{2 } )] <\infty.
$
Take $0< \eta \le 1$ so that $ M |\bar\nu| c_2 (2\eta)^{(1/H) -1} \le \beta$.
Then, we see that 
\[
\sup_{0 < \ve \le \eta} 
{\mathbb E} \bigl[
\exp \bigl(   
M \la \bar\nu,  R^{2+, \ve}_1 \ra /\ve^2
\bigr)
I_{ \{  \|\ve w\|_{m,\alpha -B} \le \eta\} }
\bigr]  <\infty.
\] 
Note that, 
if $\|\ve w\|_{m,\alpha -B} \le \eta$ and $\eta \le \ve \le 1$, 
then $\| R^{2+, \ve} \|_{\alpha -hld} /\ve^2$ is bounded.
This completes the proof.
\QED

Next we consider 
\begin{align}
\frac{ R^{1+, \ve} }{\ve}
&=  
\frac{1}{\ve} (\tilde{y}^{\ve} -   \phi^0  -\ve    \phi^1)
= \ve^{(1/H )-1}    \phi^{1/H} +  \ve^1 \phi^2
+\cdots.
\nn
\end{align}

\begin{lm}\label{lm.integ2}
Assume {\bf (A2)}.
For any $M>0$, there exists $\eta >0$ such that 
\[
\sup_{0 < \ve \le 1} 
{\mathbb E} \bigl[
\exp \bigl(   
M  \|  R^{1+, \ve} \|_{\alpha -hld}^2 /\ve^2
\bigr)
I_{ \{  \|\ve w\|_{m,\alpha -B} \le \eta \} }
\bigr]  <\infty.
\] 
\end{lm}

\Proof
By Lemma \ref{lm.tay.deter}, if $\|\ve w\|_{ \alpha -hld} \le 1$, then 
there exists a constant $c_1 >0$ such that
$$
\| R^{1+, \ve} \|_{\alpha -hld} \le c_1 (\ve +  \|\ve w\|_{ \alpha -hld} )^{1/H }
\le 
c_2 (\ve +  \|\ve w\|_{m,\alpha -B})^{1/H }.
$$
Hence, if  $\|\ve w\|_{m,\alpha -B} \le \eta \le 1$, then 
$$
\| R^{2+, \ve} \|^2_{\alpha -hld} /\ve^2 
\le 
c_2 (1 +  \| w\|_{m,\alpha -B})^{2 }(\ve + \eta )^{(2/H) -2} .
$$
Then, we can prove the lemma in the same way as in Lemma \ref{lm.integ1}.
\QED


From now on we assume {\bf (A1)} and {\bf (A2)}.
In addition, we introduce the following assumption;
\vspace{3mm}
\\
{\bf (A3)':}  \qquad\qquad
${\mathbb E}[  \exp (  \la\bar\nu ,  \phi^2_1 \ra )  ~|~ \phi^1_1  =0 ]
< \infty.$
\vspace{3mm}
\\

For all $1 \le j \le n$, $\phi^{1, j}_1 \in {\cal W}_B^* \subset {\cal H}^*$.
When we regard $\phi^{1, j}_1$ as an element of ${\cal H}$ by Riesz isometry, 
we write ${}^{\sharp} \phi^{1, j}_1 \in {\cal H} \subset {\cal W}_B$.
We have an orthogonal decomposition
${\cal H} =  
 \ker \phi^1_1  \oplus  ( \ker \phi^1_1 )^{\bot}$.
We denote by $\pi$ the orthogonal projection from ${\cal H}$ onto $\ker \phi^1_1$.
Note that
$ ( \ker \phi^1_1 )^{\bot}$ is an $n$-dimensional linear subspace spanned by 
$\{ {}^{\sharp} \phi^{1, 1}_1,  \ldots,  {}^{\sharp} \phi^{1, n}_1\}$.
Since $\dim  ( \ker \phi^1_1 )^{\bot} <\infty$,  the abstract Wiener space splits into two;
$ {\cal W}_B =    \overline{\ker \phi^1_1}^{ \| \,\cdot\,\|_{m,\alpha -B} }  \oplus  ( \ker \phi^1_1 )^{\bot}$.
The projection $\pi$ naturally extends to the one from 
${\cal W}_B$ onto $\overline{\ker \phi^1_1}^{ \| \,\cdot\,\|_{m,\alpha -B} } $,
which is again denoted by the same symbol.
There exist Gaussian measures $\mu_1$ and $\mu_2$ such that 
$(  \overline{\ker \phi^1_1}^{ \| \,\cdot\,\|_{m,\alpha -B} }, \ker \phi^1_1, \mu_1 )$
and 
$( (\ker \phi^1_1)^{ \bot} ,  (\ker \phi^1_1)^{ \bot}  , \mu_2)$ are abstract Wiener spaces.
Naturally, $\mu_1 = \pi_* \mu$, $\mu_2 = \pi_*^{\bot} \mu$
and $\mu = \mu_1 \times \mu_2$ (the product measure).
One may think $\mu_1$ is the definition of the conditional measure 
${\mathbb P}[   \,\cdot\, |~ \phi^1_1  =0 ]$ in {\bf (A3)'} above.

Therefore, {\bf (A3)'} is equivalent to the following;
\begin{equation}\label{equia3.ineq}
{\mathbb E}[  \exp (  \la\bar\nu ,  \phi^2_1 \circ \pi  \ra ) ]
< \infty.
\end{equation}
Set 
\begin{align}
\psi (w, w') 
&= 
\frac12 \tilde{J}(\bar\gamma)_1 \int_0^1 \tilde{J}(\bar\gamma)_t^{-1}  
\{
\nabla \sigma ( \phi^0_t)\la \phi^1_t (w'),  dw_t \ra 
+
\nabla \sigma ( \phi^0_t)\la \phi^1_t (w),  dw'_t \ra 
\}
\nn\\
&\quad 
+ \frac12
 \tilde{J}(\bar\gamma)_1 \int_0^1 \tilde{J}(\bar\gamma)_t^{-1}  
  \nabla^2 \sigma ( \phi^0_t)\la \phi^1_t (w),   \phi^1_t (w'), d\bar\gamma_t \ra,
\label{equia6.ineq}
\end{align}
where $\phi^1_T (w) = \tilde{J}(\bar\gamma)_T
 \int_0^T \tilde{J}(\bar\gamma)_t^{-1} \sigma ( \phi^0_t) dw_t$.
Then, $\psi$ is  a bounded bilinear mapping  on ${\cal W}_B$
and so is $\psi \la \pi  \cdot, \pi \cdot \ra$.
Clearly, $\psi (w,w) = \phi^2_1(w)$ and $\psi (\pi w, \pi w) = \phi^2_1(\pi w)$. 
By Goodman's theorem (see Theorem 4.6, p. 83, \cite{kuo}), 
restricted on ${\cal H} \times {\cal H}$,
$\la \bar\nu, \psi \la \pi  \cdot, \pi \cdot \ra \ra$ 
is of trace class and, in particular,   Hilbert-Schmidt.  
The corresponding trace class operator on ${\cal H}$
and corresponding element of the second Wiener chaos are 
  denoted by $A$ and  $\Xi_A$, respectively.
Then, $\la\bar\nu ,  \phi^2_1( \pi w) \ra = \Xi_A (w) +
 {\rm Tr} ( A)$.
Hence, (\ref{equia3.ineq}) is equivalent to 
${\mathbb E}[  \exp (  \Xi_A  ) ]< \infty$, which in turn is equivalent to 
$\sup {\rm Spec }(A) <1/2$.
Since the inequality is strict, 
there exists $r >1$ such that $\sup {\rm Spec }(rA) <1/2$.
This implies ${\mathbb E}[  \exp (  \Xi_{r A}  ) ]= {\mathbb E}[  \exp ( r \Xi_{ A}  ) ]< \infty$.
Summing it up, we have seen that  {\bf (A3)'} is
equivalent to the following;
\begin{equation}\label{equia7.ineq}
{\mathbb E}[  \exp ( r \la\bar\nu ,  \phi^2_1 \circ \pi  \ra ) ]
< \infty
\qquad \mbox{  for some $r >1$.}
\end{equation}
%
%


Let us check here that {\bf (A3)} and {\bf (A3)'}
are equivalent under {\bf (A1)}, {\bf (A2)}.
\begin{pr}\label{a3equi.pr}
Under {\bf (A1)} and {\bf (A2)}, the two conditions {\bf (A3)} and {\bf (A3)'}
are equivalent.
\end{pr}

\Proof
As is explained above, ${\bf (A3)'}$ is equivalent to $\sup {\rm Spec }(A) <1/2$.
Keep in mind that the only accumulation point of ${\rm Spec }(A)$ is $0$,
since $A$ is of trace class.
Let $(- \ve_0, \ve_0) \ni u \mapsto f(u) \in K_a^{a'}$ be a smooth curve in $K_a^{a'}$
such that $f(0) = \bar\gamma$ and $f^{\prime} (0)  \neq 0$ as in {\bf (A3)}.
Then, a straight forward calculation shows that
\begin{align}
\lefteqn{
\frac{d^2}{du^2} \Big|_{u=0} \frac{ \|   f(u) \|^2_{{\cal H}} }{2}
=
\frac{d^2}{du^2} \Big|_{u=0}
\Bigl(  \frac{ \|   f(u) \|^2_{{\cal H}} }{2}  - \la \bar\nu, \phi^0_1(f_u) -a' \ra  \Bigr)
}
\nn\\
&
=
\|   f^{\prime}(0) \|^2_{{\cal H}}+ \la  f^{\prime\prime}(0),  \bar\gamma \ra_{{\cal H}}
 - \bigl\la \bar\nu, D\phi^0_1(\bar\gamma) \la f^{\prime\prime}(0) \ra \bigr\ra 
 - \bigl\la \bar\nu, D^2\phi^0_1(\bar\gamma) \la f^{\prime}(0), f^{\prime}(0) \ra \bigr\ra 
    \nn\\
&
=
\|   f^{\prime}(0) \|^2_{{\cal H}}
-  \bigl\la \bar\nu, D^2\phi^0_1(\bar\gamma) \la \pi f^{\prime}(0), \pi f^{\prime}(0) \ra \bigr\ra
\nn\\
&=
\|   f^{\prime}(0) \|^2_{{\cal H}}
-  2 \bigl\la \bar\nu, \psi \la \pi f^{\prime}(0), \pi f^{\prime}(0) \ra \bigr\ra,
 \label{2kaibi.eq}
    \end{align}
where we used (\ref{Lmul1.eq})--(\ref{Lmul2.eq}) and the fact that $f^{\prime}(0)$
is tangent to the submanifold $K_a^{a'}$.
Since  $f^{\prime}(0)$ can be any non-zero element in ${\rm Im}~ \pi$,
$\sup {\rm Spec }(A) <1/2$ is equivalent to 
that right hand side of (\ref{2kaibi.eq})
is strictly positive, that is {\bf (A3)}.
\QED


The following is a key technical lemma. 
It states that, restricted on a sufficiently small subset,
 $\exp(\la \bar\nu , R_1^{2, \ve} \ra /\ve^2) \in \cup_{1< q <\infty} L^q$ uniformly in $\ve$.
\begin{lm}\label{lm.integ3}
Assume {\bf (A1)}, {\bf (A2)} and {\bf (A3)}.
Then, there exists $r_1 >1$ and $\eta >0$ such that
\[
\sup_{0 < \ve \le 1}
{\mathbb E} \bigl[
\exp \bigl(
    r_1  \la \bar\nu ,  R_1^{2, \ve}   \ra /\ve^2
           \bigr) 
I_{  \{  \| \ve w \|_{m, \alpha -B} \le \eta  \}}
I_{  \{   | R_1^{1, \ve}   /\ve | \le \eta_1  \}}
\bigr]  <\infty
\]
for any $\eta_1 >0$.
\end{lm}

\Proof
By Lemma \ref{lm.integ1} and the relation $R_1^{2, \ve}  /\ve^2 = \phi^2_1 +R_1^{2+, \ve}  /\ve^2 $,
it is sufficient to show that
\begin{equation}
\sup_{0 < \ve \le 1}
{\mathbb E} \bigl[
\exp \bigl(
    r_1  \la \bar\nu ,   \phi^2_1 \ra 
           \bigr) 
I_{  \{  \| \ve w \|_{m, \alpha -B} \le \eta  \}}
I_{  \{   | R_1^{1, \ve}   /\ve | \le \eta_1  \}}
\bigr]  <\infty.
\label{integ3suff.eq}
\end{equation}

We give an explicit expression for the projection $\pi$.
Set $C_{jj'} = \la \phi_1^{1,j} ,  \phi_1^{1,j'}\ra_{{\cal H}^*}$ 
and $C= ( C_{jj'})_{1 \le j, j' \le n} \in {\rm GL}(n, {\bf R})$.
The components of its inverse is denoted by $C^{-1}= ( D_{jj'})_{1 \le j, j' \le n}$.
By straight forward calculation, $\pi : {\cal H} \to \ker \phi_1^1$ 
is given by 
\[
\pi h = h - \sum_{j,j'}  {}_{{\cal H}^*}\la \phi_1^{1,j} , h \ra_{{\cal H}}
D_{j j'} 
\cdot  {}^{\sharp} \phi_1^{1,j'}.
\] 
From this, it is easy see that $\pi : {\cal W}_B \to \overline{\ker \phi_1^1}$ 
is given by 
\begin{equation}
\pi w = w - \sum_{j,j'}   \phi_1^{1,j} (w)
D_{j j'} 
\cdot  {}^{\sharp} \phi_1^{1,j'}.
\label{proj1.eq}
\end{equation}
Then,  we have
\begin{align}
 \phi^2_1 (w)  &=  \psi \la w,w \ra 
     =\phi^2_1 ( \pi w)  + 2    \sum_{j,j'}   \phi_1^{1,j} (w)
D_{j j'} 
\cdot   \psi  \la w,   {}^{\sharp} \phi_1^{1,j'} \ra
   \nn\\
      & \qquad  + 
          \sum_{j,j' , k, k'}  \phi_1^{1,j} (w) \phi_1^{1,k} (w)D_{j j'} D_{k k'} \cdot   \psi  \la   {}^{\sharp} \phi_1^{1,j'}  ,   {}^{\sharp} \phi_1^{1,k'}  \ra
=: J_1+ J_2+ J_3.
                                     \label{proj2.eq}
                                        \end{align}

Exponential integrability of the first term $J_1$ 
on the right hand side of (\ref{proj2.eq}) is given in (\ref{equia7.ineq}).
So, we estimate the second term $J_2$.
Since $\ve \phi_1^{1} (w) = R_1^{1+, \ve} (w)- R_1^{1, \ve} (w)$,
\begin{align}
|\phi_1^{1,j} (w) \psi  \la w,   {}^{\sharp} \phi_1^{1,j'} \ra|
&\le
c_1 \Bigl\{
\Bigl|\frac{ R_1^{1+, \ve} (w)  }{\ve}\Bigr| + \Bigl| \frac{ R_1^{1, \ve} (w) }{\ve}  \Bigr|   \Bigr\}
\|w\|_{m,\alpha -B}  
\nn\\
&\le
c_1 \Bigl\{  \Bigl|\frac{ c' R_1^{1+, \ve} (w)  }{\ve}\Bigr|^2 
+ \frac{\|w\|_{m,\alpha -B}^2 }{4 c'^2}   \Bigr\}
+
c_1 \Bigl| \frac{ R_1^{1, \ve} (w) }{\ve}  \Bigr|    \|w\|_{m,\alpha -B}  
\nn
\end{align}
for any $c' >0$.

Set $c_2 =2c_1n^2 \sup_{j,j'} |D_{j,j'}|$ 
and let $M>0$. 
Then, by H\"older's inequality,
\begin{align}
{\mathbb E} \bigl[
e^{M |J_2| }
I_{  \{  \| \ve w \|_{m, \alpha -B} \le \eta  \}}
I_{  \{   | R_1^{1, \ve}   /\ve | \le \eta_1  \}}
\bigr]  
\le
{\mathbb E} \bigl[
\exp \bigl( 3M c_2 c'^2  |R_1^{1+, \ve}   /\ve|^2 \bigr)
I_{  \{  \| \ve w \|_{m, \alpha -B} \le \eta  \}}
\bigr]^{1/3} 
\nn\\
 \quad  \times
{\mathbb E} \bigl[
e^{  3M c_2 \|w\|_{m,\alpha -B}^2 /(4 c')}
\bigr]^{1/3} 
{\mathbb E} \bigl[
e^{  3M c_2 \eta_1 \| w \|_{m, \alpha -B}  }
\bigr]^{1/3}.
\nn
\end{align}
For any $M>0$ and $\eta_1 >0$, the third factor is integrable.
If $c'$ is chosen sufficiently large, 
then the second factor is also integrable by Fernique's theorem.
By Lemma \ref{lm.integ2},  there exists $\eta >0$ such that $\sup_{\ve}$ of the first 
factor is finite and, hence, 
\begin{equation}\label{proj3.eq}
\sup_{0< \ve \le 1}
{\mathbb E} \bigl[
e^{M |J_2| }
I_{  \{  \| \ve w \|_{m, \alpha -B} \le \eta  \}}
I_{  \{   | R_1^{1, \ve}   /\ve | \le \eta_1  \}}
\bigr]  
< \infty.
\end{equation}

Since $\phi_1^{1,j} (w) \phi_1^{1,k} (w) 
= \ve^{-1} \{ R_1^{1+, \ve} (w)^j- R_1^{1, \ve} (w)^j \} \phi_1^{1,k} (w)$,
we can deal with $J_3$ in the same way.
For any $M>0$ and $\eta_1 >0$,  there exists $\eta >0$ such that
\begin{equation}\label{proj4.eq}
\sup_{0< \ve \le 1}
{\mathbb E} \bigl[
e^{M |J_3| }
I_{  \{  \| \ve w \|_{m, \alpha -B} \le \eta  \}}
I_{  \{   | R_1^{1, \ve}   /\ve | \le \eta_1  \}}
\bigr]  
< \infty.
\end{equation}

Let $r>1$ be as in (\ref{equia7.ineq}).
Set $r_1=(1+r)/2 >1$, $q =2r/(1+r) >1$, and $1/q +1/q' =1$.
Then, from H\"older's inequality and (\ref{equia7.ineq}), (\ref{proj2.eq})--(\ref{proj4.eq}),
we can easily see that
\begin{align}
\lefteqn{
{\mathbb E} \bigl[
\exp \bigl(
    r_1  \la \bar\nu ,   \phi^2_1 \ra 
           \bigr) 
I_{  \{  \| \ve w \|_{m, \alpha -B} \le \eta  \}}
I_{  \{   | R_1^{1, \ve}   /\ve | \le \eta_1  \}}
\bigr]  
}
\nn\\
&\le
{\mathbb E} \bigl[  \exp \bigl(
    r \la \bar\nu ,   \phi^2_1 \circ \pi \ra 
           \bigr) \bigr]^{1/q}
\prod_{i=1}^2
{\mathbb E} \bigl[
e^{ 2q' r_1 |\bar\nu| |J_i | }
I_{  \{  \| \ve w \|_{m, \alpha -B} \le \eta  \}}
I_{  \{   | R_1^{1, \ve}   /\ve | \le \eta_1  \}}
\bigr]^{1/(2q')}.
\nn
\end{align}
From this, (\ref{integ3suff.eq}) is immediate.
This completes the proof.
\QED


\subsection{Proof of off-diagonal short time asymptotics}

In this subsection we prove Theorem \ref{thm.MAIN.off}, namely,
off-diagonal short time asymptotics
of the density of the solution $(y_t) = (y_t (a))$ 
of Young ODE (\ref{ygSDE.eq}) 
driven by fBm $(w_t)$ with $1/2 <H<1$ under Assumptions {\bf (A1)}--{\bf (A3)}. 

%
%

First, let us calculate the kernel $p(t, a, a')$.
Take  $\eta >0$ as in Lemma \ref{lm.integ3}.
Then, we see
\begin{align}
p(\ve^{1/H}, a, a') 
&=
{\mathbb E} \bigl[  
\delta_{a'} ( y_1^{\ve})   
 \bigr]
 \nn\\
 &= 
   {\mathbb E} \bigl[  
\delta_{a'} ( y_1^{\ve} )   \chi_{\eta} (\ve, w)
 \bigr]
   + 
    {\mathbb E} \bigl[  
\delta_{a'} ( y_1^{\ve} )  \bigl\{ 1- \chi_{\eta} (\ve, w)  \bigr\}
 \bigr]
 = : I_1 +I_2.
 \nn
  \end{align}
As we have shown in Lemma \ref{lm.ldpcut}, the second term $I_2$ on the right hand side 
does not contribute to the asymptotic expansion.
So, we have only to calculate the first term $I_1$.
By Cameron-Martin formula, 
\[
I_1 
=
   {\mathbb E} \bigl[  
   \exp \bigl(  -\frac{ \|  \bar\gamma \|^2_{{\cal H}}}{2\ve^2}   - \frac{1}{\ve} \la \bar\gamma, w \ra \bigr)
\delta_{a'} ( \tilde{y}_1^{\ve} )   \chi_{\eta} (\ve, w + \frac{\bar\gamma}{\ve})
 \bigr].
 \]
Recall that $\la \bar\gamma, w \ra = \la \bar\nu , \phi^1_1 (w) \ra$ for all $w$.
Hence, noting that $\phi^{1/H}$ is non-random, we have
\begin{align}
I_1  
&=
 \exp \bigl(  -\frac{ \|  \bar\gamma \|^2_{{\cal H}}}{2\ve^2}  \bigr)
  {\mathbb E} \bigl[  
   \exp \bigl(
    - \frac{1}{\ve} \la \bar\nu , \phi^1_1  \ra
           \bigr)
\delta_{a'} ( a' +\ve \phi_1^1 + \ve^{1/H} \phi_1^{1/H}   +R_1^{2, \ve} )   \chi_{\eta} (\ve, w + \frac{\bar\gamma}{\ve})
 \bigr]
    \nn\\
      &=  \frac{1}{\ve^n}  \exp \bigl(  -\frac{ \|  \bar\gamma \|^2_{{\cal H}}}{2\ve^2}  \bigr)
  {\mathbb E} \bigl[  
   \exp \bigl(
    - \frac{1}{\ve} \la \bar\nu , \phi^1_1  \ra
           \bigr)
\delta_{0} (  \phi_1^1 + \ve^{(1/H)-1} \phi_1^{1/H}+ \ve^{-1} R_1^{2, \ve} )   \chi_{\eta} (\ve, w + \frac{\bar\gamma}{\ve})
 \bigr]
    \nn\\
      &=    
        \frac{1}{\ve^n}  \exp \Bigl(  -\frac{ \|  \bar\gamma \|^2_{{\cal H}}}{2\ve^2}  
     +  \frac { \la \bar\nu ,  \phi_1^{1/H} \ra  }  { \ve^{2 - (1/H)} }   \Bigr)
           \nn\\
           &  \qquad \qquad 
\times 
  {\mathbb E} \bigl[  
   \exp \bigl(
      \la \bar\nu ,  R_1^{2, \ve}   \ra /\ve^2
           \bigr)
\delta_{0} (  \phi_1^1 + \ve^{(1/H)-1} \phi_1^{1/H}+ \ve^{-1} R_1^{2, \ve} )   \chi_{\eta} (\ve, w + \frac{\bar\gamma}{\ve})
 \bigr]
    \nn\\
      &= 
     \frac{1}{\ve^n}  \exp \Bigl(  -\frac{ \|  \bar\gamma \|^2_{{\cal H}}}{2\ve^2}  
    + \frac { \la \bar\nu ,  \phi_1^{1/H} \ra  }  { \ve^{2 - (1/H)} }   \Bigr)
        {\mathbb E} \bigl[      F(\ve, w)
\delta_{0} \bigl(  \frac{ \tilde{y}_1^{\ve} - a'}{\ve } \bigr) 
 \bigr],
 \nn
       \end{align}
where 
\begin{align}
F(\ve, w) =  \exp \bigl(  \ve^{-2}
      \la \bar\nu ,  R_1^{2, \ve}   \ra 
           \bigr) 
           \chi_{\eta} (\ve, w + \frac{\bar\gamma}{\ve})       
\psi \Bigl(  \frac{1}{\eta_1^2}  \Bigl| \frac{\tilde{y}_1^{\ve} - a'}{\ve}   \Bigr|^2 \Bigr)
\label{Fve.def}
\end{align}
for any positive constant $\eta_1$. 
It is easy to see that {\rm (i)}~
$\chi_{\eta} (\ve, w +  \bar\gamma /\ve)  $ and its derivatives vanish 
outside $\{   \|\ve w\|_{m, \alpha -B} \le \eta \}$
and  {\rm (ii)}~
$\psi \bigl( \eta_1^{-2}  \bigl| (\tilde{y}_1^{\ve} - a' )/\ve   \bigr|^2  \bigr)$ 
and its derivatives vanish 
outside $\{  |R_1^{1, \ve} /\ve | \le \eta_1 \}$.
Hence, by Lemma \ref{lm.integ3}, $F(\ve, w) \in \tilde{\bf D}_{\infty}$ and 
$F(\ve, w) =O(1)$ with respect to that topology.
Roughly speaking,
since $\delta_{0} (  ( \tilde{y}_1^{\ve} - a')/\ve  )$ admits an asymptotic expansion
in $\tilde{\bf D}_{-\infty}$,
the problem reduces to whether 
$F(\ve, w) $  admits an asymptotic expansion
in $\tilde{\bf D}_{\infty}$.


\begin{lm}\label{lm.sikihen1}
Assume {\bf (A1)}--{\bf (A3)}.
For any $M \in {\bf N}$,  we have
\[
  {\mathbb E} \bigl[      F(\ve, w)
\delta_{0} \bigl(  \frac{ \tilde{y}_1^{\ve} - a'}{\ve}  \bigr) 
 \bigr]
=
  {\mathbb E} \bigl[      F(\ve, w)  \psi ( | \phi_1^1 / \eta_1|^2)
\delta_{0} \bigl(  \frac{ \tilde{y}_1^{\ve} - a'}{\ve}  \bigr) 
 \bigr]
+
O(\ve^M)
\]
as $\ve \searrow 0$.
\end{lm}

\Proof
By using Taylor expansion for $\psi$, we see that, for given $M$,
there exist $m \in {\bf N}$ and $G_j (\ve, w) \in {\bf D}_{\infty} ~(1 \le j \le m)$
such that
\begin{align}
\psi \Bigl(  \frac{1}{\eta_1^2}  \Bigl| \frac{\tilde{y}_1^{\ve} - a'}{\ve}   \Bigr|^2 \Bigr)
&=
\psi \Bigl( \bigl| \frac{\phi_1^1}{\eta_1}  \bigr|^2 \Bigr)
+ \sum_{j=1}^m \psi^{(j)} \Bigl( \bigl| \frac{\phi_1^1}{\eta_1} \bigr|^2   \Bigr) G_j (\ve, w)
+O(\ve^M)
\label{siki2.eq}
\end{align}
in ${\bf D}_{\infty}$ as $\ve \searrow 0$.
$G_j (\ve, w) =O(1)$,  but its explicit form is not important.
Note that $\psi^{(j)} (| \phi_1^1/\eta_1|^2 ) T(\phi_1^1)=0$ 
if $j \ge 1$ and ${\rm supp} (T) \subset \{ a \in {\bf R}^n~|~|a| <  \eta_1 /2\}$.

By Theorem \ref{tm.asym.plbk} and Proposition \ref{pr.unif.nondeg},
$\delta_{0} (  ( \tilde{y}_1^{\ve} - a')/\ve )$ admits an asymptotic expansion
in $\tilde{\bf D}_{-\infty}$ as follows.
As before, we set $\{ 0= \nu_0 < \nu_1 < \nu_2 <\cdots\}$ to be all the elements of $\Lambda_3$
in increasing order. 
For given $M$, let $l \in {\bf N}$ be the smallest integer such that $M \le \nu_{l+1}$.
Then, for some 
$\Phi_{\nu_j}  \in \tilde{\bf D}_{-\infty} ~(1 \le j \le l)$,
it holds that
\begin{align}
\delta_{0} (  ( \tilde{y}_1^{\ve} - a')/\ve )
=
\delta_{0} ( \phi_1^1 )
+
\ve^{\nu_1} \Phi_{\nu_1} +\cdots + \ve^{\nu_l} \Phi_{\nu_l}
+O(\ve^{ \nu_{l+1}})
\label{siki3.eq}
\end{align}
in $\tilde{\bf D}_{-\infty}$ as $\ve \searrow 0$.
Here, $\Phi_{\nu_j}$ is a finite linear combination of terms of the form
\[
\partial^{\alpha} \delta_0 (\phi_1^1) \times 
\{
\mbox{a polynomial of the components of $\phi_1^{\kappa_i}$\,'s}
\}.
\]
Hence, $\psi^{(j')} (| \phi_1^1/\eta_1|^2 )  \Phi_{\nu_j}$ vanish for all $j, j'$.

Now, using (\ref{siki2.eq}) and (\ref{siki3.eq}), we prove the lemma. 
\begin{align}
\lefteqn{ 
 {\mathbb E} \bigl[      F(\ve, w)
\delta_{0} (  ( \tilde{y}_1^{\ve} - a')/\ve ) 
 \bigr]
}
\nn\\
&=
 {\mathbb E} \bigl[      F(\ve, w)
\psi \Bigl(  \frac{1}{\eta_1^2}  \Bigl| \frac{\tilde{y}_1^{\ve} - a'}{\ve}   \Bigr|^2 \Bigr)
\delta_{0} (  ( \tilde{y}_1^{\ve} - a')/\ve ) 
 \bigr]
\nn\\
&=
 {\mathbb E} \bigl[      F(\ve, w)
\psi ( | \phi^1_1/\eta_1|^2 )
\delta_{0} (  ( \tilde{y}_1^{\ve} - a')/\ve ) 
 \bigr]
\nn\\
&  \quad
+  {\mathbb E} \bigl[      F(\ve, w)
\Bigl(  
\sum_{j=1}^m \psi^{(j)} \Bigl( \bigl| \frac{\phi_1^1}{\eta_1} \bigr|^2   \Bigr) G_j (\ve, w)
 \Bigr)
\delta_{0} (  ( \tilde{y}_1^{\ve} - a')/\ve ) 
 \bigr]
+
O(\ve^M)
\nn\\
&=
 {\mathbb E} \bigl[      F(\ve, w)
\psi ( | \phi^1_1/\eta_1|^2 )
\delta_{0} (  ( \tilde{y}_1^{\ve} - a')/\ve ) 
 \bigr]
\nn\\
&  \quad
+  {\mathbb E} \bigl[      F(\ve, w)
\Bigl(  
\sum_{j=1}^m \psi^{(j)} \Bigl( \bigl| \frac{\phi_1^1}{\eta_1} \bigr|^2   \Bigr) G_j (\ve, w)
 \Bigr)
\bigl(  \delta_{0} ( \phi_1^1 )
+
\cdots + \ve^{\nu_l} \Phi_{\nu_l} \bigr)
 \bigr]
+
O(\ve^M)
\nn\\
&=
  {\mathbb E} \bigl[      F(\ve, w)  \psi ( | \phi_1^1 / \eta_1|^2)
\delta_{0} (  ( \tilde{y}_1^{\ve} - a')/\ve ) 
 \bigr]
+
O(\ve^M).
\nn
\end{align}
Thus, we have shown the lemma.
\QED


Set 
$
\Lambda'_2 = \{ \kappa -2 ~|~ \kappa \in \Lambda_1 \setminus \{0, 1, 1/H\} \}
=
\{ 0 <  H^{-1} -1 <
   \bigl(  3 \wedge 2H^{-1} \bigr)-2  <\cdots \}.
$
Next we set 
$
\Lambda'_3 = 
\{  a_1 +a_2 + \cdots + a_m ~|~  \mbox{$m \in {\bf N}_+$ and $a_1 ,\ldots, a_m \in \Lambda'_2$} \}$.
In the following lemma, 
$\{ 0=\rho_0 <\rho_1<\rho_2 <\cdots \}$ stands for all the elements of 
$\Lambda'_3$ in increasing order.

\begin{lm}\label{lm.asyFve}
Assume {\bf (A1)}--{\bf (A3)} and let $F(\ve, w) \in \tilde{\bf D}_{\infty}$ as in (\ref{Fve.def}).
Then, for every $k=1,2,3, \ldots$,
\begin{align}
\lefteqn{
F(\ve, w) \psi (  |  \phi^1_1 (w)  / \eta_1 |^2 )
}
\nn\\
 &= \exp \bigl(
      \la \bar\nu ,  \phi^2_1 (w)
           \bigr) 
\psi (  |  \phi^1_1 (w)  / \eta_1 |^2 )^2
 \{1+ \ve^{\rho_1} \gamma_{\rho_1}(w)  +\cdots + \ve^{\rho_k} \gamma_{\rho_k} (w)\}
+
F_{k+1} (\ve, w),
\nn
\end{align}
where $F_{k+1} (\ve, w) \in \tilde{\bf D}_{\infty}$ satisfies that
\[
F_{k+1} (\ve, w) T ( \phi^1_1 ) = O(\ve^{\rho_{k+1}})  \qquad
\mbox{ in ${\bf D}_{-\infty}$ as $\ve \searrow 0$}
\]
for any $T \in {\cal S}^{\prime} ({\bf R}^n)$ with ${\rm supp}(T) 
\subset \{ a\in {\bf R}^n ~|~ |a| \le \eta_1/2 \}$.
Moreover, $\gamma_{\rho_j} \in {\bf D}_{\infty} ~(j=1,2,\ldots)$ are determined by the following 
formal expansion ($\kappa_4 =H^{-1} +1$);
\begin{align}
\sum_{m=0}^{\infty}  \frac{
\la  \bar\nu,  R_1^{2+ , \ve}  /\ve^2\ra^m }{m!}
&=
\sum_{m=0}^{\infty}
\frac{1}{m!}
\Bigl\{
  \ve^{\kappa_4 -2} \la  \bar\nu,  \phi_1^{\kappa_4} \ra
+  \ve^{\kappa_5 -2} \la  \bar\nu, \phi_1^{\kappa_5}  \ra +\cdots   
\Bigr\}^m
\nn\\
&=
 1+ \ve^{\rho_1} \gamma_{\rho_1} + \ve^{\rho_2} \gamma_{\rho_2}  +\cdots.
\nn
\end{align}
\end{lm}

\Proof
Let $r_1 >1$ be as in Lemma \ref{lm.integ3}.
First we show that, for any $\eta_1 >0$, 
\begin{equation}
{\mathbb E} \bigl[
\exp \bigl(
    r_1  \la \bar\nu ,   \phi^2_1 \ra 
           \bigr) 
I_{  \{   | \phi_1^1 | \le \eta_1  \}}
\bigr]  <\infty.
\label{siki1.eq}
\end{equation}
We can choose a subsequence $\{ \ve_k \}$
such that, as $k \to \infty$, $\ve_k \searrow 0$ and 
$R_1^{1, \ve_k} / \ve_k  \to \phi^1_1$ a.s.
To prove (\ref{siki1.eq}), we apply Fatou's lemma
to  (\ref{integ3suff.eq}) with $\eta_1$ replaced by $2\eta_1$.
\begin{align}
\infty &> \liminf_{k \to \infty}
{\mathbb E} \bigl[
\exp \bigl(
    r_1  \la \bar\nu ,   \phi^2_1 \ra 
           \bigr) 
I_{  \{  \| \ve_k w \|_{m, \alpha -B} \le \eta  \}}
I_{  \{   | R_1^{1, \ve_k}   /\ve_k | \le 2\eta_1  \}}
\bigr] 
\nn\\
&\ge 
{\mathbb E} \bigl[
\exp \bigl(
    r_1  \la \bar\nu ,   \phi^2_1 \ra 
           \bigr) 
\liminf_{k \to \infty}
I_{  \{   | R_1^{1, \ve_k}   /\ve_k | \le 2\eta_1  \}}
\bigr] 
\ge 
{\mathbb E} \bigl[
\exp \bigl(
    r_1  \la \bar\nu ,   \phi^2_1 \ra 
           \bigr) 
I_{  \{   |\phi^1_1 | \le \eta_1  \}}
\bigr]. 
\nn
\end{align}
From (\ref{siki1.eq}), it is easy to check that 
$\exp \bigl(
      \la \bar\nu ,  \phi^2_1 (w) \ra
           \bigr) 
\psi (  |  \phi^1_1 (w)  / \eta_1 |^2 ) \in \tilde{\bf D}_{\infty}$.

Now we expand $\exp ( \la \bar\nu ,  R_1^{2, \ve}   \ra /\ve^2 ) 
= \exp (    \la \bar\nu ,  \phi^2_1 (w) \ra ) \exp ( \la \bar\nu ,  R_1^{2+ , \ve}   \ra /\ve^2 ) $ in $\ve$.
Set $Q_{l+1} : {\bf R} \to {\bf R}$ by
\[
Q_{l+1} (u) = e^u -\Bigl(  1+ u +\frac{u^2}{2!}  +\cdots +  \frac{u^l}{l!}  \Bigr) = 
u^{l+1} 
\int_0^1 \frac{(1-\theta)^l }{l!} e^{\theta u}  d\theta 
\qquad
(u \in {\bf R}).
\]
We will prove that, for sufficiently large $l \in {\bf N}$, as $\ve \searrow 0$,
\begin{equation}\label{siki4.eq}
e^{ \la \bar\nu ,  \phi^2_1  \ra }
Q_{l+1}( \la  \bar\nu, R_1^{2+ , \ve} \ra / \ve^2 )   \chi_{\eta} (\ve, w + \frac{\bar\gamma}{\ve})  \psi (  |  \phi^1_1 (w)  / \eta_1 |^2 ) 
=
O( \ve^{\rho_{k+1}})
\quad
\mbox{in $\tilde{\bf D}_{\infty}$. }
\end{equation}
Note that $  \chi_{\eta} (\ve, w + \frac{\bar\gamma}{\ve})  =O(1)$ in ${\bf D}_{\infty}$ as $\ve \searrow 0$ by (\ref{chi_asy.eq}).
By Proposition \ref{pr.tay.Dinf},
$R_1^{2+ , \ve} / \ve^2 =O(\ve^{(1/H) -1})$ in ${\bf D}_{\infty}$.
So, if $l +1 \ge \rho_{k+1}/  \{(1/H) -1\}$, 
then 
$( \la \bar\nu, R_1^{2+ , \ve} \ra / \ve^2 )^{l+1} =O(\ve^{ \rho_{k+1}})$ in ${\bf D}_{\infty}$.
Therefore, in order to verify (\ref{siki4.eq}), it is sufficient to show that, as $\ve \searrow 0$,
\begin{equation}\label{siki5.eq}
\int_0^1 (1-\theta)^l  
e^{ \la     \bar\nu ,  \phi^2_1 +   \theta   R_1^{2+ , \ve} / \ve^2 \ra }  d\theta \cdot
 \chi_{\eta} (\ve, w + \frac{\bar\gamma}{\ve})  \psi (  |  \phi^1_1 (w)  / \eta_1 |^2 ) 
=
O( 1)
\quad
\mbox{in $\tilde{\bf D}_{\infty}$. }
\end{equation}
To verify the integrability of this Wiener functional, 
note that $e^{\theta u}  \le 1+ e^u$ for all $u \in {\bf R}$
and $0 \le \theta \le 1$.
This implies that the first factor on the left hand side of (\ref{siki5.eq})
 is dominated by $e^{ \la \bar\nu ,  \phi^2_1  \ra } + e^{\la \bar\nu, R_1^{2 , \ve} \ra / \ve^2}$.
From Lemma \ref{lm.integ3}
and (\ref{siki1.eq}),  we see that the left hand side of  (\ref{siki5.eq}) is 
 $O(1)$ in any $L^q ~(1<q<\infty)$.
In the same way, 
the Malliavin derivatives of the left hand side of (\ref{siki5.eq}) are $O(1)$ in any $L^q$.

It is easy to see that, as $\ve \searrow 0$,
\begin{align}
\sum_{k=0}^l
\frac{  \{ \la  \bar\nu, R_1^{2+ , \ve} \ra  / \ve^2 \}^k  }{k!}  
=
 1+ \ve^{\rho_1} \gamma_{\rho_1} +\cdots + \ve^{\rho_k} \gamma_{\rho_k}  
 +
O( \ve^{\rho_{k+1}})
\quad
\mbox{in ${\bf D}_{\infty}$. }
\end{align}
From this and (\ref{chi_asy.eq}), we see that 
\begin{align}
\lefteqn{
F(\ve, w) \psi (  |  \phi^1_1 (w)  / \eta_1 |^2 )
}
\nn\\
 &= \exp \bigl(
      \la \bar\nu ,  \phi^2_1 (w)
           \bigr) 
\psi (  |  \phi^1_1 (w)  / \eta_1 |^2 )  
 \psi \Bigl(  \frac{1}{\eta_1^2}  \Bigl| \frac{\tilde{y}_1^{\ve} - a'}{\ve}   \Bigr|^2 \Bigr) \{1+ \ve^{\rho_1} \gamma_{\rho_1}(w)  +\cdots + \ve^{\rho_k} \gamma_{\rho_k} (w)\}
\nn\\
& \qquad\qquad 
 +
O( \ve^{\rho_{k+1}})
\quad
\mbox{in $\tilde{\bf D}_{\infty}$. }
\nn
\end{align}
Using (\ref{siki2.eq}), we finish the proof.
\QED


{\it Proof of Theorem \ref{thm.MAIN.off}}~
Here we prove our main theorem in this paper.
We set 
$$\Lambda_4 =\Lambda_3  +\Lambda'_3 =\{ \nu + \rho ~|~ \nu   \in  \Lambda_3  , \rho \in  \Lambda'_3\}.$$
We denote by $\{ 0 = \lambda_0 < \lambda_1 < \lambda_2 < \cdots \}$ all the elements of $\Lambda_4$  in increasing order.
There is no mystery why this index set  appears in the short time expansion of the kernel 
because, very formally speaking,  the problem reduces to finding asymptotic behavior  of
$ {\mathbb E}[  \exp(  \la \bar\nu, R_1^{2 , \ve} \ra  / \ve^2 )  \cdot \delta_0 (R_1^{1, \ve} / \ve) ]$, as we have seen.
Now, by (\ref{Fve.def}), Lemma \ref{lm.sikihen1},
Lemma \ref{lm.asyFve}, and (\ref{siki3.eq}), we can easily prove Theorem \ref{thm.MAIN.off}.
(First, expand the Watanabe distribution by (\ref{siki3.eq}),
then expand $F$ by Lemma \ref{lm.asyFve}.)
\toy

\section{Sufficient condition for {\bf (A2)} and {\bf (A3)}}
\label{sec.cmp}

In this final section we give a sufficient condition for
our main result (Theorem \ref{thm.MAIN.off}) 
on the off-diagonal asymptotics 
and compare it 
with a preceding result by Baudoin and Ouyang (Theorem 1.2, \cite{bo}),
which is probably the only paper on this kind of problem.

\begin{pr}\label{pr.boin_new}
Assume {\bf (A1)} at the starting point $a \in {\bf R}^n$.
If $a'$ is sufficiently near $a$,
then {\bf (A2)} and  {\bf (A3)}  are satisfied and, 
in particular, Theorem \ref{thm.MAIN.off} holds for such $a'$.
\end{pr}

In the latter half of this section,
we will prove this proposition in a rather general setting so that 
the same argument applies to a wider class of Gaussian processes.
(To obtain Proposition \ref{pr.boin_new},
just set $F=\phi^0_1$ and $x=a'$ in Proposition \ref{pr.impl}.)

Before doing so, we first recall the result in \cite{bo} and compare.
They set $n=d$ and assume {\bf (A1)} for {\it any} starting point $a \in {\bf R}^d$
and, moreover, the following assumption {\bf (H)}:
\\
\\
{\bf (H):}~ 
There exist smooth and bounded real-valued functions $\omega_{ij}^l$ such that 
\[
\omega_{ij}^l = - \omega_{il}^j  \qquad \mbox{and} \qquad
[V_i, V_j] = \sum_{i=1}^d \omega_{ij}^l V_l
 \qquad \mbox{for all $1 \le i,j,l \le d$.}
\]
\\
Note that $V_0$ does not appear in this condition.
Under {\bf (A1)} for any $a$, 
$\sigma (a)\sigma (a)^*$ is a $d \times d$ positive symmetric matrix,
where $\sigma (a) = [ V_1(a), \ldots, V_d (a)]$ as before.
As a result, a Riemannian metric tensor $(g_{ij}(a))_{1 \le i,j \le d}$ is defined on ${\bf R}^d$ by 
$g^{ij}(a) = [\sigma (a)\sigma (a)^*]^{ij}$.
The distance with respect to this Riemannian structure is denoted by $d(a,a')$.
In terms of Riemannian geometry, {\bf (H)} is equivalent to the condition that
$\nabla^{LC}_X Y =[X,Y]$ for all smooth vector fields $X, Y$, 
where $\nabla^{LC}$ is the Levi-Civita connection for this metric.
From this, one can guess that this assumption may not be very mild.

They proved short time kernel asymptotics 
under these assumptions when $a$ and $a'$ are sufficiently near.
The following is Theorem 1.2, \cite{bo} (Notations are adjusted):
\begin{tm}\label{tm.bomain}
Assume that $n=d$, $V_0 \equiv 0$, {\bf (H)},
and  {\bf (A1)} for any starting point $a \in {\bf R}^d$.
Then, in a neighborhood $U$ of $a$, we have
\begin{eqnarray*}
p(t, a, a') 
&=&
\frac{1}{ t^{Hn}}   \exp(- \frac{ d(a,a')^2 }{2t^{2H}}) 
\\
&&
\times \Bigl( \sum_{i=0}^N \alpha_{2i} (a,a') t^{2iH} + r_{N+1} (t,a,a') t^{2(N+1)H}
\Bigr),
\qquad a' \in U
\end{eqnarray*}
near $t=0$ for any $N =1,2,\ldots$.
Moreover, $U$ can be chosen so that $\alpha_{2i}$ are smooth on $U \times U$ and 
for all multi-indices $\beta, \beta'$ 
\[
\sup_{t \le t_0} \sup_{a,a' \in U \times U} 
 | \partial_a^{\beta}   \partial_{a'}^{\beta'}  r_{N+1}(t,a,a')| <\infty,
\qquad
(\mbox{for some $t_0 >0$}).
\]
\end{tm}

Now we compare the two results.
The most important issue is of course whether the asymptotic expansion holds or not.
Concerning this point, we observe  {\bf (i)-(ii)} below;
\\
\\
\noindent
{\bf (i)}~The conditions on the dimension ($n=d$), and 
on vector fields ($V_0 \equiv 0$ and {\bf (H)}) in \cite{bo} are 
much stronger than ours.
Moreover, the ellipticity condition {\bf (A1)} is assumed at any $a$ in \cite{bo}.
So we believe that our result is "basically" better than Theorem 1.2, \cite{bo}.
\\
\\
\noindent
{\bf (ii)}~
In our paper we did not give a quantitative estimate of 
how near $a$ and $a'$ should be in order for the asymptotics to hold 
(neither in \cite{bo}).
Therefore, we could not say our result completely includes Theorem 1.2, \cite{bo}.
\\
\\
The following {\bf (iii)} may not be a  major issue, 
but Theorem 1.2, \cite{bo} is better than ours concerning this point.
\\
\\
\noindent
{\bf (iii)}~
In Theorem \ref{tm.bomain}, or Theorem 1.2, \cite{bo}, 
they proved smoothness of the coefficient and 
gave an uniform estimate of (derivatives of)  the remainder terms.
However, we did not.


\begin{re}
If we assume {\bf (A1)} everywhere, then a Riemannian structure on ${\bf R}^n$ is 
naturally induced as we explained above. 
If the case of the usual stochastic analysis (i.e., $H=1/2$),
{\bf (A2)} and {\bf (A3)} have a geometric meaning. (See Remark 3.2, \cite{wa},
which was originally in \cite{mol, bi}.) 
First, {\bf (A2)} means that there is a unique shortest geodesics between $a$ and $a'$.
Second, {\bf (A3)} or {\bf (A3)'} means that these two points are not conjugate along the geodesics.
So, Assumptions {\bf (A1)--(A3)} are very mild and cover a lot of examples.

It seems natural to guess from this that, in our case (i.e., $1/2 <H<1$), too, 
 Assumptions {\bf (A1)--(A3)} are not bad.
 At this moment, however, the author is not aware of 
 a nice example except Proposition \ref{pr.boin_new}.
 \end{re}

For the rest of this section,  we discuss in a general setting.
Our goal here is to prove a generalized version of
Proposition \ref{pr.boin_new}.
The key is the implicit function theorem.

Let ${\cal H}$ be a real separable Hilbert space
and let $F: {\cal H} \to {\bf R}^n$ be a Fr\'echet smooth map
such that $F(0)=a$ and the tangent map $DF(h): {\cal H} \to {\bf R}^n$ 
is surjective at any $h \in {\cal H}$.
Necessarily, 
$F$ is a surjection onto a certain neighborhood of $a$ in ${\bf R}^n$.
By a well-known application of the inverse/implicit function theorem, 
$F^{-1}(x) \subset {\cal H}$ is a Hilbert submanifold for any $x \in {\bf R}^n$
if it is not empty.
We define
\begin{equation}
d(a, x) = \inf\{ \|h\|_{{\cal H}}  ~|~ h \in F^{-1}(x) \}.
\nn
\end{equation}
If $x$ is sufficiently near $a$, then $d(a, x) <\infty$.

\begin{pr}\label{pr.impl}
Let the notations be as above.
Furthermore, we assume that, for any $x$ sufficiently near $a$,
 the minimum in the definition of $d(a, x)$ above is actually attained.
Then, for any $x$ sufficiently near $a$,
we have the following;
\\
{\rm (i)}~
There exists a unique $h_x \in F^{-1}(x)$ such that $d(a, x) =  \|h_x\|_{{\cal H}}$.
\\
{\rm (ii)}~
The mapping $x \mapsto d(a, x)^2$ is smooth.
\\
{\rm (iii)}~The Hessian of $F^{-1}(x) \ni h \mapsto \|h\|_{{\cal H}}^2 /2$ 
at $h_x$ is non-degenerate in the sense in {\bf (A3)}.
\end{pr}

\Proof
Set ${\cal K}=\ker DF(0)$.
This is a closed linear subspace in ${\cal H}$
which is tangent to $F^{-1}(a)$ at $0$.
We denote by $\hat{D}$ and $\hat{D}^{\bot}$
the gradient operator on ${\cal K}$ and ${\cal K}^{\bot}$, respectively.
Then, $D = \hat{D} + \hat{D}^{\bot}$.
We often write $h =(k,l)$, where $k$ and $l$ are the orthogonal projections 
onto ${\cal K}$ and ${\cal K}^{\bot}$, respectively.

Consider the following function $G:{\cal K} \times {\cal K}^{\bot} \times {\bf R}^n 
(={\cal H}\times {\bf R}^n)
\to {\bf R}^n$ defined by 
$
G(k,l; x) = F(k, l) -x$.
Then, 
$G(0,0;a)=0$.
By the assumption, $(\hat{D}^{\bot} G)(0,0; a) = (\hat{D}^{\bot} F)(0,0)$ is a
linear isomorphism from ${\cal K}^{\bot}$ to ${\bf R}^n $.

Hence, we can use the implicit function theorem near $(0,0;a)$
to have the following;
There exist open neighborhoods $V \subset {\bf R}^n$ of $a$,
$W \subset {\cal K}$ of $0 \in {\cal K}$, and
$U \subset {\cal K}^{\bot}$ of $0 \in {\cal K}^{\bot}$ such that
a unique implicit function 
$l=l(k;x)$ for $G=0$
from $W \times V$ to $U$ exists. Moreover, $l$ is smooth.
Therefore, if $F^{-1}(x) \cap ( W \times U) \neq \emptyset$,
any element of the set is of the form $(k, l(k;x))$ for some $k \in W$.
Note that 
$l(0;a) =0$ and $\hat{D} l(0;a) =0 \in L({\cal K}, {\cal K}^{\bot})$ since $F^{-1}(a)$ and ${\cal K}$
are tangent at $0 \in {\cal H}$.

Next, consider $(k, x) \to \|(k, l(k;x))\|^2_{{\cal H}} /2 =(\|k\|^2 + \|l(k;x)\|^2)/2$.
Take $\hat{D}$ of this function and we get 
\[
\hat{G}(k, x) :=
\la k, \,\cdot\, \ra_{{\cal K}} + \la l(k;x), \hat{D} l(k;x)\ra_{{\cal K}^{\bot}},
\]
which is a smooth map from $W \times V$ to ${\cal K}^*$.
Note that $\hat{G}(0, a) =0$ and
\[
\hat{D}\hat{G}(k, x) =
\la \,\cdot\, , \,\cdot\, \ra_{{\cal K}} 
+ \la l(k;x), \hat{D}^2 l(k;x)\ra_{{\cal K}^{\bot}}
+ \la \hat{D} l(k;x), \hat{D} l(k;x)\ra_{{\cal K}^{\bot}}.
\]
This takes values in $L({\cal K}, {\cal K}^*) =L^{(2)}({\cal K}\times {\cal K};{\bf R})$,
where the latter space is the space of bounded bilinear maps from ${\cal K}\times {\cal K}$ to ${\bf R})$.
Since $\hat{D}\hat{G}(0,a) = \la \,\cdot\, , \,\cdot\, \ra_{{\cal K}}$,
which is clearly a linear isomorphism when regarded as an element of $L({\cal K}, {\cal K}^*)$,
we can use the implicit function theorem again.
If we retake $V$ and $W$ smaller,
then there exists a unique  implicit function $k =k(x)$ for $\hat{G}=0$
from $V$ to $W$.
Moreover, $k$ is smooth in $x$.

Take $r >0$ small enough so that
the open ${\cal H}$-ball $B_r$ of radius $r$ centered at $0 \in {\cal H}$
is contained in $W \times U$.
Assume $F^{-1}(x) \cap B_r \neq \emptyset$.
Then, the minimum is the definition of $d(x,a)$ must be achieved inside $B_r$.
That point can be written as $(k_0, l(k_0,x))$ in a unique way.
Any point of $F^{-1}(x)$ near  
$(k_0, l(k_0,x))$ can also be expressed using the implicit function like this.
As a result, 
this point must be a critical point of $k \mapsto \|(k, l(k;x))\|^2_{{\cal H}} /2$
and hence 
$\hat{G}(k_0 ,x)=0$.
Therefore, such $k_0$ must be unique, namely, $k_0 =k(x)$.
Note that $k(a)=0$.
Thus, we have seen $h_x =( k(x), l(k(x), x ))$
and shown {\rm (i)} and {\rm (ii)}.

We now show {\rm (iii)}.
Let $f: (-\ve_0, \ve_0) \to F^{-1}(x)$ such that $f(0) = k(x)$
and $f^{\prime}(0) \neq 0$.
Then, 
$(d/du)^2 \vert_{u=0} \|   f(u) \|^2_{{\cal H}} /2$
depends only on $f^{\prime}(0)$, i.e., $f^{\prime\prime}(0)$ is irrelevant.
(We can check this 
by using the Lagrange multiplier method in the same way 
as in (\ref{2kaibi.eq}) in the proof of Proposition \ref{a3equi.pr}.)
So, we have only to consider $f(u) = (k(x) + u\xi; l(k(x) + u\xi ;x))$
for any non-zero $\xi \in {\cal K}$.
By straight forward computation, we have
\[
\Bigl(\frac{d}{du} \Bigr)^2     \Big|_{u=0}    \frac{\|   f(u) \|^2_{{\cal H}} }{2}
=
\|\xi\|^2 + \bigl( l(k(x); x) , ( \hat{D}_{\xi})^2 l(k(x);x) \bigr) 
+
\| \hat{D}_{\xi} l(k(x);x)\|^2.
\]
By the smoothness of $l$ and $k$, the right hand side is larger than $\|\xi\|^2/2$
if $x$ is sufficiently near $a$. This proves {\rm (iii)}.
\QED


\noindent
{\bf Acknowledgement}~The author thanks an anonymous referee 
for suggesting how to improve Section \ref{sec.cmp}.


\end{document}